\documentclass[12pt,a4paper,twoside,final,notitlepage,leqno,english]{article}
\usepackage[babel=true]{csquotes}
\usepackage[T1]{fontenc}
\usepackage[latin1]{inputenc}
\usepackage{epsfig, graphicx, amssymb}
\usepackage{amsmath,amsthm,epsfig,amsfonts}
\usepackage{float}
\usepackage{color}
\usepackage[mathscr]{eucal}

\def\br {\break}

\newcommand{\moneq}{\vspace*{-7pt} \begin{equation} \displaystyle }
\newcommand{\moneqstar}{\vspace*{-6pt} \begin{equation*} \displaystyle }
\newcommand{\monendstar}{\vspace*{-6pt} \end{equation*}   }
\newcommand{\monend}{\vspace*{-7pt} \end{equation}   }
\newcommand{\moneqarraystar}{ \begin{eqnarray*} \displaystyle }
\newcommand{\monendarraystar}{ \end{eqnarray*}   }

\newcommand{\N}{\mathbb{N}}

\newcommand{\ZZ}[0]{\mathbb{Z}}

\def\ib#1{_{_{\scriptstyle{\!#1}}}}
\def \page #1{\unskip\leaders\hbox to 1.3 mm {\hss.\hss}\hfill   {$\,\,${#1}}}

\definecolor{vertfonce}{rgb}{0.0, 0.5, 0.0}

\def\section*#1{}

\usepackage{fancyhdr}
\renewcommand{\headrulewidth}{0pt}

\begin{document}

\fancypagestyle{plain}{ \fancyfoot{} \renewcommand{\footrulewidth}{0pt}}
\fancypagestyle{plain}{ \fancyhead{} \renewcommand{\headrulewidth}{0pt}}

~

  \vskip 2.1 cm


\centerline {\bf \LARGE  Pourquoi la gamme a sept notes ?}

\bigskip

\centerline {\bf \LARGE  Une math\'ematique des notes et des gammes}

\bigskip  \bigskip \bigskip

 \centerline {\large    Fran\c{c}ois Dubois}

\smallskip
\centerline {math\'ematicien appliqu\'e, amateur de   musique}
\centerline {Conservatoire National des Arts et M\'etiers, Paris}


\bigskip  \bigskip

\centerline {01 juin 2023 {\footnote{\rm $\,$ Ce texte a \'et\'e publi\'e dans le journal {\it Matapli}, \'edit\'e par la Soci\'et\'e
  de Math\'ematiques Appliqu\'ees et Industrielles, num\'ero 131, pages 51-89, en juin 2023 ; \'edition du 4 octobre 2023.}}}




\bigskip  \bigskip
\noindent {\bf \large R\'esum\'e}

\smallskip \noindent
Nous proposons une construction alg\'ebrique des notes de musique et montrons
comment les associer de diverses fa\c{c}ons  pour former des gammes.
Ainsi \'emerge toute une famille de gammes avec un nombre fix\'e de notes :
deux, trois, cinq, sept, douze, dix-sept, {\it etc.}
La classification des   gammes est simple avec le concept  de   structure
multiplicative.  L'action d'une structure multiplicative sur une note permet de
d\'efinir une tonalit\'e.
Le  classement des structures multiplicatives  est alors facile avec  les notions de type et de mode.
Nous d\'etaillons le cas des gammes \`a cinq et sept notes 
et proposons des exemples musicaux  pour toute une  vari\'et\'e de modes. 

\bigskip \smallskip   
{\bf \large Abstract}

\smallskip \noindent
We present an algebraic construction of music notes and  show how to associate them
in several ways to construct  music ranges.
Then a family of ranges emerge with a fixed number of notes~: two, three, five, seven, twelve, seventeen,  {\it etc.}
A classification  of these scales is simple with the concept 
of multiplicative structure. The action of such a multiplicative structure on a music note
introduces  the definition of tonality.
Multiplicative structure classification is then straightforward with the notions of type and mode.
We detail the case of music scales with five and seven notes and propose musical examples
for a wide range of modes.

\bigskip \smallskip  
\noindent {\bf  Mots-clefs :} Pythagore, type, mode, tonalit\'e, transposition. 

\newpage 
\bigskip  \bigskip
\noindent {\bf \large Introduction}

\smallskip \noindent
L'\^etre humain a toujours cherch\'e \`a imiter la musique offerte par la nature :
chant des oiseaux, souffle du vent, clapotis des vagues, {\it etc}.
La voix humaine permet d'ailleurs une infinit\'e de modulations et d'expressions sonores.
La musique peut aussi \^etre fabriqu\'ee avec des objets~:
un os perc\'e de trous 
bouch\'es \`a volont\'e permet d'\'emettre  diff\'erents sons  reproductibles,
des notes de musique.
Les notes de musique  diff\`erent par leur hauteur, leur puissance, leur timbre, leur dur\'ee.
De plus, des sons sp\'ecifiques sont parfois \'emis par des outils. Ainsi,
un simple arc met en tension une corde qui vibre lors du d\'epart de la fl\`eche
et peut devenir musical pour peu qu'on le dote d'un r\'esonateur,
comme on en constate encore l'usage aujourd'hui dans des soci\'et\'es dites primitives...
La hauteur du son \'emis d\'epend de la tension du fil et de la longueur de la corde.
Cet instrument de musique fondamental reste une r\'ef\'erence
pour la construction formelle de la musique avec le monocorde, un instrument \`a une seule corde.
La recherche d'instruments qui produisent des sons harmonieux, de la musique,
est une qu\^ete qui a travers\'e toutes les \'epoques et toutes les civilisations.

\bigskip \noindent
Une fois cr\'e\'ee, une  musique agr\'eable \`a l'oreille doit pouvoir \^etre reproduite.
D'o\`u une recherche d'un langage musical, avec une codification propre \`a
chaque culture.
Structurer cette musique \`a l'aide de mod\`eles math\'ematiques a probablement commenc\'e
avant m\^eme l'\'ecriture  et Marc Chemiller (2007) a explor\'e les difficult\'es
pour  expliciter ces math\'ematiques naturelles. 
Ainsi, Pythagore de Samos ($-582$-$-496$), en plus de son c\'el\`ebre th\'eor\`eme,
nous a transmis une codification des notes de musique et de leur assemblage, une gamme.
La ~\guillemotleft~gamme de Pythagore~\guillemotright~
permet de d\'evelopper un vocabulaire sp\'ecifique pour la musique
qui sert de r\'ef\'erence \`a l'\'ecriture musicale
durant tout le moyen-\^age.
Elle  est fond\'ee sur l'intervalle de quinte (entre le {\textsl {do}} et le {\textsl {sol}} typiquement)
et  fait appara\^itre la fraction 3/2. Nous y reviendrons.
Si aucun \'ecrit de Pythagore ne nous est parvenu,
Aristox\`ene de Tarente, philosophe grec du 4e si\`ecle avant  J\'esus-Christ, actif vers $-330$,
 est l'auteur d'un des premiers ouvrages sur l'harmonie musicale qui nous ait \'et\'e transmis.

\bigskip \noindent
Gioseffo Zarlino (1517-1590) publie en 1558 un trait\'e de musique 
qui fait la synth\`ese des connaissances
th\'eoriques et pratiques de son \'epoque. Surtout,
il remet en cause l'intervalle de tierce (intervalle entre le {\textsl {do}} et le {\textsl {mi}}  par exemple)
utilis\'e depuis Pythagore au b\'en\'efice de la tierce majeure caract\'eris\'ee par un rapport de hauteurs \'egal \`a~5/4.
%
Zarlino propose ainsi  une ~\guillemotleft~gamme naturelle~\guillemotright~ qui utilise la tierce majeure et la quinte.
Mais  l'\'emergence des instruments \`a clavier,
o\`u le musicien a acc\`es \`a un ensemble pr\'ed\'efini de notes de musique,
 a demand\'e  des compromis entre l'esth\'etique des accords et la complexit\'e de la gamme naturelle.

\fancyhead[EC]{\sc{Fran\c{c}ois Dubois}}
\fancyhead[OC]{\sc{Pourquoi la gamme a sept notes ?}} 
\fancyfoot[C]{\oldstylenums{\thepage}}

\bigskip \noindent
Ce compromis d\'elicat, ce choix du   ~\guillemotleft~temp\'erament~\guillemotright,
a donn\'e lieu \`a de tr\`es nombreux travaux th\'eoriques et pratiques depuis plusieurs si\`ecles.
Dans leur ouvrage didactique,
Claude  Abromont et  Eug\`ene de Montalembert (2001) pr\'esentent pr\`es de 600 trait\'es  sur la
th\'eorie de la musique en occident, certains propos\'es par des compositeurs majeurs ou des math\'ematiciens
tr\`es c\'el\`ebres par ailleurs.
Il est bien s\^ur hors de question d'\^etre exhaustif ici. Nous avons simplement
choisi quelques ouvrages qui nous semblent tr\`es importants.
Ainsi,
Marin Mersenne  (1588-1648) a l'intuition de la n\'ecessit\'e d'un temp\'erament \'egal (1636),
Simon Stevin (1548-1620), avec son trait\'e sur l'art de chanter retrouv\'e au 19e si\`ecle,
propose une division de l'octave en douze intervalles \'egaux,
tout comme son contemporain en Orient Zhu Zaiyu (1536-1611) avec un  ouvrage sur le temp\'erament \'egal 
(voir par exemple la th\`ese de Woo Shingkwan, 2017).
Enfin, Andreas Werckmeister (1645-1706),
propose (1691) le d\'eveloppement d'un temp\'erament \'egal
qui permet de transposer sans difficult\'e un morceau de musique d'une hauteur \`a une autre,
mais au d\'etriment de tierces souvent  bien loin de l'\'etat pur et de
quintes parfois tr\`es  fortement r\'eduites. 
Pour les liens entre  divers temp\'eraments et   ~\guillemotleft~Le clavier bien temp\'er\'e~\guillemotright~
de Jean-S\'ebastien Bach, nous renvoyons au travail de Johan Broekaert (2021).

\bigskip \noindent
Au  18e si\`ecle,
Jean-Philippe Rameau (1683-1764), c\'el\`ebre compositeur fran\c cais,  propose en 1722 son {\it  Trait\'e de l'harmonie}
qui sera suivi de plusieurs autres ouvrages.
La structure des accords de tierce et de quinte sert d'hypoth\`ese pour d\'evelopper une th\'eorie
de l'harmonie qui aura une grande influence.
Et Friedrich  Marpurg (1718-1795) lui r\'epond dans son trait\'e de  1776.
Au 19e si\`ecle,
Robert Bosanquet  (1841-1912) propose (1876) une distinction entre
le temp\'erament \'egal et un  temp\'erament bien temp\'er\'e.
Enfin, Hugo Riemann (1849-1919) publiera plus de 60 ouvrages, dont son c\'el\`ebre {\it Lexique musical} (1882).

\bigskip \noindent
Au 20e si\`ecle, les recherches sur la mise en forme de l'\'ecriture musicale
se poursuivent.
Le compositeur Arnold Sch\"onberg (1874-1951), inventeur du dod\'ecaphonisme, \'ecrit aussi un trait\'e
sur l'harmonie classique (1911) et
Ivan Wyschnegradsky (1893-1979) formalise l'emploi de quarts de tons (1932).
%
Alain Dani\'elou (1907-1994),  \'etudie la musique indienne sur le terrain et publie un trait\'e de
musicologie (1959),
Jacques Chailley (1910-1999) pr\'esente l'imbroglio des modes (1960) et
Iannis Xenakis (1922-2001) introduit de nombreux mod\`eles de la physique et des math\'ematiques
dans la formalisation de l'\'ecriture musicale~; un de ses ouvrages majeurs est publi\'e en 1963.
%

\bigskip \noindent
Dans son article sur les gammes naturelles (1999), Yves Hellegouarch (1936-2022) donne une formalisation math\'ematique
de la notion de temp\'erament avec des \'el\'ements de th\'eorie des groupes,
\'egalement pr\'esente dans la th\`ese de Moreno Andreatta (2003).
\`A partir d'une approche volontairement \'el\'ementaire, 
Michel Brou\'e (2002) montre que la gamme classique \'emerge des douze notes du syst\`eme temp\'er\'e
de la gamme occidentale. 
\`A l'aide d'une m\'ethode d'analyse spectrale originale, 
Thomas H\'elie a d\'evelopp\'e avec  ses coll\`egues de l'Institut de Recherche et Coordination Acoustique/Musique
le~\guillemotleft~snail~\guillemotright,
un outil logiciel et mat\'eriel d'analyse et de visualisation du son (2016).
Enfin, Andr\'e  Calvet propose (2020) une vaste fresque sur l'histoire de la musique et des temp\'eraments.

\bigskip \noindent
Suite \`a tous ces d\'eveloppements tr\`es \'erudits et subtils, il ne semble
pas raisonnable de proposer une autre  fa\c con de structurer divers \'el\'ements de l'\'ecriture de la musique avec un syst\`eme formel.
Mais sans avoir connaissance de tout ce corpus musicologique, 
nous avons propos\'e il y a une vingtaine d'ann\'ees 
une fa\c con assez simple de construire \`a la fois des notes et des gammes.
Nous la reproduisons dans les paragraphes qui suivent,  \`a  quelques variantes 
pr\`es de pr\'esentation et de r\'edaction.

\bigskip \noindent
Nous nous donnons d'abord une toute premi\`ere note de r\'ef\'erence (paragraphe 1).
\`A partir du premier  harmonique, de la quinte,
une premi\`ere gamme tr\`es primitive appa\-ra\^it, avec seulement deux notes (paragraphe 2).
Avec ces deux premi\`eres notes, la construction it\'erative de la suite des quintes au paragraphe 3 permet d'expliciter
une suite infinie de notes, en retrouvant essentiellement la construction initiale de Pythagore.
Repartant de la gamme \`a deux notes, de l'intervalle  de quinte et de son intervalle compl\'ementaire, la quarte,
il est alors naturel de diviser la quinte par la quarte
afin de faire appara\^itre une nouvelle gamme qui comporte trois notes (paragraphe 4).

\bigskip \noindent
On recommence ensuite ce  processus de  ~\guillemotleft~brisure de ton~\guillemotright~
et l'on coupe le plus grand des deux intervalles d'une gamme \`a trois
notes par le plus petit. On obtient  au paragraphe 5 une  famille de gammes \`a cinq notes.
Nous pr\'esentons ensuite la  construction r\'ecurrente non  lin\'eaire d'une famille infinie
de gammes (paragraphe 6).
\`A partir du principe de d\'ecoupage en deux parties du plus grand des intervalles d'une gamme
donn\'ee par le plus petit,
nous mettons en \'evidence quelques propri\'et\'es math\'ematiques
au paragraphe 7 et tentons de mettre un peu d'ordre dans la multiplicit\'e
de choix possibles aux paragraphes 8 et 9.
De fa\c con incidente, nous trouvons que le nombre de notes de la famille
de gammes qui suit les gammes pentatoniques n'est pas \'egal \`a 8,
nombre suivant 5 dans la suite de Fibonacci 1, 2, 3, 5, 8, {\it etc.}
et  serait apparu avec un algorithme lin\'eaire  de d\'ecoupage. 
Nous explicitons au dernier paragraphe  les 21 gammes \`a sept notes qui commencent par une  note de r\'ef\'erence,
avant une conclusion.
En Annexe, nous avons regroup\'e 
la preuve d'une proposition un peu technique
et une tentative de construction d'une gamme \`a huit notes. 


%
\bigskip \bigskip    \noindent {\bf \large    1) \quad  Note de r\'ef\'erence}

\smallskip   \noindent 
La musique est \'ecrite avec des notes. Du point de vue de la mod\'elisation
acoustique, une  ~\guillemotleft~note~\guillemotright~
est une variation  de pression au cours du temps
caract\'eris\'ee par  une fen\^etre temporelle qui d\'efinit la dur\'ee de la
note et par une fr\'equence caract\'eristique de la 
~\guillemotleft~hauteur~\guillemotright~  du son. Du point de vue
math\'ematique, notre  mod\`ele de note de musique
est  une fonction p\'eriodique du temps, de pulsation
$\, \omega \,$ appel\'ee abusivement dans la suite ~\guillemotleft~{fr\'equence $\omega$}~\guillemotright,
de p\'eriode $\, T = 2 \pi / \omega \, $  dupliqu\'ee quelques centaines de fois
pour les fr\'equences sonores habituelles. 
Nous ne nous int\'eressons pas ici au timbre de
l'instrument,   c'est-\`a-dire au  spectre du motif int\'erieur \`a la p\'eriode,
mais simplement \`a la fr\'equence de la note.

\bigskip \noindent
Pour nous, la premi\`ere note, la note de r\'ef\'erence, est le  {\textsl {do}}. Il s'agit d'une
convention et le  {\textsl {do}}  repr\'esente ici une note de fr\'equence $\, \omega_0 $,
avec $\, \omega_0 \,$ fix\'e pour toute la suite de cette contribution.  Nous
introduisons aussi un instrument de musique primitif,  la corde
vibrante de longueur~$\, L \,$ tenue \`a ses deux extr\'emit\'es. Nous supposons
que le   {\textsl {do}} correspond au premier mode de vibration de cette corde 
et   $\, \omega_0  \,$ est la fr\'equence fondamentale. 
Si une corde de longueur $\, L \,$ vibre \`a la fr\'equence $\, \omega_0 $,
son  second 
harmonique vibre \`a la fr\'equence $\, 2 \, \omega_0 \,$ et
correspond \`a une longueur de $\, {{L}\over{2}}  $.
Le troisi\`eme   harmonique
vibre \`a la fr\'equence $\, 3 \, \omega_0 \,$ et est associ\'e \`a une longueur de
$\, {{L}\over{3}} $.  De fa\c{c}on g\'en\'erale, pour $ \, n \, $
entier sup\'erieur ou \'egal \`a 1,
le  $\, n^{\rm o} $ harmonique  est de
fr\'equence $\, n \, \omega_0 \,$ (voir par exemple
le livre 
de  Laurent Schwartz, 1965) et 
correspond \`a une longueur de corde \'egale \`a $\, {{L}\over{n}} $.

\bigskip \noindent
L'oreille humaine est sensible aux rapports de fr\'equences.
Si ce rapport est \'egal \`a 2, les deux notes sont  \`a l'octave.
L'octave, ou  la  fr\'equence $2 \, \omega_0  $, 
correspond au second  harmonique 
pour la vibration d'une corde  si la
fr\'equence fondamentale est \'egale \`a~$\, \omega_0 $.
Par convention, nous appelons {\textsl {do}}$^* \, $
la note de fr\'equence   $2 \,  \omega_0 \, $:
\moneqstar
{\textsl {do}}^* \,  \,\, : \,\, 2 \,  \omega_0  \, .
\monendstar
Ensuite, il y a tous les autres {\textsl {do}}.
Tous ces {\textsl {do}} sont multiples du
pr\'ec\'edent  via une puissance de $2$  en termes de fr\'equence, le {\textsl {do}}$^*$
\`a l'octave, le suivant  \`a deux octaves  ($4 \, \omega_0 $) {\it etc.}    ou bien \`a
l'octave plus grave (${1\over2}  \, \omega_0 $), deux octaves (${1\over4} \, \omega_0
$), {\it etc.}   Nous retenons  
que la note de r\'ef\'erence 
de l'octave d'ordre $ \, k \, $ a une fr\'equence 
d\'efinie par
\moneqstar
\omega_k  =  2^k \, \omega_0 \,, \quad  k \in \ZZ \,.
\monendstar
%

\bigskip \noindent
On \'etudie ensuite  les autres notes et par convention on se place
dans l'intervalle $\, [ \, \omega_0 , \,2 \,  \omega_0  \,] , \,$ entre le   {\textsl {do}} 
et le {\textsl {do}}$^*  \, $ : le {\textsl {do}} est la premi\`ere note de l'intervalle  $\, [ \,
\omega_0 , \,2 \,  \omega_0  \,] \,  $ et le~{\textsl {do}}$^*$ la derni\`ere. On
introduit la notation $\, \Omega_0 \,$ utile dans la suite :
\moneqstar
 {\textsl {do}} \,: \,\, \Omega_0 =  \omega_0 \,.
\monendstar
Dans le syst\`eme classique qui remonte  
au moins \`a Pythagore 
et nous sert  ici de r\'ef\'erence, la note suivante est le {\textsl {sol}}. On g\'en\`ere le {\textsl {sol}} gr\^ace
au troisi\`eme harmonique,   
de fr\'equence $\, 3 \, \omega_0 $, pour une  corde de
fr\'equence fondamentale associ\'ee \`a la note  {\textsl {do}}. On ram\`ene ensuite cette nouvelle note dans
l'octave de r\'ef\'erence  $\, [ \, \omega_0 , \,2 \,  \omega_0  \,] $  par division de
la fr\'equence par deux, une transformation du boulanger multiplicative (voir par
exemple l'ouvrage de Vladimir Arnold et Andr\'e Avez, 1967). 
Le {\textsl {sol}} d\'efinit la notation~$\, \Omega_1 , \,
$ de fr\'equence \'egale \`a $ \, {3\over2} \, \omega_0 \,$ :
\moneqstar
{\textsl {sol}} \,: \,\, \Omega_1  =  {3\over2} \, \omega_0 \, . 
\monendstar

\bigskip \bigskip  \noindent {\bf \large    2) \quad  Gammes \`a deux notes}

\smallskip   \noindent 
Une ~\guillemotleft~gamme~\guillemotright~ est  l'ensemble des notes disponibles pour \'ecrire un morceau de
musique,  c'est-\`a-dire par convention l'ensemble des fr\'equences  $\, \nu_j \,$
utilisables \`a partir de la fr\'equence fondamentale $\, \omega_0 \,$
jusqu'au second harmonique $\, 2 \, \omega_0 \,$:
\moneqstar
\omega_0 \,\, \leq \,\, \nu_j \,\,  < \,\,  2 \, \omega_0 \,.
\monendstar
Une gamme contient $\,p \,$ notes si l'indice entier $\,j \,$ prend 
les valeurs 
$\, j = \, 0, \, 1 , \, 2, \, $...$ \, , \, p-1 $. On ordonne alors 
les fr\'equences $\, \nu_j \,$ en ordre croissant~:
\moneq \label{2.2}
\omega_0 = \nu_0 \, < \, \nu_1 \, < \, \nu_2 \, < \, \cdots \, < \, \nu_j \, < \,
\nu_{j+1} \, < \, \cdots \, < \,  \nu_{p-1} \, < \, 2 \omega_0 .
\monend

\smallskip  \noindent
Nous proposons dans ce qui suit une structure pour d\'ecrire toute une famille de gammes.
Nous ne parlons pas davantage d'une gamme \`a une seule note, form\'ee d'une simple octave : 
\moneqstar
{\cal G}^0_1  =  \bigl( {\textsl{do}},\,  {\textsl{do}}^* \bigr)  \, .  
\monendstar 

\smallskip \noindent 
Nous commen\c{c}ons une premi\`ere famille de gammes avec les trois notes dont nous disposons~:
le {\textsl {do}}, le {\textsl {sol}} et le  {\textsl {do}}$^*$. On fabrique ainsi une gamme ~\guillemotleft~primitive~\guillemotright~ \`a deux
notes
\moneq \label {2.3}
{\cal G}^1_1  =   \bigl(  {\textsl{do}},\,  {\textsl {sol}},\, {\textsl{do}}^* \bigr)  \,.
\monend

\smallskip \noindent
L'intervalle  $\,  {\textsl {sol}} \, /  {\textsl {do}} = {3\over2} \,$
d\'efinit une 
~\guillemotleft~{quinte}~\guillemotright~ pure 
et le second intervalle (plus petit)\br
$ {\textsl {do}}^* /  {\textsl {sol}}  =  {{2}\over{3/2}} = {4\over3} \, $ une
~\guillemotleft~{quarte}~\guillemotright~ pure.  
Il n'y a aucune raison pour r\'eduire une gamme primitive de deux
notes \`a l'unique choix $ \, {\cal G}^1_1 \, $propos\'e \`a la relation (\ref{2.3}),
 une quinte suivie d'une quarte. On peut faire aussi le contraire,  c'est-\`a-dire 
une quarte suivie d'une quinte. Ce faisant, on d\'efinit  \`a la fois une
nouvelle note et une nouvelle gamme. La nouvelle note, le {\textsl {fa}}, de fr\'equence
$\, \Omega_{-1} \, $ est issue de $\, \omega_0 \,$ par une quarte pure, ce qui
s'\'ecrit par d\'efinition~:
\moneqstar
{\textsl {fa}} \, : \, \Omega_{-1} = {4\over3} \, \omega_0 
 \approx  1.333333 \,\,  \omega_0  \, .
\monendstar
La seconde gamme primitive correspond donc \`a la s\'equence {\textsl {do}},  {\textsl {fa}},
{\textsl {do}}$^*$ :
\moneq \label {2.5}
{\cal G}^1_2 =  \bigl(  {\textsl{do}},\,  {\textsl {fa}},\, {\textsl{do}}^* \bigr) \,.
\monend

\smallskip \noindent
Nous disposons de deux ~\guillemotleft~{structures multiplicatives}~\guillemotright~   pour les gammes primitives :
la gamme~$\,  {\cal G}^1_1 \, $ avec (dans cet ordre), quinte et quarte et la
gamme $\,  {\cal G}^1_2 \, $ avec une  quarte suivie d'une quinte. Mais une
gamme, m\^eme primitive, doit-elle commencer par la note {\textsl {do}}~? La pratique de la
musique 
montre qu'il n'en est rien. Nous nous livrons donc dans le
paragraphe qui suit \`a la recherche d'une gamme ayant la m\^eme structure   que la
gamme  $\,  {\cal G}^1_1 \, $ et qui commence par un {\textsl {sol}}. La structure de
la gamme  $\,  {\cal G}^1_1 \, $ est par d\'efinition la suite des intervalles
multiplicatifs  entre les notes de la gamme, ici $\, ({\rm quinte}, {\rm quarte}) \,=
\, \bigl( {3\over2} ,\, {4\over3} \bigr) . \,$  Nous retenons la relation 
\moneq \label{2.6}
{\rm structure  \, multiplicative} \, \bigl( {\cal G}^1_1 \bigr) =
\Bigl( {3\over2} \, ,\, {4\over3} \Bigr) \,.
\monend

\bigskip \noindent
Nous fabriquons la gamme $  \, {\cal G}^1_1 \,$ qui commence par un {\textsl {sol}} en
posant $\, \nu_0 = \Omega_1 \,$ au lieu de $\, \nu_0 = \omega_0 \,$ \`a la relation
(\ref{2.2}), ce qui consiste
~\guillemotleft~{\`a changer le {\textsl {do}} en {\textsl {sol}}}~\guillemotright,
ou \`a changer la
fr\'equence de r\'ef\'erence pour la premi\`ere note de cette gamme.
Cette nouvelle gamme ~\guillemotleft~{transpos\'ee}~\guillemotright~
a la m\^eme structure multiplicative que la gamme $  \, {\cal G}^1_1 \,$
et  commence par un {\textsl {sol}}. Nous la notons~$  \, {\cal G}^1_1(\Omega_1) \,$ et
elle consiste en la succession suivante de fr\'equences~:
\moneq \label{2.7}
{\cal G}^1_1(\Omega_1)  =  \bigl( \Omega_1 , \,  {3\over2}\, \Omega_1
,\, 2 \, \Omega_1 \bigr) \,.
\monend
La note interm\'ediaire de fr\'equence $\, {3\over2}\, \Omega_1 \, $ entre le
{\textsl {sol}} et son harmonique {\textsl {sol}}$^* = 2 \, \Omega_1 \, $  est par d\'efinition
le {\textsl {r\'e}}$\,^*$ de fr\'equence $\,  {9\over4}   \, \omega_0  $.  Comme  $\, {9\over4}
\,$ n'appartient pas \`a l'intervalle $\, [1,\,2] \,$ mais \`a l'intervalle $\, [2,\,4] $,
on ram\`ene le {\textsl {r\'e}}$^*$  dans l'intervalle  $\, [\omega_0,\,2\omega_0] \,$
en divisant sa fr\'equence par deux. On d\'efinit ainsi  le {\textsl {r\'e}}, de fr\'equence $\,
\Omega_2 = {9\over8} \, \omega_0 \,$ :
\moneq \label{2.8}
{\textsl {r\'e}} \, : \,  \Omega_2 =  {9\over8} \, \omega_0 =  {{3^2}\over{2^3}} \, \omega_0
= 1.125 \,\,  \omega_0 \,\,.
\monend
La gamme de m\^eme structure que $\, {\cal G}^1_1 \,$ et qui commence par la
note {\textsl {sol}}, soit  $\, {\cal G}^1_1(\Omega_1) \, $ s'\'ecrit finalement, compte
tenu de (\ref{2.7}) et (\ref{2.8})~:
\moneq \label{2.9}
{\cal G}^1_1(\Omega_1) =  \bigl( \Omega_1 , \,  2\, \Omega_2 ,\, 2 \, \Omega_1 \bigr) \,.
\monend
Nous allons voir au paragraphe suivant que ce processus se poursuit. 

\newpage 
\bigskip \bigskip   \noindent {\bf \large    3) \quad   Suite infinie des notes }

\smallskip   \noindent 
Nous poursuivons l'id\'ee de changer de fr\'equence de r\'ef\'erence, mais \`a
partir du {\textsl {r\'e}} au lieu du {\textsl {sol}} lors de l'\'etude pr\'ec\'edente. Compte
tenu des relations (\ref{2.6}) et (\ref{2.9}), la gamme 
\moneqstar
{\cal G}^1_1(\Omega_2) =  \bigl( \Omega_2 , \,  {3\over2}\, \Omega_2 ,\, 2 \, \Omega_2 \bigr)
\monendstar
d\'efinit une nouvelle note de fr\'equence $\, \Omega_3 = {3\over2} \, \Omega_2
= {27\over16} \, \omega_0 , \,$ qui appartient \`a l'intervalle $\, [ \, \omega_0 , \,2 \,  \omega_0  \,] \,$ ;
c'est le {\textsl {la}}.
\moneqstar
{\textsl {la}} \,: \,  \Omega_3  =  {27\over16} \, \omega_0 = {{3^3}\over{2^4}} \, \omega_0 =  1.6875 \,  \omega_0 \, .
\monendstar

\smallskip \noindent
La suite des premi\`eres quintes $ \, \{ \Omega_0 ,\, \Omega_1 ,\, \Omega_2 ,\, \Omega_3 \} \, $
est d\'ej\`a d'une grande importance. Par exemple,  la famille
$ \,  \{ {\textsl {do}} ,\, {\textsl {sol}} ,\, {\textsl {r\'e}}^* ,\, {\textsl {la}}^* \} \, $
permet  d'accorder un violoncelle. 

\smallskip  \noindent
On poursuit maintenant la construction pr\'ec\'edente \`a l'infini. La gamme  $\,
{\cal G}^1_1(\Omega_k) \, $ a la m\^eme structure multiplicative que  $\,
{\cal G}^1_1 =  {\cal G}^1_1(\Omega_0) \, $ de la relation (\ref{2.6}) et commence par
la note~$\, \Omega_k \,$:  elle est {\it a priori}  de la forme $\, {\cal G}^1_1(\Omega_k) =
( \Omega_k ,\, {3\over2} \, \Omega_k ,\, 2 \, \Omega_k \, ) $.  Si
$\, {3\over2} \, \Omega_k \leq 2 \, \omega_0 ,\,$ on pose $\, \Omega_{k+1} = {3\over2}
\, \Omega_k \,$ et dans le cas contraire o\`u
$\, 2 \, \omega_0 < {3\over2} \, \Omega_k < 4 \, \omega_0 $, on d\'efinit la note $\, \Omega_{k+1} \,$ par la
relation  $\, \Omega_{k+1} = {3\over4} \, \Omega_k . \,$ De cette fa\c{c}on, on
construit une infinit\'e de notes de fr\'equence donn\'ee $\, \Omega_k \,$
(o\`u~$\,k\,$ est un entier sup\'erieur ou \'egal \`a~1)  avec une expression alg\'ebrique
de la  forme   $\, 3^k \, \omega_0 / 2^{\ell(k)} ,\,$ laquelle appartient par
d\'efinition \`a l'intervalle  $\, [ \, \omega_0 , \,2 \,  \omega_0  \,] . \,$ Pour
cr\'eer la note $\, \Omega_k ,\,$ on commence  par construire~$ \, k \, $ quintes
successives \`a partir du {\textsl {do}}  puis par mise \`a l'octave, on am\`ene cette note
dans l'intervalle de r\'ef\'erence~:
\moneq \label{3.2}
\Omega_k = {{3^k}\over{2^{\ell(k)}}} \, \omega_0 , \qquad \Omega_k
\in [ \, \omega_0 , \,2 \,  \omega_0  \,] \, .
\monend

\smallskip \noindent
Notons bien que  l'entier $\, \ell(k) \,$ est d\'efini par la condition (\ref{3.2}). On peut
encore \'ecrire\br
$\,\, 1 \leq {{3^k}\over{2^{\ell(k)}}}  < 2 \,$ ou encore de
fa\c{c}on \'equivalente $\,\, \ell(k) \, {\rm log}\,2 \, \leq \, k \, {\rm log}\,3 \, <
\,  \bigl( \ell(k) \!+ \! 1 \bigr) \,  {\rm log}\,2 $. Donc   l'entier $\,
\ell(k) \,$ est la partie enti\`ere de $\, k \, {\rm log}\,3 /  {\rm log}\,2 \,$ :
\moneq \label{3.3}
\ell(k) \,\leq \, k \, {{ {\rm log}\,3}\over{ {\rm log}\,2}} \, < \, \ell(k) + 1 \,\,; \quad
\ell(k) =  E \, \Big( k \, {{ {\rm log}\,3}\over{ {\rm log}\,2}} \Big) \,.
\monend
On introduit maintenant une  famille $\, \xi_k  \, $ de nombres rationnels  par les relations 
%
\moneq \label{7.9}
\xi_k = {{3^k}\over{2^{\ell(k)}}} \,,\,\,  k \in \ZZ \,.  
\monend
avec $ \, \ell(k) \, $ introduit \`a la relation (\ref{3.3}).
La note num\'ero $ \, k \, $ satisfait \`a la relation 
\moneqstar
\Omega_k = \xi_k \,\,  \omega_0 \,,\,\,  k \,\, {\textrm {nombre entier.}}
\monendstar

\smallskip \noindent
Nous  pouvons expliciter  les premi\`eres notes ainsi construites
au sein du tableau suivant.
Nous notons aussi la valeur en cent, \'echelle logarithmique 
propos\'ee par  Alexander John Ellis (1885), d\'efinie ici par la relation
\moneqstar 
    {\rm cents} \, (\Omega_k) = {\rm nombre \,\,  entier \,\,  le \,\,  plus \,\,  proche  \,\, de\,}
    \Big( {{1200}\over{\log(2)}} \, \log {{\Omega_k}\over{\Omega_0}} \Big)\,    
\monendstar
pour fixer les id\'ees.

\bigskip
\centerline { \begin{tabular}{|c|c|c|c|c|c|c|c|c|c|c|}    \hline
num\'ero de la note & nom de la note & $ \Omega_{\rm num \acute e ro} / \omega_0   $ & valeur approch\'ee & cents  \\   \hline
0 & {\textsl {do}}              &  $ 1        $      & 1        & 0 \\   \hline
1 & {\textsl {sol}}             &  $ 3 / 2    $     & 1.5       & 702  \\   \hline
2 & {\textsl {r\'e}}            &  $ 3^2 / 2^3 $     & 1.125    & 204  \\   \hline
3 & {\textsl {la}}              &  $ 3^3 / 2^4 $     & 1.6875    & 906    \\   \hline
4 & {\textsl {mi}}              &  $ 3^4 / 2^6 $     & 1.265625  & 408 \\   \hline
5 & {\textsl {si}}              &  $ 3^5 / 2^7 $     & 1.898437  & 1110      \\   \hline
6 & {\textsl {fa}}$\,\sharp$    &  $ 3^6 / 2^9 $     &  1.423828  & 612     \\   \hline
7 & {\textsl {do}}$\,\sharp$    &  $ 3^7 / 2^{11} $   &  1.067871  & 114   \\   \hline
8 & {\textsl {sol}}$\,\sharp$   &  $ 3^8 / 2^{12} $   & 1.601806  & 816   \\   \hline
9 &  {\textsl {r\'e}}$\,\sharp$ &  $ 3^9 / 2^{14} $   & 1.201355  & 318   \\   \hline
10 & {\textsl {la}}$\,\sharp$   &  $ 3^{10} / 2^{15} $ & 1.802032  & 1020 \\   \hline
11 &  {\textsl {mi}}$\,\sharp$  &  $ 3^{11} / 2^{17} $ & 1.351524  & 522  \\   \hline
12 &  {\textsl {si}}$\,\sharp$  &  $ 3^{12} / 2^{19} $ & 1.013643  & 23   \\   \hline
\end{tabular} }

\bigskip \noindent 
Les six premi\`eres notes sont clairement associ\'ees \`a la gamme de Pythagore
(voir Ellis, 1885).
On observe que les deux derni\`eres notes sont 
confondues dans la musique classique occidentale
avec le {\textsl {fa}} et le {\textsl {do}}
puisque d'une part $\, (\Omega_{11} / \omega_0 \approx 1.351524 \, $ (522 cents) 
alors que $ \, \Omega_{-1}  / \omega_0 \approx 1.333333  \,$ (498 cents) 
et d'autre part  $\, \Omega_{12}  / \omega_0 \approx 1.013643 \ $ (23 cents) alors que
$ \, \Omega_{0}  / \omega_0 = 1 \, $ (0 cent). 

\bigskip \noindent
La gamme $\,  {\cal G}^1_1 (\Omega_6) \,$ par exemple est \'egale \`a
$\, ( \Omega_6 ,\, 2 \Omega_7 ,\, 2 \Omega_6 ) = ( {\textsl {fa}}\sharp ,\, {\textsl {do}}\sharp^*
,\, {\textsl {fa}}\sharp^* ) \,$ (avec les notations  
habituelles). De mani\`ere g\'en\'erale, on a
\moneqstar
{\cal G}^1_1(\Omega_k) =  \bigl( \Omega_k , \, \{ \Omega_{k+1} \,\,
{\rm ou} \,\, 2 \, \Omega_{k+1}  \}  ,\, 2 \Omega_k \bigr) \,.
\monendstar
Le choix $\, \Omega_{k+1} \,$ comme seconde note de la gamme $\, {\cal G}^1_1(\Omega_k) \,$
a lieu lorsque $\, \Omega_{k+1} >  \Omega_k \,$
et le choix $\,  2 \, \Omega_{k+1} \,$ correspond au cas de figure
$ \, \Omega_{k+1}  < \Omega_k \,$ et dans ce cas  $\,  2 \, \Omega_{k+1} >  2 \, \omega_0 >  \Omega_k  $.

\bigskip \noindent
Il manque encore des notes \`a  notre catalogue. Nous repartons de la gamme
primi\-tive\br
de deux notes $\, {\cal G}^1_2 =  {\cal G}^1_2(\Omega_0) =  ( {\textsl {do}} ,\,  {\textsl {fa}}  ,\,  {\textsl {do}}^* )
 = ( \Omega_0 ,\, \Omega_{-1} ,\, 2\, \Omega_0 ) \,$
 propos\'ee \`a la relation (\ref{2.5}).
Au lieu de faire commencer la premi\`ere note de ce type de gamme
par un {\textsl {do}}, on peut la faire commencer par un {\textsl {fa}}. On d\'efinit pour cela
la structure multiplicative de $\,  {\cal G}^1_2 \,$ \`a savoir une quarte
puis une quinte :
\moneqstar
{\rm structure  \, multiplicative} \, \bigl( {\cal G}^1_2 \bigr) =
\Bigl( {4\over3} ,\,{3\over2}  \Bigr) \,.
\monendstar
La note qui suit $\, \Omega_{-1} \,$ dans la gamme $\,  {\cal G}^1_2(\Omega_{-1})
\, $ est $\, \Omega_{-2} = {4\over3} \,\Omega_{-1} \,$ : c'est le {\textsl {si}} b\'emol
ou   {\textsl {si}}$\,\flat$~:
\moneqstar
{\textsl {si}}\,\flat  \,: \, \Omega_{-2} =  {16\over9} \, \omega_0 =
{{2^4}\over{3^2}} \, \omega_0   \approx  1.777778 \,   \omega_0
\monendstar
et
\moneqstar
{\cal G}^1_2(\Omega_{-1}) = \bigl( {\textsl {fa}} ,\,  {\textsl {si}}\,\flat ,\,  {\textsl {fa}}^* \bigr)
=  \bigl( \Omega_{-1} ,\,  \Omega_{-2} ,\, 2 \, \Omega_{-1} \bigr) \,.
\monendstar
%

\bigskip \noindent
La construction de nouvelles notes entra\^{\i}ne celle de nouvelles gammes et
r\'ecipro\-quement les nouvelles gammes $\, {\cal G}^1_2(\Omega_{-k}) \,$
permettent d'introduire de nouvelles notes $\, \Omega_{-(k \!+ \!1)} . \,$ 
L'algorithme se poursuit  avec les quartes comme avec les quintes. La note $\, \Omega_{-k} \,$ ($k$
entier sup\'erieur ou \'egal \`a un) est issue de la note fondamentale {\textsl {do}} par
une succession de $\,k\,$ quartes ensuite ramen\'ees dans l'intervalle   $\, [ \,
\omega_0 , \,2 \,  \omega_0  \,[\,$  \`a l'aide d'une multiplication  par la puissance de $2$
qui convient.  La note $\, \Omega_{-k} \,$ est de la forme
$\, 2^{-\ell(-k)} \, \omega_0 / 3^k \,$ ce qui revient \`a g\'en\'eraliser la relation~(\ref{3.2})
aux entiers $k$ n\'egatifs et la relation (\ref{3.3}) s'applique encore.

\bigskip \noindent
Les nouvelles notes d'indice n\'egatif sont aussi utiles que celles d'indice
positif. 
Nous  pr\'esentons les premi\`eres dans un tableau.

\bigskip
\centerline { \begin{tabular}{|c|c|c|c|c|c|c|c|c|c|c|}    \hline
num\'ero de la note & nom de la note & $ \Omega_{\rm num \acute e ro} / \omega_0   $ & valeur approch\'ee & cents  \\   \hline
0 & {\textsl {do}}              &  $ 1        $      & 1          & 0 \\   \hline
-1 & {\textsl {fa}}              &  $ 4/3     $      & 1.333333   & 498 \\   \hline
-2 & {\textsl {si}}$\,\flat$     &  $ 2^4/3^2 $      & 1.77778   & 996 \\   \hline
-3 & {\textsl {mi}}$\,\flat$   &  $ 2^5 / 3^3 $     & 1.185185   & 294    \\    \hline
-4 & {\textsl {la}}$\,\flat$   &  $ 2^7 / 3^4 $     & 1.580247   & 792  \\   \hline
-5 & {\textsl {r\'e}}$\,\flat$ &  $ 2^8 / 3^5 $     &  1.053498   & 90  \\   \hline
-6 & {\textsl {sol}}$\,\flat$  &  $  2^{10} / 3^6 $    &  1.404664 & 588 \\   \hline
-7 & {\textsl {do}}$\,\flat$  &  $  2^{12} / 3^7 $    &  1.872885  & 1086 \\   \hline
-8 & {\textsl {fa}}$\,\flat$  &  $  2^{13} / 3^8 $    &  1.248590 & 384 \\   \hline
\end{tabular} }

\bigskip  \noindent
Les deux derni\`eres notes de cette suite, \`a savoir  
%
{\textsl {do}}$\,\flat$ et {\textsl {fa}}$\,\flat$
sont \'elimin\'ees classiquement au profit des notes  {\textsl {si}} et {\textsl {mi}}.  
En effet, on a d'une part $\, \Omega_{-7} / \omega_0 \approx 1.872885 \, $ (1086 cents) 
alors que $  \Omega_5  / \omega_0 \approx 1.898437 \,$ (1110 cents) 
et d'autre part
$\, \Omega_{-8}  / \omega_0 \approx 1.248590 \,$ (384 cents) 
se confronte \`a $ \, \Omega_4 / \omega_0 \approx 1.265625 \, $ (408 cents). 

\bigskip \noindent
On remarque que pour $\, k \in \{-1 ,\, \cdots ,\, 5 \} \, ,$ on a
\moneqstar
\Omega_k \, \sharp =  \Omega_{k+7} \,\,\,\,\, {\rm et} \,\,\,\,\, \Omega_k \, \flat = \Omega_{k-7} \, .
\monendstar
%
%
Il est donc naturel de poser ici
\moneqstar
\Omega_k \, \sharp \sharp =  \Omega_{k+14} \,\,\,\,\, {\rm et} \,\,\,\,\, \Omega_k \, \flat \flat =  \Omega_{k-14}
\,\,\,\,\,\,\, {\rm pour} \,\,  -1 \leq k \leq 5 \,.
\monendstar
%
%
Ces relations fournissent en pratique 14 nouvelles notes de musique, 7 avec des indices positifs~:

\bigskip
\centerline { \begin{tabular}{|c|c|c|c|c|c|c|c|c|c|c|}    \hline
num\'ero de la note & nom de la note & $ \Omega_{\rm num \acute e ro} / \omega_0   $ & valeur approch\'ee & cents  \\   \hline
13  & {\textsl {fa}}$\,\sharp \sharp $   &  $ 3^{13} / 2^{20} $     & 1.520465   &  735 \\   \hline
14  & {\textsl {do}}$\,\sharp \sharp $   &  $ 3^{14} / 2^{22} $     & 1.140349   &  227  \\   \hline
15  & {\textsl {sol}}$\,\sharp \sharp $  &  $ 3^{15} / 2^{23} $     & 1.710523   &  923   \\   \hline
16  & {\textsl {r\'e}}$\,\sharp \sharp $ &  $ 3^{16} / 2^{25} $   &   1.282892   &  431   \\   \hline
17  & {\textsl {la}}$\,\sharp \sharp $   &  $ 3^{17} / 2^{26} $   &  1.924338   &  1133   \\   \hline
18  &  {\textsl {mi}}$\,\sharp \sharp $  &  $ 3^{18} / 2^{28} $   &  1.443254   &  635   \\   \hline
19  &  {\textsl {si}}$\,\sharp \sharp $  &  $ 3^{19} / 2^{30}  $ &  1.082440   &  137    \\   \hline
\end{tabular} }

\bigskip

\newpage
et 7 avec des indices n\'egatifs : 

\bigskip
\centerline { \begin{tabular}{|c|c|c|c|c|c|c|c|c|c|c|}    \hline
num\'ero de la note & nom de la note & $ \Omega_{\rm num \acute e ro} / \omega_0   $ & valeur approch\'ee & cents  \\   \hline
-9  &  {\textsl {si}}$\,\flat \flat $   &  $ 2^{15} / 3^{9}   $ &  1.664787    &  882    \\   \hline
-10 &  {\textsl {mi}}$\,\flat \flat $   &  $ 2^{16} / 3^{10}  $ &  1.109858    &  180\\   \hline
-11 &  {\textsl {la}}$\,\flat \flat $   &  $ 2^{18} / 3^{11}   $ &  1.479810   &  678   \\   \hline
-12 &  {\textsl {r\'e}}$\,\flat \flat $ &  $ 2^{20} / 3^{12}  $ &  1.973081    &  1177  \\   \hline
-13 &  {\textsl {sol}}$\,\flat \flat $  &  $ 2^{21} / 3^{13}  $ & 1.315387     &  475   \\   \hline
-14 &  {\textsl {do}}$\,\flat \flat $   &  $ 2^{23} / 3^{14}   $ & 1.753849    &  973   \\   \hline
-15 &  {\textsl {fa}}$\,\flat \flat $   &  $ 2^{24} / 3^{15}  $ & 1.169233     &  271   \\   \hline
\end{tabular} }

\bigskip \noindent
Nous retenons que de fa\c con g\'en\'erale, pour $ \, k \, $ nombre entier positif ou n\'egatif,
la note $ \, \Omega_k \, $ est associ\'ee \`a une fr\'equence proportionnelle \`a $\, 3^k $,
divis\'ee par une puissance de 2 convenable de sorte que le quotient $ \,  {{\Omega_k}\over{\Omega_0}} \, $
appartienne \`a l'intervalle $ \, [1,\, 2[ $. 
Quand on range les notes $\, \Omega_{-15} \,$ \`a $ \, \Omega_{19} \,$ par
fr\'equences croissantes et non  plus par puissances de 3 croissantes, on trouve
\moneqstar \left \{ \begin {array}{cccccccccccccc}
{\textsl {do}} &    <     & {\textsl {si}}\,\sharp &     <    & 
{ \textsl {r\'e}}\,\flat  &    <     &  {\textsl {do}} \,\sharp &    <     & 
{\textsl {si}}\,\sharp \sharp  &    <     &  {\textsl {mi}} \,\flat  \flat &  
<     &    { \textsl {r\'e}} \\
{ \textsl {r\'e}}   &     <   &  {\textsl {do}} \,\sharp \sharp  &  
<    & {\textsl {fa}} \,\flat  \flat  &     <   &  {\textsl {mi}} \,\flat  &
  <   & { \textsl {r\'e}} \,\sharp &    <  &    {\textsl {fa}} \,\flat  &   <    &  {\textsl {mi}} \\
 {\textsl {mi}} &    <   &  { \textsl {r\'e}} \, \sharp \sharp  & 
<  &  {\textsl {sol}} \,\flat \flat   &   <    &   {\textsl {fa}} \\ 
  {\textsl {fa}}   &    <     &   {\textsl {mi}} \,\sharp   &    <     &
{\textsl {sol}} \, \flat  &    <     &   {\textsl {fa}}\, \sharp   &    < 
&  {\textsl {mi}}\,  \sharp \sharp  &    <     &  {\textsl {la}} \, \flat \flat  
&    <     & {\textsl {sol}} \\
 {\textsl {sol}}   &    <     &  {\textsl {fa}} \,\sharp \sharp   &    < 
&  {\textsl {la}}\, \flat &    <    &  {\textsl {sol}}\, \sharp  &    < 
&  {\textsl {si}}\, \flat  \flat  &    <    &  {\textsl {la}} \\
 {\textsl {la}}   &    <     &  {\textsl {sol}} \,\sharp \sharp   &    < 
&  {\textsl {do}}\, \flat \flat &    <    &  {\textsl {si}} \,\flat  &  
<    &  {\textsl {la}}\, \sharp   &    <    &  {\textsl {do}}\, \flat   &  
<    &  {\textsl {si}} \\
{\textsl {si}}   &    <     &  {\textsl {la}} \,\sharp \sharp   &    < 
& { \textsl {r\'e}} \, \flat \flat &    <    &  {\textsl {do}}^* .
 \end {array} \right. \monendstar 

\bigskip  \bigskip    \noindent {\bf \large    4) \quad   Brisure de ton et gammes \`a trois notes }

\smallskip   \noindent 
Nous disposons des deux gammes \`a deux
notes $\, {\cal G}^1_1 \, $ et  $\, {\cal G}^1_2 $. Comme nous avons
construit toute la famille de notes $\, \bigl(\Omega_k\bigr)  \, \ib{k \in \ZZ} $, 
nous  pouvons apr\`es transposition utiliser \'egalement les gammes $\,{\cal G}^1_1(\Omega_{k}) \,$
et $\,{\cal G}^1_2(\Omega_{k}) \, (k \in \ZZ) $. 
Toutefois, on peut avoir  envie  de
construire des gammes avec plus de deux notes. Nous proposons dans la suite  le
proc\'ed\'e r\'ecursif de ~\guillemotleft~{brisure de ton}~\guillemotright~ 
\`a partir des gammes \`a deux notes. 
%

\bigskip \noindent
Nous d\'efinissons le   ~\guillemotleft~{ton}~\guillemotright~ des gammes \`a deux notes et nous le notons $\,
\theta_1 \,$ ; il s'agit d'une quinte pure, comme propos\'e au paragraphe 2. Nous
avons
\moneqstar
\theta_1 =  {3\over2} \, .
\monendstar
Le ~\guillemotleft~{demi-ton}~\guillemotright~ des gammes \`a deux notes, not\'e $\, \delta_1 ,\,$ est par d\'efinition une quarte pure :
\moneqstar
\delta_1 = {4\over3} \,  .
\monendstar
Nous cherchons \`a former une nouvelle gamme issue d'une gamme \`a deux notes et
telle que le rapport de  fr\'equences de deux notes successives
$\, \nu_{j+1} \,/\, \nu_j \,$ arbitraires vaille soit un ton $\, \theta_2 \,$ de la nouvelle gamme, soit un
demi-ton $\, \delta_2 \,$ de celle ci :
\moneqstar
{{\nu_{j+1}} \over{\nu_j}} \, \in \, \{ \theta_2 ,\, \delta_2 \} \,.
\monendstar
%
Nous g\'en\'eralisons pour les nouvelles gammes ce qui existe pour la gamme \`a deux
notes~$\, {\cal G}^1_1 $~:
\moneqstar
\setbox21=\hbox { {\textsl {sol}} }
\setbox22=\hbox { {\textsl {do}} }
\setbox31=\hbox { {\textsl {do}}$^*$ }
\setbox32=\hbox { {\textsl {sol}} }
{{\box21} \over{{\box22} }}  , \,  {{\box31} \over{{\box32} }} \, \in \, \{ \theta_1 ,\, \delta_1 \} \,.
\monendstar

\bigskip \noindent
Nous brisons multiplicativement  le ton $\, \theta_1 \,$ \`a l'aide du demi-ton $\,
\delta_1 \,$ et formons ainsi un nouveau ton avec le demi-ton pr\'ec\'edent :
\moneqstar
\theta_2 = \delta_1  = {4\over3} \,. 
\monendstar
Un nouveau demi-ton $ \, \delta_2 \, $ est obtenu avec le rapport de fr\'equences li\'e \`a la brisure, 
 c'est-\`a-dire  le rapport  de l'ancien ton sur le nouveau ton :
\moneqstar
\delta_2 =  {{\theta_1}\over{\theta_2}} =
{{\theta_1}\over{\delta_1}} =  {9\over8} = {{3^2}\over{2^3}} \,.
\monendstar
En proc\'edant de cette mani\`ere, le ton primitif (une quinte pure !) est bris\'e en
un nouveau ton (une quarte pure) et un nouveau demi-ton  $\, \delta_2 \,$ qu'on
appelle classiquement une ~\guillemotleft~seconde~\guillemotright. 

\setbox21=\hbox { $  {\textsl {sol}}  $ }
\setbox22=\hbox{ $  {\textsl {do}}  $ }

\bigskip \noindent
La gamme $\, {\cal G}^1_1  \,$ donne naissance par brisure de son unique ton \`a
deux gammes de seconde g\'en\'eration~$\, {\cal G}^2_1 \,$ et $\, {\cal G}^2_2 \,$ selon la
fa\c{c}on dont la quinte $\, {\textsl {sol}} /  {\textsl {do}}  \,$  est d\'ecoup\'ee en
~\guillemotleft~{un demi-ton plus un ton}~\guillemotright.
Dans le premier cas, on a $\, {{\nu_1}\over{\nu_0}} =
\delta_2 = {9\over8} , \,$   $\, {{\nu_2}\over{\nu_1}} = \theta_2 = {4\over3} , \,$
 $\, {{\nu_3}\over{\nu_2}} = \theta_2 = {4\over3}  \,$;  on vient de cr\'eer une
gamme \`a trois notes, de seconde g\'en\'eration, 
qui comporte deux tons $ \, \theta_2 $ (une quarte pure) et un
demi-ton $\, \delta_2 \,$ (une  seconde)~:
\moneqstar
{\cal G}^2_1  =  \bigl( {\textsl {do}} , \,{\textsl {r\'e}}, \, {\textsl {sol}} ,\, {\textsl {do}}^* \bigr)
=   \bigl( \Omega_0 ,\, \Omega_2 ,\, \Omega_1 ,\, 2 \, \Omega_0 \bigr) \,.
\monendstar
La seconde gamme \`a trois notes $\, {\cal G}^2_2 \,$ est issue de la gamme
$\, {\cal G}^1_1, \,$ ce qui signifie que toutes les notes de la gamme
$\, {\cal  G}^2_2 \,$  sont \'egalement des notes de la gamme $\, {\cal G}^1_1 $.
On d\'ebute cette gamme par un ton~:
$\, {{\nu_1}\over{\nu_0}} = \theta_2 = {4\over3} $,
puis on la continue par un demi-ton~: $\, {{\nu_2}\over{\nu_1}} =
\delta_2 = {9\over8}  \,$ pour la  terminer  par un ton~:
$\, {{\nu_3}\over{\nu_2}} = \theta_2 = {4\over3} $. Le dernier ton de l'ancienne gamme devient le demi-ton de
la nouvelle. La gamme  $\, {\cal G}^2_2 \,$ est finalement donn\'ee par la suite
\moneq \label {4.9}
{\cal G}^2_2   =  \bigl( {\textsl {do}} , \,{\textsl {fa}}  , \, {\textsl {sol}} ,\, {\textsl {do}}^* \bigr)
=   \bigl( \Omega_0 ,\, \Omega_{-1} ,\, \Omega_1 ,\, 2 \, \Omega_0 \bigr) \,.
\monend

\bigskip \noindent
La gamme  $\, {\cal G}^1_2 \,$  
(relation (\ref{2.5}))  peut elle-aussi servir de
point de d\'epart \`a deux nouvelles gammes \`a trois notes : on coupe le ton
$\, \theta_1 = {\textsl {do}}^* / {\textsl {fa}} \,$ en un demi-ton $\, \delta_2 \,$ suivi d'un
ton  $\, \theta_2 \,$ ou bien en un ton $\,  \theta_2 \,$ suivi d'un demi-ton $\,
\delta_2 .\,$ Dans le premier cas, on obtient \`a nouveau la gamme $\, {\cal   G}^2_2 \,$
avec un demi-ton, un ton et un ton  
et dans le second, on a
 $\, {{\nu_1}\over{\nu_0}} = \theta_2 = {4\over3} , \,$  $\, {{\nu_2}\over{\nu_1}} =
 \theta_2 = {4\over3} , \,$  $\, {{\nu_2}\over{\nu_1}} = \delta_2 = {9\over8} . \,$
La note {\textsl {si}} b\'emol  est utile pour fabriquer $\, {\cal G}^2_3 \,$:  
\moneq \label {4.10}
{\cal G}^2_3  =  \bigl( {\textsl {do}} , \,{\textsl {fa}}  , \,{\textsl {si}}\,\flat ,\, {\textsl {do}}^* \bigr)
=  \bigl( \Omega_0 ,\, \Omega_{-1} ,\, \Omega_{-2} ,\, 2 \, \Omega_0 \bigr) \,.
\monend

\noindent
Nous retenons la structure multiplicative des trois gammes $\, {\cal G}^2_i\,$
$\, (i = 1, \, 2 ,\, 3)  \,$ :
%
\moneq \label{4.12}  \left\{  \begin {array}{l}
{\rm structure  \, multiplicative} \, \bigl( {\cal  G}^2_1  \bigr) =
\bigl( \delta_2 ,\, \theta_2 ,\,  \theta_2  \bigr) =
\bigl(  {9\over8} ,\,  {4\over3} ,\,  {4\over3}   \bigr) \\
{\rm structure  \, multiplicative} \, \bigl( {\cal  G}^2_2  \bigr) =
\bigl( \theta_2 ,\, \delta_2 ,\, \theta_2  \bigr) =
\bigl(  {4\over3} ,\,  {9\over8} ,\,  {4\over3}   \bigr) \\
{\rm structure  \, multiplicative} \, \bigl( {\cal  G}^2_3  \bigr) =
\bigl( \theta_2 ,\,  \theta_2  ,\, \delta_2  \bigr) =
\bigl(  {4\over3}  ,\,  {4\over3}  ,\,  {9\over8}  \bigr) \, .
\end{array} \right. \monend

\bigskip \noindent
Les gammes \`a trois notes permettent-elles de composer de la musique, m\^eme tr\`es primitive~?
Avec une gamme \`a trois notes, il est naturel de penser \`a sa mise en musique  avec une lyre
\`a quatre cordes.
Mais selon  Jean-Ren\'e Jannot (1979),  il n'existe pas de gamme 
t\'etracorde et on peut supposer qu'une lyre \`a quatre cordes \'etait un instrument d'accompagnement
pour une musique avec plus de notes. 
%
%
Par contre, le proc\'ed\'e de division d'une quinte par une quarte  propos\'e dans ce paragraphe
semble directement reli\'e \`a l'accord des instruments  hourrites dans la M\'esopotamie du second mill\'enaire avant J.-C.,  comme le propose 
Marcelle Duchesne-Guillemin (1980). En effet, \`a partir d'une   quinte  ascendandante
et d'une quarte  descendante, on construit simplement  une seconde par brisure  de ces deux intervalles.
Ainsi, la brisure de ton a donc probablement une origine pratique tr\`es ancienne.  
Nous poursuivons ce travail de brisure des tons par les demi-tons avec la troisi\`eme famille  de gammes,
les gammes pentatoniques.

\bigskip   \bigskip   \noindent {\bf \large    5) \quad   Gammes \`a cinq notes }

\smallskip   \noindent 
Le proc\'ed\'e de brisure de ton permet de construire  de nouvelles gammes~:
nous d\'ecoupons  le ton $\, \theta_2 = {4\over3} \,$ en un nouveau ton $\, \theta_3 \,$
et un  demi-ton $\, \delta_3 \,$ de sorte que
\moneqstar
\theta_3 \,\,  \delta_3 = \theta_2 \, . 
\monendstar
Si on proc\`ede   {\it a priori}  comme pour les gammes \`a trois notes,  on obtient $\, \theta_3 =
\delta_2 = {9\over8} \,$ et on obtient ainsi un nouveau demi-ton 
$\, \delta_3 = \theta_2 / \theta_3  = {{4 / 3 } \over{9 / 8}} = {32\over27} = {{2^5}\over{3^3}} $.
Mais les valeurs num\'eriques sont cruelles : on trouve $\, \theta_3 = {9\over8} = 1.125 \,$
et un demi-ton $ \, \delta_3 = {32\over27} \approx 1.1852 \,$
sup\'erieur strictement au ton  $\, \theta_3 $. 
Nous venons de former un ton qui est plus petit que le demi-ton
qui lui est associ\'e.  
On d\'ecide donc d'\'echanger les r\^oles du ton et du demi-ton,
d'accepter un algorithme de construction ``non lin\'eaire''  
pour que le ton  $ \, \theta_3 \, $ de la nouvelle famille soit plus grand que le demi-ton $ \, \delta_3 \, $ associ\'e :
\moneqstar
\delta_3  = \delta_2 =  {{3^2}\over{2^3}} \,\,\,\,\, {\rm et} \,\,\,\,\,
\theta_3 \,  = {{\theta_2}\over{\delta_3}} = {{\theta_2}\over{\delta_2}} =  {{2^5}\over{3^3}} \, .
\monendstar

On poursuit ensuite le proc\'ed\'e de brisure des deux tons des gammes $\, {\cal  G}^2_i  \,$
$\, (i = 1, \, 2 ,\, 3)  \,$ pour former la troisi\`eme famille de  gammes
$\, {\cal G}^3_j  $, \`a cinq notes.

\bigskip \noindent
Les gammes \`a cinq notes ont toutes deux tons et trois demi-tons :
\moneqstar
T_3 = 2 ,\,\qquad  D_3 = 3 \,
\monendstar
alors que les gammes \`a deux notes ont toutes un ton et un demi-ton :
\moneqstar
T_1 = 1 ,\,\qquad  D_1 = 1 \,
\monendstar
et les gammes \`a trois notes ont deux tons et un demi-ton :
\moneqstar
T_2 = 2 ,\,\qquad  D_2 = 1 \,.
\monendstar

\bigskip \noindent
Il y a quatre fa\c{c}ons de briser   les deux tons de la gamme $\, {\cal G}^2_1
\,$ :  le premier demi ton $\, {\textsl {r\'e}} / {\textsl {do}} = {{\nu_1}\over{\nu_0}}
\,$ est toujours un demi-ton $\, \delta_3 = \delta_2 . \,$ Si on d\'ecide que $\,
{{\nu_2}\over{\nu_1}} = \delta_3 , \,$ on introduit  la note {\textsl {mi}}. De m\^eme,
si on casse le second ton $\, \theta_2 \,$ en commen\c{c}ant par un demi-ton $\,
\delta_3 ,\,$ on utilise le {\textsl {la}} de fr\'equence $\, \Omega_3 .\,$ La gamme $\,
{\cal G}^3_1 \,$ \`a cinq notes ainsi obtenue s'\'ecrit :
\moneqstar
    {\cal G}^3_1 =  \bigl( {\textsl {do}} ,\, {\textsl {r\'e}} ,\, {\textsl {mi}}
    ,\, {\textsl {sol}} ,\, {\textsl {la}} ,\, {\textsl {do}}^*   \bigr) \,\,
 = \,\,  \bigl( \Omega_0 ,\, \Omega_2 ,\, \Omega_4 ,\, \Omega_1 ,\, \Omega_3
,\, 2 \, \Omega_0 \bigr) \,.
\monendstar
On l'appelle souvent gamme pentatonique majeure car on peut l'\'etendre en une gamme
de mode majeur avec sept notes, comme nous verrons plus loin. 
Si on brise le premier ton $\, {\textsl {sol}} / {\textsl {r\'e}} = \theta_2 \,$ de la
gamme  $\, {\cal G}^2_1  \,$ en commen\c{c}ant par un ton $\, \theta_3 ,\,$ on
retrouve la relation  $\, {{\nu_2}\over{\nu_1}} = \theta_3 = {32\over27} , \,$ soit
$\, {{\nu_2}\over{\omega_0}} = {9\over9} \, {32\over27} ={4\over3} = {{\Omega_{-1}}
\over{\omega_0}} , \,$ ce qui introduit \`a nouveau le~{\textsl {fa}}. De m\^eme, si on
brise le second ton $\,  {\textsl {do}}^* / {\textsl {sol}} \, $ de la gamme $\, {\cal
G}^2_1  \,$ en commen\c{c}ant par un ton nouveau, on a $\, {{\nu_4}\over{\nu_3}}
= \theta_3 \,$ donc $\, {{\nu_4}\over{\omega_0}} = {3\over2} \, {32\over27} ={8\over9}
= {{\Omega_{-2}} \over{\omega_0}}  \,$ et on reconna\^\i t le {\textsl {si}}$\,\flat$. Les
trois autres gammes \`a cinq notes issues de  $\, {\cal G}^2_1  \,$ s'en
d\'eduisent ais\'ement :
\moneqstar \left\{  \begin {array}{l}
  {\cal G}^3_2 =  \bigl( {\textsl {do}} ,\, {\textsl {r\'e}} ,\, {\textsl {fa}}
  ,\, {\textsl {sol}} ,\, {\textsl {la}} ,\, {\textsl {do}}^*   \bigr) \, \,
 = \,\,  \bigl( \Omega_0 ,\, \Omega_2 ,\, \Omega_{-1} ,\, \Omega_1 ,\,
\Omega_3 ,\, 2 \, \Omega_0 \bigr) \\
{\cal G}^3_3 =  \bigl( {\textsl {do}} ,\, {\textsl {r\'e}} ,\, {\textsl {fa}}
,\, {\textsl {sol}} ,\, {\textsl {si}}\,\flat ,\, {\textsl {do}}^*   \bigr) \, \,
 = \,\,  \bigl( \Omega_0 ,\, \Omega_2 ,\, \Omega_{-1} ,\, \Omega_1 ,\,
\Omega_{-2} ,\, 2 \, \Omega_0 \bigr) \\
{\cal G}^3_4 =  \bigl( {\textsl {do}} ,\, {\textsl {r\'e}} ,\, {\textsl {mi}}
,\, {\textsl {sol}} ,\, {\textsl {si}}\,\flat ,\, {\textsl {do}}^*   \bigr) \,\,
 = \,\,  \bigl( \Omega_0 ,\, \Omega_2 ,\, \Omega_{4} ,\, \Omega_1 ,\,
\Omega_{-2} ,\, 2 \, \Omega_0 \bigr) \,.
\end{array} \right. \monendstar

\bigskip \noindent
Combien y a-t-il de gammes \`a cinq notes ? Autant que de fa\c{c}ons de placer
arbitrairement les deux ($=T_3$) tons $\, \theta_3 \,$ et les trois  ($=D_3$)
 demi-tons $\, \delta_3 \,$ parmi cinq objets, puisque 
\moneqstar
\theta_3^{  \, T_3} \,\, \delta_3^{  \, D_3} =
\biggl( {{2^5}\over{3^3}} \biggr)^{ \scriptstyle \!\! 2}  \,  \biggl(
{{3^2}\over{2^3}} \biggr)^{ \scriptstyle \!\! 2} = 2 \,.
\monendstar
On compte par cons\'equent un total de
$ \, \big( ^{5} _{3} \big) = \big( ^{5} _{2} \big) = 10 \, $ gammes \`a cinq notes
qui d\'ebutent par la note {\textsl {do}}.

\bigskip \noindent
Examinons les gammes nouvelles issues de la brisure du ton de la gamme  $\, {\cal
G}^2_2  \,$ (relation (\ref{4.9})) o\`u  les deux tons $\, \theta_2 \,$ sont aux
extr\'emit\'es. Si on commence par un demi-ton $\, \delta_3 ,\,$ on a\br
$\, {{\nu_1}\over{\omega_0}} = \delta_3 =  {9\over8} = {{\Omega_{2}} \over{\omega_0}}  \,$
et on retrouve le {\textsl {r\'e}} comme premi\`ere note alors que si on commence par un ton
$\, \theta_3 ,\,$ on obtient  $\, {{\nu_1}\over{\omega_0}} = \theta_3 =
{{2^5}\over{3^3}} = {{\Omega_{-3}} \over{\omega_0}}  \,$ et le {\textsl {mi}}$\,\flat$  doit
\^etre utilis\'e ici. De m\^eme si on brise le second ton $\, {\textsl {do}}^* / {\textsl {sol}}
\,$ de la gamme  $\, {\cal G}^2_2  \,$ en commen\c{c}ant par un demi-ton $\,
\delta_3 ,\,$ on retrouve le {\textsl {la}}, alors que si on commence par un ton $ \, \theta_3 $,
on doit utiliser le  {\textsl {si}}$\,\flat$ : c'est analogue \`a la brisure du  second
ton de la gamme  $\, {\cal G}^2_1  \, $ et nous l'avons d\'ej\`a vu plus haut. On dispose
maintenant de deux nouvelles gammes \`a cinq notes :
\moneqstar \left\{  \begin {array}{l}
{\cal G}^3_5 =  \bigl( {\textsl {do}} ,\, {\textsl {mi}}\,\flat ,\,   {\textsl {fa}} 
,\, {\textsl {sol}} ,\, {\textsl {la}} ,\, {\textsl {do}}^*   \bigr) \, \,
 = \,\,  \bigl( \Omega_0 ,\, \Omega_{-3} ,\, \Omega_{-1} ,\, \Omega_1 ,\,
\Omega_3 ,\, 2 \, \Omega_0 \bigr) \\
{\cal G}^3_6 =  \bigl( {\textsl {do}} ,\, {\textsl {mi}}\,\flat ,\,  {\textsl {fa}} 
,\, {\textsl {sol}} ,\, {\textsl {si}}\,\flat  ,\, {\textsl {do}}^*   \bigr) =  \bigl(
\Omega_0 ,\, \Omega_{-3} ,\, \Omega_{-1} ,\, \Omega_1 ,\, \Omega_{-2} ,\, 2
\, \Omega_0 \bigr) \,.
\end{array} \right. \monendstar
%
La gamme $ \, {\cal G}^3_6  \, $ s'appelle aussi gamme pentatonique mineure. 
%

\bigskip \noindent
Nous continuons avec les gammes \`a cinq notes issues de la gamme
$\, {\cal G}^2_3  \,$ (relation (\ref{4.10})) qui a ses deux tons $\, \theta_2 \,$ au d\'ebut.
Pour briser le premier ton, le processus est analogue \`a la gamme  $\, {\cal G}^2_2  \,$
et pour briser le second ton $\, {\textsl {si}}\,\flat / {\textsl {fa}} ,\,$ on peut commencer par un
demi-ton $\, \delta_3 \,$ et on obtient un  {\textsl {sol}} et par cons\'equent une gamme
d\'ej\`a  construite (${\cal G}^3_3  \,$ ou $\, {\cal G}^3_6 $)
ou bien on commence par un ton $\, \theta_3 \,$ et on a alors $\,
{{\nu_4}\over{\nu_3}} = \theta_3 = {32\over27} ,\,$ donc $\, {{\nu_4}\over{\omega_0}}
=  {4\over3} \, {32\over27} = {64\over81} = {{2^7}\over{3^4}}  = {{\Omega_{-4}}
\over{\omega_0}}  \,$ ; on on ne coupera pas au $\,{\textsl {la}}\,\flat \,$ cette fois. On a
en d\'efinitive :
\moneqstar \left\{  \begin {array}{l}
{\cal G}^3_7 =  \bigl( {\textsl {do}} ,\, {\textsl {r\'e}} ,\,   {\textsl {fa}} 
,\, {\textsl {la}}\,\flat ,\, {\textsl {si}}\,\flat ,\, {\textsl {do}}^*   \bigr) =
\bigl( \Omega_0 ,\, \Omega_2 ,\, \Omega_{-1} ,\, \Omega_{-4} ,\, \Omega_{-2}
,\, 2 \, \Omega_0 \bigr) \\
{\cal G}^3_8 =  \bigl( {\textsl {do}} ,\, {\textsl {mi}}\,\flat  ,\,   {\textsl {fa}} 
,\, {\textsl {la}}\,\flat ,\, {\textsl {si}}\,\flat ,\, {\textsl {do}}^*   \bigr) =
\bigl( \Omega_0 ,\, \Omega_{-3} ,\, \Omega_{-1} ,\, \Omega_{-4} ,\,
\Omega_{-2} ,\, 2 \, \Omega_0 \bigr) \,.
\end{array} \right. \monendstar

\bigskip \noindent
Le d\'ecompte des gammes comportant cinq notes n'est pas complet. Les huit gammes
d\'ej\`a vues sont issues des trois gammes \`a trois notes. Nous d\'etaillons
dans le tableau qui suit la succession des deux tons et des trois demi-tons,  c'est-\`a-dire la structure
multiplicative de ces huit gammes.

\bigskip
\centerline { \begin{tabular}{|c|c|c|}    \hline
nom de la gamme  & alternance des tons et des demi-tons  & valeurs num\'eriques  \\   \hline
${\cal G}^3_1$ & $ \delta_3 ,\, \delta_3 ,\, \theta_3, \, \delta_3 ,\, \theta_3 $
&  $ {{9}\over{8}} ,\, {{9}\over{8}} ,\, {{32}\over{27}} ,\, {{9}\over{8}} ,\, {{32}\over{27}}$  \\   \hline
${\cal G}^3_2$ & $ \delta_3 ,\, \theta_3 ,\, \delta_3, \, \delta_3 ,\, \theta_3 $
&  $ {{9}\over{8}} ,\, {{32}\over{27}} ,\, {{9}\over{8}} ,\,  {{9}\over{8}} ,\, {{32}\over{27}}$  \\   \hline
${\cal G}^3_3$ & $ \delta_3 ,\, \theta_3 ,\, \delta_3 ,\, \theta_3 ,\, \delta_3  $
&  $ {{9}\over{8}} ,\, {{32}\over{27}} ,\, {{9}\over{8}} ,\, {{32}\over{27}} ,\,  {{9}\over{8}} $  \\   \hline
${\cal G}^3_4$ & $ \delta_3 ,\, \delta_3 ,\, \theta_3 ,\, \theta_3 ,\, \delta_3  $
&  $ {{9}\over{8}} ,\, {{9}\over{8}} ,\, {{32}\over{27}} ,\, {{32}\over{27}} ,\,  {{9}\over{8}}  $  \\   \hline
${\cal G}^3_5$ & $ \theta_3 ,\, \delta_3 ,\, \delta_3 ,\, \delta_3 ,\, \theta_3  $
&  $ {{32}\over{27}} ,\, {{9}\over{8}} ,\, {{9}\over{8}} ,\,{{9}\over{8}} ,\, {{32}\over{27}} $  \\   \hline
${\cal G}^3_6$ & $ \theta_3 ,\, \delta_3 ,\, \delta_3 ,\, \theta_3 ,\, \delta_3  $
&  $ {{32}\over{27}} ,\, {{9}\over{8}} ,\, {{9}\over{8}} ,\, {{32}\over{27}} ,\, {{9}\over{8}} $  \\   \hline
${\cal G}^3_7$ & $ \delta_3 ,\, \theta_3 ,\, \theta_3 ,\,  \delta_3 ,\, \delta_3  $
&  $ {{9}\over{8}} ,\,  {{32}\over{27}} ,\, {{32}\over{27}} ,\, {{9}\over{8}} ,\, {{9}\over{8}}  $  \\   \hline
${\cal G}^3_8$ & $ \theta_3 ,\, \delta_3 ,\, \theta_3 ,\,  \delta_3 ,\, \delta_3  $
&  $ {{32}\over{27}} ,\,  {{9}\over{8}} ,\,  {{32}\over{27}} ,\, {{9}\over{8}} ,\, {{9}\over{8}}  $  \\   \hline
\end{tabular} }

\bigskip \noindent
Il manque clairement deux gammes qui ne sont pas issues de gammes \`a trois notes,
mais qu'on peut construire simplement en compl\'etant les structures multiplicatives
pr\'ec\'edentes. On trouve ainsi :
\setbox30= \vbox {\halign{#&#&#&#&#&#&#&#&#&#&#&#&#&#&#&#&#&#&#&#&#&#&#&#&#&#&#  \cr
$ \!\! (5.24)  $ & $\qquad  {\rm structure} \, \bigl( {\cal G}^3_9 \bigr) \,$
& $\,\, \, = \,\,\, $ &   $ \,\, \Bigl( \, \theta_3 $&$,\,$&$ \theta_3 $&$,\,$&$
 \delta_3 $&$,\,$&$ \delta_3 $&$,\,$&$ \delta_3 \, \Bigr) \,\,\,$ & = &$ \,\,\,
\Bigl( \,$&${{32}\over{27}} $&$,\,$&$ {{32}\over{27}}  $&$,\,$&$  {{9}\over{8}}
$&$,\,$&$  {{9}\over{8}}  $&$,\,$&$  {{9}\over{8}}  \, $&$ \Bigr) $ & \cr
$ \!\! (5.25)  $ & $\qquad  {\rm structure} \, \bigl( {\cal G}^3_{10} \bigr) \,$
& $\,\, \, = \,\,\, $ &   $ \,\, \Bigl( \, \delta_3 $&$,\,$&$\delta_3  $&$,\,$&$
\delta_3  $&$,\,$&$ \theta_3  $&$,\,$&$ \theta_3 \, \Bigr)  \,\,\,$ & = &$ \,\,\,
\Bigl( \,$&$ {{9}\over{8}} $&$,\,$&$ {{9}\over{8}} $&$,\,$&$  {{9}\over{8}}
$&$,\,$&$ {{32}\over{27}}  $&$,\,$&$  {{32}\over{27}}  \, $&$ \Bigr) $ & . &
\cr }}
\setbox31= \hbox{ $\vcenter {\box30} $}
\setbox44=\hbox{\noindent $\box31  $}
%
\moneq \label {5.24-et-5.25}  \left\{  \begin {array}{l}
{\cal G}^3_9 \, : \, \big(  \delta_3 ,\, \delta_3 ,\, \delta_3 ,\, \theta_3 ,\, \theta_3  \big) =
\Big( {{32}\over{27}} ,\, {{32}\over{27}} ,\, {{9}\over{8}} ,\, {{9}\over{8}} ,\, {{9}\over{8}}  \Big) \\
{\cal G}^3_{10} \, : \, \big( \theta_3 ,\, \theta_3 ,\, \delta_3 ,\, \delta_3 ,\, \delta_3 \big) =
\Big( {{9}\over{8}} ,\, {{9}\over{8}} ,\, {{9}\over{8}}  ,\, {{32}\over{27}} ,\, {{32}\over{27}}  \Big) \,.
\end {array} \right.
\monend
Les  successions de notes des deux gammes de structures d\'ecrites \`a la relation  (\ref{5.24-et-5.25}) qui
commencent par un {\textsl {do}} sont faciles \`a d\'eterminer. Nous y rencontrons le {\textsl {sol}}
b\'emol ($ \Omega_{-6} $) 
et le  {\textsl {fa}} di\`ese (ou $\,\Omega_6 $) 
\`a nouveau :
\moneqstar  \left\{  \begin {array}{l}
{\cal G}^3_9 =  \bigl( {\textsl {do}} ,\, {\textsl {mi}}\,\flat ,\, {\textsl {sol}}\,\flat
,\, {\textsl {la}}\,\flat ,\,  {\textsl {si}}\,\flat ,\, {\textsl {do}}^*   \bigr) =
\bigl( \Omega_0 ,\, \Omega_{-3} ,\, \Omega_{-6}  ,\, \Omega_{-4} ,\,
\Omega_{-2} ,\, 2 \, \Omega_0 \bigr) \\
{\cal G}^3_{10}  =   \bigl( {\textsl {do}} ,\, {\textsl {r\'e}} ,\, {\textsl {mi}} 
,\, {\textsl {fa}}\,\sharp ,\, {\textsl {la}} ,\, {\textsl {do}}^*   \bigr) =
\bigl( \Omega_0 ,\, \Omega_{2} ,\, \Omega_{4}  ,\, \Omega_{6} ,\,
\Omega_{3} ,\, 2 \, \Omega_0 \bigr)  \, .
\end {array} \right. \monendstar

\bigskip \noindent
Les gammes pentatoniques sont bien r\'epertori\'ees 
depuis de nombreuses ann\'ees.
Ainsi, Hermann  von Helmholtz pr\'esente  dans son trait\'e  (1868)  les cinq gammes suivantes :
\moneqstar  \left\{  \begin {array}{l}
 {\cal H}_1 = \bigl(  {\textsl {do}} ,\, {\textsl {r\'e}} ,\,  {\textsl {fa}}   ,\,  {\textsl {sol}}  ,\,  {\textsl {si}}\,\flat ,\,  {\textsl {do}}^* \bigr)  \\
 {\cal H}_2 = \bigl(  {\textsl {fa}} ,\, {\textsl {sol}}  ,\,  {\textsl {si}}\,\flat  ,\, {\textsl {do}}^*   ,\,   {\textsl {r\'e}}^* ,\,  {\textsl {fa}}^* \bigr)  \\
 {\cal H}_3 = \bigl(  {\textsl {sol}} ,\,  {\textsl {si}}\,\flat ,\,  {\textsl {do}}^*  ,\, {\textsl {r\'e}}^* ,\,  {\textsl {fa}}^* ,\, {\textsl {sol}}^*  \bigr)  \\
 {\cal H}_4 = \bigl( {\textsl {si}}\,\flat ,\,  {\textsl {do}}^*  ,\, {\textsl {r\'e}}^* ,\,  {\textsl {fa}}^* ,\, {\textsl {sol}}^*  ,\,    {\textsl {si}}\,\flat^*  \bigr) \\ 
 {\cal H}_5 = \bigl( {\textsl {r\'e}} ,\,  {\textsl {fa}} ,\, {\textsl {sol}}  ,\,    {\textsl {si}}\,\flat^* ,\,  {\textsl {do}}^* ,\, {\textsl {r\'e}}^*  \bigr) \, . 
\end {array} \right. \monendstar
Nous pouvons les relier par transposition aux gammes pr\'esent\'ees plus haut et on a
\moneqstar 
 {\cal H}_1 = {\cal G}^3_3 ,\, 
 {\cal H}_2 = {\cal G}^3_2 ( {\textsl {fa}}) ,\, 
 {\cal H}_3 = {\cal G}^3_6 ( {\textsl {sol}}) ,\, 
 {\cal H}_4 = {\cal G}^3_1 ( {\textsl {si}}\,\flat ) ,\, 
 {\cal H}_5 = {\cal G}^3_8 (  {\textsl {r\'e}} ) .
 \monendstar
 Nous verrons plus loin dans cette contribution que ces cinq gammes sont toutes du m\^eme  
~\guillemotleft~{type}~\guillemotright. 
%
%
Dans son article sur le pentatonisme, Fran\c cois Picard (2001) pr\'esente cinq structures
associ\'ees \`a la musique chinoise.
Nous reprenons l'essentiel de son tableau et notons en regard la gamme correspondante avec notre classification. 

\bigskip 
\centerline { \begin{tabular}{|c|c|c|}    \hline
nom chinois  & gamme pr\'esent\'ee par F. Picard  & gamme associ\'ee    \\   \hline
gong diaoshi &  {\textsl {do,  r\'e,   mi,  sol,  la, do$^*$}}  & $ {\cal G}_1^3 $ \\ 
shang diaoshi   &  {\textsl {r\'e,   mi,  sol,  la, do$^*$, r\'e$^*$}}  & $ {\cal G}_3^3 ( {\textsl {r\'e}}) $ \\ 
jue diaoshi     &  {\textsl {mi,  sol,  la, do$^*$, r\'e$^*$, mi$^*$}}  & $ {\cal G}_8^3 ( {\textsl {mi}}) $ \\ 
zhi diaoshi     &  {\textsl {sol,  la, do$^*$, r\'e$^*$, mi$^*$,  sol$^*$ }}  & $ {\cal G}_2^3 ( {\textsl {sol}}) $ \\ 
yu  diaoshi     &  {\textsl {la, do$^*$, r\'e$^*$, mi$^*$,  sol$^*$, la$^*$  }}  & $ {\cal G}_6^3 ( {\textsl {la}}) $  \\ 
  \hline \end{tabular} }

\bigskip
Nous constatons que ces cinq gammes sont des transpositions des cinq gammes propos\'ees par Helmholtz. 
Les diverses games pentatoniques sont pr\'esentes dans de nombreuses musiques (voir Fran\c cois Picard, 2001), 
y compris dans la musique de Blues, avec de nombreuses varia\-tions.
Nous formalisons au paragraphe suivant l'existence d'une 
suite infinie de  familles  de gammes.

\bigskip  \bigskip    \noindent {\bf \large    6) \quad    Construction r\'ecurrente}

\smallskip \noindent 
Nous d\'etaillons   le processus  
pour passer des gammes de
g\'en\'eration $\,k \,$ (avec une notation g\'en\'erique $\,{\cal G}^k_\star \,$)
aux gammes de la g\'en\'eration suivante d'indice $\, k\!+\!1 .\,$ Une gamme  $\,{\cal G}^k_\star
\,$ contient $\,T_k \,$ tons et $\, D_k\,$ demi-tons, de valeurs respectives $\,
\theta_k \,$ et $\, \delta_k . \,$ L'ensemble des notes construites
se situent au sein de  l'octave de base  $\, [ \, \omega_0 , \,2 \,  \omega_0  \,] , \,$ et on
en d\'eduit :
\moneq \label {6.1}
\theta_k^{  \, T_k} \,\, \delta_k^{  \, D_k}  = 2 \, .
\monend
On appelle $\, p_k \,$ le nombre de notes d'une gamme de la $k^{\rm \,i\grave eme}\,$
famille et l'on a~:
\moneq \label {6.2}
p_k =    T_k \,+\, D_k  \, .
\monend

\smallskip  \noindent
La construction consiste toujours \`a briser le ton $\, \theta_k \,$ en un nouveau
ton $\, \theta_{k+1} \,$ \guillemotleft~{plus}~\guillemotright~  
un nouveau demi-ton $\, \delta_{k+1} \,$
avec l'un des deux morceaux \'egal au demi-ton~$\, \delta_k \,$ :
\moneq \label {6.3}
\theta_{k+1}  \,\,  \delta_{k+1}  =  \theta_k \,.
\monend

\bigskip \noindent
Si 
on se donne les deux relations suivantes
$ \,  \theta_{k+1}  = \delta_k \, $ et $ \, \delta_{k+1} =  {{\theta_k}\over{\delta_k}} $,
les plus   simples d'un point de vue alg\'ebrique, 
alors on peut aussi \'ecrire
$ \, \delta_k =  \theta_{k+1}  \, $ et $ \, \theta_k =  \theta_{k+1}  \, \,  \delta_{k+1} $.
Dans ces conditions,  la relation (\ref{6.1}) qui entra\^ine $ \, p_k = T_k + D_k \, $  s'\'ecrit maintenant 
$ \, (\theta_{k+1})^{T_k + D_k} \,\, (\delta_{k+1})^{T_k} = 2 $.
On en d\'eduit $\,  T_{k+1} = T_k + D_k $,  $\,  D_{k+1} = T_k \, $ et 
$ \, p_{k+1} = 2 \, T_k + D_k $.
\`A l'ordre suivant, on a donc
$\,  T_{k+2} =  T_{k+1} + D_{k+1} =  2 \, T_k + D_k \, $
et $ \,  D_{k+2} = T_{k+1} = T_k + D_k $.
Donc $ \,  p_{k+2} =  T_{k+2} +  D_{k+2} =  3 \, T_k + 2 \, D_k \, $
est exactement \'egal \`a $ \,  p_k + p_{k+1} $. Comme $\, p_0 = 1 $, $ \, p_1 = 2 $ et $ \, p_2 = 3 $,
la suite~$ \, p_k \, $ co\"incide avec la  suite de Fibonacci (voir par exemple Michel Demazure, 1997)
et on a $ \, p_4 = 8 \, $ avec cette hypoth\`ese.
Nous proposons un exemple de ce cas de figure, une gamme avec huit notes, \`a l'Annexe A.

\bigskip \noindent
Il faut bien prendre garde
au fait que le ton doit toujours rester sup\'erieur au demi-ton.
On choisit donc de garder comme nouveau ton le demi-ton $\, \delta_k \,$ si $\,
{{\theta_k}\over{\delta_k}} < \delta_k \,$ ou bien on pose $\, \theta_{k+1} =
{{\theta_k}\over{\delta_k}} \,$ dans le cas contraire,
%
\moneq \label {6.4}
 \theta_{k+1}  = \left \{ \begin {array}{l}
{\delta_k} \quad   {\rm si} \quad  {{\theta_k}\over{\delta_k}} < \delta_k     \\  \vspace{-.4 cm}  \\
 {{\theta_k}\over{\delta_k}}   \quad  {\rm si} \quad  {{\theta_k}\over{\delta_k}} >  \delta_k  \, . \end {array} \right.
\monend
%
Le demi-ton compl\'ementaire $\, \delta_{k+1} \,$ est issu  simplement des
relations (\ref{6.3}) et (\ref{6.4}) :
\moneq \label {6.5}
\delta_{k+1} = \left \{ \begin {array}{l}
 {{\theta_k}\over{\delta_k}}   \quad    {\rm si } \quad   {{\theta_k}\over{\delta_k}} < \delta_k      \\ \vspace{-.4 cm}  \\
 \delta_k   \quad  {\rm si } \quad  {{\theta_k}\over{\delta_k}} > \delta_k \, .
\end {array} \right. \monend
%

\bigskip \noindent
L'algorithme de brisure de ton est donc non lin\'eaire. En particulier, 
le nombre $\, T_{k+1} \,$ de nouveaux tons d\'epend du cas de figure. On a
$\, T_{k+1} = T_k + D_k \,$ si les anciens demi-tons~$\, \delta_k \,$ deviennent des
nouveaux tons $\,  \theta_{k+1} \,$ ou alors $\, T_{k+1} = T_k \,$ dans
le cas contraire o\`u l'on conserve les demi-tons. De mani\`ere analogue, le nombre
$\, D_{k+1} \,$ de nouveaux demi-tons vaut $\, D_{k+1} =  T_k\,$ dans le
premier cas et $\, D_{k+1} = T_k + D_k \,$ dans le second~:
\moneq \label {6.6}
T_{k+1} = \left \{ \begin {array}{l}
T_k + D_k  \quad    {\rm si } \quad    {\theta_k}  < {\delta_k}^2      \\
T_k  \qquad \quad \,\,\,  {\rm si } \quad   {\theta_k}  > {\delta_k}^2
\end {array} \right. \monend
%
\moneq \label {6.7}
D_{k + 1}  = \left \{ \begin {array}{l}
T_k  \qquad  \quad \,\,\,\,\,  {\rm si } \quad  {\theta_k}  < {\delta_k}^2   \\
T_k \,+\,  D_k   \quad  {\rm si } \quad   {\theta_k}  > {\delta_k}^2 \, .
\end {array} \right. \monend

\bigskip \noindent
La relation (\ref{6.1}) s'\'etend clairement \`a l'ordre $\, k \!+\! 1 \,$ ainsi que la
relation (\ref{6.2}) qui prend aussi la forme :
\moneq \label {6.7}
p_{k \!+\! 1}  =    p_k \,+\, T_k \,.
\monend
On peut alors recommencer \`a l'ordre suivant et poursuivre la mise en \oe uvre de l'algorithme.

\bigskip \noindent
Le nombre $\, N_k \,$ de gammes
$\,{\cal G}^k_\star \,$ avec $\, p_k \,$ notes et $\, T_k \,$ tons  est
simplement \'egal au nombre de fa\c{c}ons d'agencer $\, T_k \,$ tons et $\, D_k \,$
demi-tons, soit le nombre $\, N_k  \, $ de combinaisons :
\moneq \label{6.9}
N_k =  
\big( \,^{p_k} _{T_k} \big) =   \big( \,^{p_k} _{\! D_k} \big) =  {{p_k !}\over{T_k ! \, D_k !}} \,.
\monend

\smallskip \noindent
On peut  \'evoquer une ~\guillemotleft~{gamme-octave}~\guillemotright~
qui ne comporterait qu'une seule note : 
\moneqstar  
p_0 = 1 ,\quad  T_0 = 1 ,\quad    D_0 = 0 ,\quad   \theta_0 = 2  ,\quad
\delta_0 = 1  ,\quad N_0 =1 \,.
\monendstar 
Mais il est alors impossible de briser cette octave initiale par un demi-on \'evanescent ;  
la relation de r\'ecurrence d\'ebute avec les gammes  comportant   
deux notes et met en exergue la quinte : 
%
%
\moneq \label{6.10} 
p_1 = 2 ,\quad  T_1 = 1 ,\quad    D_1 = 1 ,\quad   \theta_1 = {3\over2} = 1.5 ,\quad
\delta_1 = {4\over3} \approx 1.333333 ,\quad N_1 =2 \,.
\monend
On  poursuit avec les gammes \`a trois notes
\moneqstar
p_2 = 3 ,\quad  T_2 = 2 ,\quad    D_2 = 1 ,\quad   \theta_2 = {4\over3} \approx 1.333333 ,\quad
\delta_2 = {9\over8} = 1.125 ,\quad N_2 = 3
\monendstar
et nous venons  de mettre en \'evidence au paragraphe pr\'ec\'edent toutes les gammes de troisi\`eme
g\'en\'eration \`a cinq notes~:
\moneqstar
p_3 = 5 ,\quad  T_3 = 2 ,\quad    D_3 = 3 ,\quad   \theta_3 = {{2^5}\over{3^3}} \approx  1.185185
,\quad \delta_3 = {9\over8}  = 1.125 ,\quad N_3 = 10 \,.
\monendstar
Les gammes de quatri\`eme g\'en\'eration comportent sept notes :
\moneqstar
p_4 = 7 ,\quad  T_4 = 5 ,\quad    D_4 = 2 ,\quad   \theta_4 = {{3^2}\over{2^3}}  = 1.125
,\quad \delta_4 = {{2^8}\over{3^5}}  \approx 1.053498  ,\quad N_4 = 21 \,
\monendstar
et contiennent toutes les gammes de la musique  classique occidentale, plus quelques autres
moins utilis\'ees. Nous y reviendrons.
Nous retenons que la construction r\'ecurrente des gammes avec une brisure
de ton, d\'efinie au d\'ebut de ce paragraphe par les relations (\ref{6.1}) \`a (\ref{6.9}),
permet de faire \'emerger  \`a la quatri\`eme \'etape une gamme de sept notes. 
Voil\`a selon nous pourquoi la gamme a sept notes !

\bigskip \noindent
Nous poursuivons  avec les  gammes
comprenant  
un nombre de notes plus important~:
\moneqstar
p_5 = 12 ,\quad  T_5 = 5 ,\quad    D_5 = 7 ,\quad   \theta_5 = {{3^7}\over
{2^{11}}}  \approx  1.067871 ,\quad \delta_5 = {{2^8}\over{3^5}}  \approx 1.053498 ,\quad N_5 = 792.
\monendstar
Le faible \'ecart entre le ton $ \,\, \theta_5  \approx 1.0678711 \,\,$ et le demi-ton
$\, \,\delta_5  \approx 1.053498 \,\, $ puisque\br
$\, {{\theta_5}\over{\delta_5}} = {{3^{12}}\over{2^{19}}} \approx 1.013643 \,$
soit 23 cents, 
pose la question de sa perception et
motive la gamme temp\'er\'ee
propo\-s\'ee d'abord par Simon Stevin (1548-1620) 
puis perfectionn\'ee par  Andreas Werckmeister (1645-1706)    
avant de se propager \`a  l'ensemble de la musique occidentale.
On les confond tous les deux au
b\'en\'efice de l'irrationnel racine douzi\`eme de 2.

\bigskip \noindent
Nous terminons ce paragraphe en montrant quelques gammes avec de plus en plus de
notes et poursuivons l'algorithme (de Pythagore ?). Nous omettons
d'expliciter ces diverses gammes,
compte tenu 
de leur nombre important,
sans outil plus \'elabor\'e de classi\-fication. Nous avons~:
\moneqstar \left \{ \begin {array}{l}
p_6 = 17 ,\quad T_6 = 12  ,\quad   D_6 = 5  ,\quad   \theta_6 =
{{2^8}\over{3^5}} \approx  1.053498 ,\quad  \delta_6 = {{3^{12}}\over{2^{19}}} \approx 1.013643 \\
 p_7 = 29 ,\quad  T_7 = 12 ,\quad   D_7 = 17 ,\quad    \theta_7 =
 {{2^{27}}\over{3^{17}}}  \approx 1.039318 ,\quad  \delta_7 = {{3^{12}}\over{2^{19}}} \approx  1.013643 \\
 p_8 = 41,\quad  T_8 = 12 ,\quad   D_8 = 29 ,\quad    \theta_8 =
 {{2^{46}}\over{3^{29} }}  \approx 1.025329  ,\quad  \delta_8 =
 {{3^{12}}\over{2^{19}}}  \approx  1.013643  \,.
\end {array} \right. \monendstar

\smallskip \noindent 
On retrouve l\`a une famille de  gammes de Pythagore \`a 41 degr\'es, dite de  Paul von Jank\'o (1901).
La famille suivante comporte 53 notes et
est connue  
avec le nom de  Nikolaus Mercator (1620-1687)  
(voir le musicologue William  Holder, 1694) : 
\moneqstar
p_9 = 53 ,\quad T_9 = 41  ,\quad   D_9 = 12  ,\quad   \theta_9 =
{{3^{12}}\over{2^{19}}}\approx  1.013643 ,\quad  \delta_9 = {{2^{65}}\over{3^{41}}} \approx 1.011529   \,.
\monendstar
Pour ces gammes, le ton et le demi-ton diff\`erent \`a nouveau tr\`es peu puisque\br
$\, {{\theta_9}\over{\delta_9}} = {{3^{53}}\over{2^{84}}} \approx 1.002090 $, soit 4 cents.
On a ainsi
$ \, | {{\theta_9}\over{\delta_9}} -1 | \approx 0.002  \ll  | \theta_9 -1 | \approx 0.014 \, $
et $ \,  | {{\theta_9}\over{\delta_9}} -1 |  \ll   | \delta_9 -1 | \approx 0.011 $.   
%
%
Il est alors concevable de confon\-dre les tons et les demi-tons pour former une
seule gamme de cinquante trois degr\'es.
Enfin, la dixi\`eme famille de gammes a les caract\'eristiques suivantes 
\moneqstar
p_{10} = 94  ,\quad T_{10} = 53  ,\quad   D_{10} = 41  ,\quad
\theta_{10} = {{2^{65}}\over{3^{41}}} \approx 1.011529
,\quad  \delta_{10} = {{3^{53}}\over{2^{84}}} \approx 1.002090    \,.
\monendstar
%

\bigskip \noindent
Si on \'etudie le comportement des diverses
familles de gammes pour~$ \, k \, $ tendant vers l'infini, on trouve un syst\`eme dynamique
discret qui a des ruptures lorsque
la diff\'erence entre le nombre irrationnel $\, k \, {{{\rm log}3}\over {{\rm log}\,2}} \,$ et
sa partie enti\`ere
saute d'une valeur proche de 1 \`a une valeur
proche de 0. 
Dans un de ses articles, Franck Jedrzejewski (2007) \'evoque  l'\'etude des gammes pythagori\-ciennes et introduit
une famille de gammes comportant successivement 
3, 5, 7, 12, 17, 29, 41, 53 et 94 notes. Il observe que cette suite de nombres  
a une parent\'e avec les valeurs des fractions r\'eduites du d\'eveloppement en fractions continues
(voir Michel Demazure, 1997)
du nombre $ \, \log_2(3/2) $. Pourtant, il n'y a pas identit\'e entre
la suite~$ \, p_k \, $ du nombre de notes des diverses familles et les
nombres entiers qui apparaissent dans les fractions  r\'eduites du d\'eveloppement
de $ \, \log_2(3/2) \, $  en fractions continues.
On a en effet
\moneqstar
\log_2 \Big({3\over2} \Big) \,  \simeq \, 1 \,,\,\, {1\over2}  \,,\,\, {3\over5}  \,,\,\, {7\over12}  \,,\,\, {24\over41}  \,,\,\,
{31\over53}  \,,\,\, {179\over306}  \,,\,\, {\it etc.}
\monendstar
Comprendre les propri\'et\'e fines de l'alg\`ebre des gammes de Pythagore reste donc pour nous
un champ  de recherche potentiel. Au  paragraphe suivant, nous mettons  en \'evidence quelques
propri\'et\'es abstraites.

\bigskip  \bigskip     \noindent {\bf \large    7) \quad    Quelques propri\'et\'es math\'ematiques}

\smallskip \noindent 
L'algorithme d\'ecrit au paragraphe pr\'ec\'edent  (relations (\ref{6.2}) \`a ({\ref{6.9})) joint
\`a la ~\guillemotleft~{condi\-tion initiale}~\guillemotright~  
(\ref{6.10}) des gammes \`a deux notes satisfait plusieurs
propri\'et\'es math\'ematiques. Nous les  \'enon\c{c}ons et d\'emontrons dans ce
paragraphe.

\bigskip  \noindent    {\bf Proposition 1. $\quad$  \'Equilibre harmonique }

\noindent
Le nombre de tons $\,T_k \,$ et le nombre de demi-tons $\,D_k \,$ sont, pour $k$
entier sup\'erieur ou \'egal \`a 1, premiers entre eux :
\moneq  \label{7.1}
\bigl( \, T_k ,\, D_k \, \bigr) = 1 \, .
\monend
On a une propri\'et\'e analogue pour le nombre de tons  $\,T_k \,$ et le nombre $\,
p_k = T_k + D_k \,$ de notes d'une gamme d'ordre $k$ :
\moneq  \label{7.2}
\bigl( \, T_k ,\, p_k \, \bigr) = 1 \, .
\monend

\bigskip \noindent
La preuve de la Proposition 1 et en particulier celle de la  relation (\ref{7.1}) s'effectue
par r\'ecurrence sur l'entier $k$, en utilisant
l'identit\'e de B\'ezout (voir par exemple Demazure, 1997)~: les entiers $ \, x \, $ et $\, y \, $ sont
premiers entre eux si et seulement si il existe deux entiers $ \, u \, $ et $\, v\, $
de sorte que $ \, u \, x + v \, y =  1 $.
%
%
On a dans un premier temps $\, T_1 = 1 \,$ et $\, D_1 = 1 ,\,$ donc ces deux unit\'es
sont clairement premi\`eres entre elles. Si la relation  (\ref{7.1}) est vraie jusqu'\`a
l'ordre~$ \, k$, on a
\moneq  \label{7.4}
u_k \,T_k +  v_k \, D_k =  1
\monend
pour deux entiers $ \, u_k \, $ et $ \, v_k $. Mais alors, $\, T_k \,$ et  $\, T_k + D_k \,$ sont
encore premiers entre eux, puisque la relation (\ref{7.4}) s'\'ecrit aussi :
\moneq  \label{7.5}
(u_k - v_k) \,T_k +  v_k \, (T_k + D_k ) = 1 \, .
\monend
Par suite, compte tenu des relations (\ref{6.6}) et (\ref{6.7}), les nombres de nouveaux
tons $\, T_{k+1} \,$ et de  nouveaux  demi-tons $\, D_{k+1} \,$ sont toujours premiers entre eux. De plus, compte
tenu de la relation (\ref{6.2}) et de l'identit\'e de B\'ezout,  la relation (\ref{7.5}) exprime
exactement la relation (\ref{7.2}) et la proposition est \'etablie \`a l'ordre suivant.
 \hfill $\square$

\bigskip  \noindent    {\bf Proposition 2. $\quad$  Factorisation   de la structure multiplicative}

\noindent
Pour $k$ entier sup\'erieur \`a 1  
arbitraire, le nombre $\, N_k \,$ de gammes \`a $
\,p_k \,$ notes (voir les relations (\ref{6.2}) et (\ref{6.9})) est divisible par $\, p_k \, $:
\moneqstar
\exists \,\,  \tau_k \in \N \,,\,\, N_k = p_k \, \tau_k \,.
\monendstar
%

\bigskip \noindent
Pour d\'emontrer la  Proposition 2, rappelons d'abord que  la relation (\ref{6.9}) s'\'ecrit\br
$ \,  N_k = \begin{pmatrix} p_k \cr T_k \end{pmatrix}  $. Donc elle  entra\^{\i}ne
%
%
$ \, T_k \, N_k = p_k \,  \begin{pmatrix} p_k - 1  \cr T_k - 1 \end{pmatrix}  \, $
et l'entier  $ \, p_k \, $ divise le produit~$\, T_k \, N_k .\, $ Comme $\, p_k \, $ et $\, T_k \, $ sont
premiers entre eux d'apr\`es la proposition pr\'ec\'edente, $\, p_k \, $ divise $\, N_k \, $
compte tenu du lemme de Gauss.
La proposition en r\'esulte.   \hfill $\square$

\bigskip \noindent
On rappelle que  la famille de rationnels $\, \xi_k \,$ d\'ej\`a utilis\'ee
pour d\'efinir la famille des notes $\, \Omega_k $
est d\'efinie \`a la relation (\ref{7.9}). 
Compte tenu de la relation (\ref{3.2}),
on a
\moneqstar
\Omega_k = \xi_k \,\,   \omega_0 \,,\,\,  k \in \ZZ \,
\monendstar
%
%
%
On remarque qu'on dispose \'egalement de la relation 
\moneqstar
\xi_k \,\, \xi_\ell = \xi_{k+\ell} \,\, {\rm ou} \,\, {1\over2} \, \xi_{k+\ell} 
\monendstar
selon que $\,  \xi_{k+\ell} \, $ est sup\'erieur ou inf\'erieur \`a 2. 
On a ensuite la propri\'et\'e suivante :

\bigskip  \noindent    {\bf Proposition 3. $\quad$  Notes d'une gamme et d\'ecalage d'indices   }

\noindent
Soit $\,  {\cal G}^k \,$ (pour $ \, k \geq 1$) l'ensemble des gammes $\,  {\cal G}^k_j \, (1 \leq j \leq N_k) \,$
d\'efinies aux relations (\ref{6.2})  \`a (\ref{6.9}) du  paragraphe 6
et $ \, \xi_m \, $  (pour $ \, m \in \ZZ $ )  l'ensemble des hauteurs relatives des notes
d\'efinie aux relations (\ref{3.3})  et  (\ref{7.9}).
Alors il existe un signe $\, \varepsilon_k \in \{-1 ,\, +1\} \,$ tel que
\moneq \label{7.10-7.11}
\theta_k = \xi\ib{\, \varepsilon_k \, D_k}  \,, \,\,\,
\delta_k = \xi\ib{ \,-\varepsilon_k \, T_k} \, .
\monend
%


\bigskip \noindent
Le ton $\, \theta_k \,$ 
correspond pour la famille $\, \Omega \,$ de toutes les notes \`a
un d\'ecalage de param\`etre  \'egal (au signe pr\`es) au nombre de demi-tons et le
demi-ton $\, \delta_k \,$ \`a un d\'ecalage \'egal (au signe pr\`es) au nombre de
tons. On peut r\'e-\'ecrire de cette fa\c{c}on quelques gammes vues lors des
paragraphes pr\'ec\'edents :
\moneq  \left \{ \begin {array}{rcl}  \label{7.12-a-7.14}
{\cal G}^1_1 &=& \bigl( \, \Omega_0 ,\, \Omega_{0+1} ,\, 2 \,
\Omega_{0+1-1} \, \bigr) \\
{\cal G}^2_3  &=& \bigl( \, \Omega_0 ,\, \Omega_{0-1} ,\,
\Omega_{0-1-1} ,\, 2 \, \Omega_{0-1-1+2} \, \bigr) \, \\
{\cal G}^3_{10} &=& \bigl( \, \Omega_0 ,\, \Omega_{0+2} ,\,
\Omega_{0+2+2} ,\, \Omega_{0+2+2+2} ,\, \Omega_{0+2+2+2-3} ,\,   2 \,
\Omega_{0+2+2+2-3-3} \, \bigr) \, . 
\end {array} \right. \monend
%
On dispose  ainsi  d'un tableau de correspondances  entre
les valeurs du ton et du demi-ton des diff\'erentes gammes et 
les num\'eros des notes. 
%
%

\smallskip 
\moneq \label{7.15-a-7.17}  \begin {array}{l}
\begin{tabular}{|c|c|c|c|}    \hline
famille de gammes num\'ero  $\, k $  & nombre de notes  $\, p_k $  &  ton $ \, \theta_k $ & demi-ton $ \, \delta_k $    \\   \hline
 1 &  2 & $ \xi_1 $ & $ \xi_{-1} $ \\ 
 2 &  3 & $ \xi_{-1} $ & $ \xi_{2} $ \\
 3 &  5 & $ \xi_{2} $ & $ \xi_{-3} $ \\ 
 4 &  7 & $ \xi_{2} $ & $ \xi_{-5} $ \\ 
 5 &  12 & $ \xi_{7} $ & $ \xi_{-5} $ \\ 
 6 &  17 & $ \xi_{5} $ & $ \xi_{12} $ \\ 
 7 &  29 & $ \xi_{-17} $ & $ \xi_{12} $ \\ 
 8 &  41 & $ \xi_{-29} $ & $ \xi_{12} $ \\ 
 9 &  53 & $ \xi_{12} $ & $ \xi_{-41} $ \\ 
 10 & 94 &  $ \xi_{-41} $ & $ \xi_{53} $ \\ 
\hline \end{tabular}
\end {array} \monend

\bigskip \noindent
La Proposition 3 s'\'etablit par
r\'ecurrence sur $ \, k$.
Comme la preuve est assez technique, nous l'avons renvoy\'ee en Annexe B.

\bigskip  \bigskip    \noindent {\bf \large    8) \quad   Structure multiplicative et tonalit\'e }

\smallskip \noindent 
On utilise dans ce paragraphe une gamme $\, {\gamma} \,$
de la famille \'etudi\'ee ici, \`a $\, p \, $ notes et qui commence
par un {\textsl {do}}. C'est une suite $\, \bigl( \nu_j \bigr)\ib{0 \leq j \leq p} \,$
qui satisfait \`a la relation (\ref{2.2}) : 
\moneqstar 
\omega_0 = \nu_0 ,\, \nu_1 ,\, \cdots ,\, \nu_j <
\nu_{j+1} ,\, \cdots ,\, \nu_{p-1} ,\, \nu_p = 2 \, \omega_0 
\monendstar 
De plus, le rapport $\, \nu_{j+1} / \nu_j \,$ de deux notes successives est soit un
ton $\, \theta \,$ soit un demi-ton~$\, \delta \,$:
\moneqstar
{{\nu_{j+1}} \over {\nu_j}} \in \{ \, \theta ,\, \delta \, \} , \qquad 0 \leq j
\leq p-1 \,.
\monendstar
Par d\'efinition, la ~\guillemotleft~{structure multiplicative}~\guillemotright~
 $\, {\cal S}(\gamma) \,$ de la
gamme $\gamma$ est la suite des rapports de fr\'equences entre deux notes
successives,   c'est-\`a-dire~:
\moneqstar
 {\cal S}(\gamma) = \Bigl( \, {{\nu_{1}} \over {\nu_0}} ,\,
 {{\nu_{2}} \over {\nu_1}} ,\, \cdots ,\, {{\nu_{j+1}} \over {\nu_j}} ,\,
\cdots ,\, {{\nu_{p}} \over {\nu_{p-1}}} \, \Bigr) \,\,  \in  \,\,  \{ \, \theta ,\, \delta \, \}^p \,.
\monendstar
La structure multiplicative $\,  {\cal S}(\gamma)  \,$ permet de savoir quel est l'agencement
pr\'ecis des tons et des demi-tons de la gamme $\, \gamma .\,$ Si $T$ est le nombre
de tons et $D$ le nombre de demi-tons de cette gamme, on a vu que
\moneqstar
p = T \,+\, D \, .
\monendstar
De plus (Proposition 3), il existe $\, \varepsilon \in \{-1 ,\, +1 \} \,$ tel
que
\moneqstar
\theta = \xi\ib{ \,\varepsilon \, D} , \qquad \delta =
\xi\ib{ -\varepsilon \, T}
\monendstar
ce qui permet de construire les notes $\, \nu_j \,$ de la gamme $\, \gamma \,$ par
d\'ecalage successif des indices des notes ~\guillemotleft~de base~\guillemotright~
$\, \Omega_k $, comme illustr\'e aux exemples (\ref{7.12-a-7.14}). 
Il est donc naturel de ne pas forc\'ement
commencer une gamme de structure  $\,  {\cal S}(\gamma)  \,$ par le {\textsl {do}}
(ou $\, \Omega_0 $) mais par une note~$\, \Omega_k\,$ arbitraire telle qu'introduite au
paragraphe~3.

\bigskip \noindent
On se donne une structure  multiplicative $ \, {\cal S} $. 
La gamme $\,  {\cal G} ({\cal S}) (\Omega_k) \,$ de
structure  multiplicative~$\,  {\cal S} \,$ et  de  ~\guillemotleft~{tona\-lit\'e}~\guillemotright~$\,\Omega_k \,$
est par d\'efinition la gamme obtenue \`a partir de la relation (\ref{2.2}) en gardant la m\^eme
structure
\moneqstar
    {\cal S} = \bigl( \mu_1 =   {{\nu_{1}} \over {\nu_0}} ,\, \mu_2 =  {{\nu_{2}} \over
{\nu_1}} ,\, \cdots ,\,  \mu_j =  {{\nu_{j+1}} \over {\nu_j}} ,\, \cdots ,\, \mu_p
    =   {{\nu_{p}} \over {\nu_{p-1}}} \bigr) 
\monendstar
que la gamme $\, \gamma \,$  mais en
changeant la premi\`ere note g\'en\'erique  $\, \omega_0 \,$ pour la note
particuli\`ere $\, \Omega_k .\,$ On a donc en g\'en\'eral
\moneq \label{8.6}
{\cal G} ({\cal S}) (\Omega_k)  = \bigl( \, \Omega_k ,\, \mu_1 \, \Omega_k ,\,
\mu_2 \, \mu_1 \, \Omega_k ,\, \cdots ,\, (\mu_j \mu_{j-1} \cdots \mu_1 )
 \,  \Omega_k ,\, \cdots ,\, 2 \, \Omega_k  \, \bigr) \,. 
 \monend
%
%
Nous avons vu au paragraphe 3 que les structures des deux gammes primitives \`a deux
notes, \`a savoir
\moneqstar
{\cal S}_1 = \Bigl( {3\over2} ,\, {4\over3} \Bigr) \,,\,\,
{\cal S}_2 = \Bigl( {4\over3}  ,\, {3\over2} \Bigr)  \,
\monendstar
permettent de construire de proche en proche la suite infinie $\, \bigl( \Omega_k
\bigr)\ib{k \in \ZZ} \,$ de toutes les notes.

\bigskip  \noindent 
Une fois la structure multiplicative fix\'ee, la transposition avec comme premi\`ere
note $ \, \Omega_k \, $ permet de faire varier la hauteur absolue, mais pas les intervalles
entre les notes d'une m\^eme gamme. Nous l'avons explicitement utilis\'ee avec les gammes
\`a cinq notes. 
Afin de s'autoriser une richesse tonale
arbitraire, nous ne devons retenir de la construction propos\'ee au paragraphe 6 que
la structure multiplicative des diff\'erentes gammes. Le paragraphe suivant propose
d'y mettre un peu d'ordre, avec une classification en
~\guillemotleft~{type}~\guillemotright~ et 
~\guillemotleft~{mode}~\guillemotright.

\bigskip   \bigskip   \noindent {\bf \large    9) \quad  Type et mode}

\smallskip \noindent 
On suppose dans ce paragraphe que l'on a fix\'e le nombre $\,p \,$ de notes (avec $p$
de la forme $\,p=p_k \,$ pour un certain entier $\,k \,$ que l'on n'\'ecrit pas pour
all\'eger les notations),  un nombre de tons $\,T \,$ et un nombre de demi-tons $\, D $,
la valeur du ton $\, \theta \,$ et du demi-ton $\, \delta .\,$ On s'int\'eresse \`a 
l'ensemble des gammes commen\c{c}ant
par la note {\textsl {do}} et  comportant exactement $ \, p \, $ notes.
L'ensemble  $ \, {\rm {\bf S}}^{^{\scriptstyle p}} \, $
des structures multiplicatives des gammes
\`a $ \, p \, $ notes est l'ensemble des $\,p$-uplets  \`a valeurs dans la paire
$\, \{ \theta ,\, \delta \} \,$ tels que le nombre de tons $\, \theta \,$ est \'egal \`a
$T$ et le nombre de demi-tons $\, \delta \,$ est exactement \'egal \`a $ \, D \,$  :
\moneqstar
{\rm {\bf S}}^{^{\scriptstyle p}} = \big\{ (\mu_1 ,\, \mu_2 ,\, \cdots ,\, \mu_p) ,\,
\mu_j \in \{\theta,\, \delta\} \,\, {\rm si} \,\, 1 \leq j \leq p ,\,
\sharp \{j,\, \mu_j = \theta \} = T ,\, \sharp \{j,\, \mu_j = \delta \} = D \big\} \,. 
\monendstar 
Il est facile de voir que 
l'ensemble  $ \, {\rm {\bf S}}^{^{\scriptstyle p}} \, $ comporte  exactement
$\, N = \big( ^{\,p} _{\, T} \big) \,$ \'el\'ements
distincts,
comme pour l'ensemble des gammes \`a $ \, p \, $ notes, 
puisqu'il suffit de choisir la position des   $ \, T \, $ tons parmi les $ \, p \, $
possibles rapports des notes successives.

\bigskip \noindent
On regarde maintenant comment agit la permutation circulaire $\, \sigma \,$
d\'efinie par 
la transformation 
$ \,\, 1 \longmapsto~2 ,\,  \, 2 \longmapsto 3 ,\,  ... ,\, p \!-\!1 \longmapsto p \,\, {\rm et}
\,\,  p \longmapsto 1 \,\, $ 
%
de l'ensemble $\, \{1,\, \cdots ,\, p\} \,$ :
\moneqstar
\sigma = \begin{pmatrix}
  1 & 2 & $\dots$ & p-1 & p  \\ 2 & 3 &  $\dots$ & p & 1
\end{pmatrix}
\monendstar
sur l'ensemble  $ \, {\rm {\bf S}}^{^{\scriptstyle p}}  \, $
des structures multiplicatives des gammes \`a $ \, p \, $  notes.   
\`A  partir d'une
configuration donn\'ee $\,  (  \mu_1 ,\, \mu_2 \, ,\, \cdots \, ,\, \mu_p  ) \,$
appartenant \`a  $ \, {\rm {\bf S}}^{^{\scriptstyle p}} , \, $ on pose
\moneqstar
\sigma \,  (  \mu_1 ,\, \mu_2 \, ,\, \cdots \, ,\, \mu_p  )  = \bigl( \,
\mu\ib{\sigma(1)} ,\, \mu\ib{\sigma(2)} \, ,\, \cdots \, ,\, \mu\ib{\sigma(p)} \,
\bigr)
\monendstar
qui appartient encore \`a la famille  $ \, {\rm {\bf S}}^{^{\scriptstyle p}} . \, $
Deux structures multiplicatives  $\,  (  \mu_1 ,\, \mu_2  ,\, \cdots \, \,$ $ \,
\cdots \, ,\, \mu_p  ) \,$ et $\,  (  \mu'_1  ,\, \mu'_2 \, ,\, \cdots \, ,\,
\mu'_p  ) \,$ sont dites  ~\guillemotleft~{de m\^eme  type}~\guillemotright~
si il existe une permutation circulaire
$\, \sigma \,$ de l'ensemble $\, \{1,\cdots ,\, p\} \,$ et un entier $\, i \in
\{0,\, \cdots ,\, p\!-\!1 \} \,$  de sorte que
\moneq \label{9.4}
\mu'_j = \mu\ib{ \, \sigma^{\scriptstyle i} (j)}  \,\,\, {\rm pour} \,\,
{\rm tout } \,\, j \,\, {\rm tel} \, {\rm que}  \,\,\,  1 \leq j \leq p \,.
\monend
%

\bigskip \noindent
L'ensemble  $ \, {\rm {\bf T}}^{^{\scriptstyle p}}  \, $ des types des gammes \`a $p$
notes est l'ensemble des structures multipli\-catives des gammes \`a $ \, p \, $ notes \`a une
permutation circulaire pr\`es. Le nombre de types
\moneqstar
\sharp \, {\rm {\bf T}}^{^{\scriptstyle p}}   = {{1}\over{p}} \,
\begin{pmatrix} p \cr T  \end{pmatrix} \,\,; \quad p = T + D \,
\monendstar
est bien un nombre entier, comme on l'a vu lors de la proposition 2. 
Il y a un seul type de gammes \`a deux  notes et un seul type de gammes \`a
trois notes, 2 types de gammes \`a cinq notes, 3 types de gammes \`a sept
notes, 66 types de gammes \`a douze notes, {\it etc.}  D'un point de vue math\'ematique,
l'ensemble des types $\,  {\rm {\bf T}}^{^{\scriptstyle p}}  \, $ est le quotient de
l'ensemble~$\,  {\rm {\bf S}}^{^{\scriptstyle p}}  \, $ par la relation d'\'equivalence $\, \sim \,$
suivante :
\moneqstar
\mu \, \sim \, \mu' \,\,\, {\rm  si} \,\, {\rm  et} \,\,{\rm  seulement} \,\, {\rm  si} \,\,\, \exists \, i \in \{0,\, \cdots
,\,p\!-\!1 \} ,\, {\rm la} \,\, {\rm  relation} \,\,  (\ref{9.4}) \,\, {\rm a} \,\,  {\rm lieu}
\monendstar
et on note $ \,\, {\rm {\bf T}}^{^{\scriptstyle p}}  =  {\rm {\bf S}}^{^{\scriptstyle p}}  / \sim $.
%
D'un point de vue pratique, le type d\'efinit la succession pr\'ecise des tons et des demi-tons, mais
\`a une permutation circulaire pr\`es. 

\bigskip \noindent
On suppose dans la suite qu'un repr\'esentant $\, \tau \in \widetilde{\tau}  \,$ de
chaque type $\, \widetilde{\tau} \in \, {\rm {\bf T}}^{^{\scriptstyle p}}  \, $ a
\'et\'e choisi. Il s'agit d'une structure multiplicative, \`a une
permutation pr\`es :
\moneq \label{9.8}
\tau = \bigl( \, \mu_1 ,\, \mu_2 \, ,\, \cdots ,\,
 \mu_j ,\, \cdots \, ,\, \mu_p \,  \bigr) \, \in  \, \widetilde{\tau} \, \in
\, {\rm {\bf T}}^{^{\scriptstyle p}}  \,.
\monend
Pour les gammes de moins de cinq notes, nous proposons le choix suivant :
\moneqstar \left \{ \begin {array}{rclcl}
\tau^1 &=& \bigl( \, {3\over2} ,\, {4\over3} \, \bigr) &\in& {\rm \bf S}^{^{\scriptstyle 2}} \, \\
\tau^2 &=& \bigl( \, {4\over3} ,\, {9\over8} ,\, {4\over3} \, \bigr)  &\in& {\rm \bf S}^{^{\scriptstyle 3}} \, \\
\tau_1^3 &=& \bigl( \, {{32}\over{27}} ,\, {9\over8} ,\,{9\over8} ,\,
{9\over8} ,\, {{32}\over{27}}\, \bigr) &\in&   {\rm \bf S}^{^{\scriptstyle 5}} \, \\
\tau_2^3 &=& \bigl( \, {9\over8} ,\, {{32}\over{27}} ,\, {9\over8} ,\,
{{32}\over{27}} ,\, {9\over8}  \, \bigr) &\in&    {\rm \bf S}^{^{\scriptstyle 5}} \,.
\end {array} \right.
\monendstar
On a alors, compte tenu des relations (\ref{2.6}) (\ref{4.12})
et de la structure des gammes $ \bigl( {\cal G}^3_3 \bigr) \, $ et $  \bigl( {\cal G}^3_5 \bigr) \, $:
%
\moneq \label{9.13-a-9.16}  
{\cal S}\bigl({\cal G}_1^1 \bigr) = \tau^1  \,,\,\, 
{\cal S}\bigl({\cal G}_2^2 \bigr) = \tau^2 \,,\,\,
{\cal S}\bigl({\cal G}_5^3 \bigr) =\tau_1^3 \,,\,\,
{\cal S}\bigl({\cal G}_3^3 \bigr) =\tau_2^3 \,.
\monend

\bigskip \noindent
Une fois le type fix\'e et (surtout) une structure multiplicative $\, \tau \,$
qui le repr\'esente, on construit simplement une 
~\guillemotleft~gamme fondamentale~\guillemotright~
$\, {\cal G}_{\tau} \,$ en commen\c{c}ant par convention par la note {\textsl {do}}, comme on l'a
vu aux relations (\ref{9.13-a-9.16})  
pour les gammes \`a deux, trois et cinq notes. Il
est alors naturel de fabriquer les autres gammes de m\^eme type, mais par
d\'efinition de   ~\guillemotleft~{mode}~\guillemotright~ diff\'erent, en construisant les $p$ gammes ayant les
m\^emes notes que $\, {\cal G}_{\tau} ,\,$ mais obtenues par action successive de la
permutation circulaire $\, \sigma $.  Si le type $\, \tau \,$ a un repr\'esentant
donn\'e \`a la relation  (\ref{9.8}), on a (cf. (\ref{8.6})) :
\moneqstar 
{\cal G}_{\tau} =  \bigl( \, \omega_0 ,\, \mu_1 \, \omega_0 ,\, \cdots
,\, (\mu_1  \cdots \mu_j ) \, \omega_0 ,\, \cdots ,\, (\mu_1  \cdots \mu_{p-1} )
\, \omega_0 ,\, 2 \, \omega_0 \,\bigr) \, . 
\monendstar
%
Nous notons $\, \nu(\tau,\,i) \, $ les notes $\, p \, $ d'une gamme $ \, {\cal G}_{\tau} \, $ de type $ \, \tau \, $
avec $ \, 0 \leq i < p \, $ et par convention
$ \, \nu(\tau,\,0) = \Omega_0 = {\textsl {do}} $.
De plus, si $\, \tau = ( \mu_1 ,\, \mu_2 \, ,\, \cdots ,\,  \mu_j ,\, \cdots \, ,\, \mu_p ) $,
on a les relations 
\moneqstar
\nu(\tau,\,i+1) = \mu_i \,\,  \nu(\tau,\,i) 
\monendstar
pour  $ \, 1 \leq i < p-1 $.
%
Par exemple, pour la gamme \`a une seule note $ \, {\cal G}^1_1 = ( {\textsl {do}} ,\, {\textsl {sol}} ,\,  {\textsl {do}}^* ) \, $ de type $ \, \tau^1 $,
on a simplement 
\moneqstar 
\nu(\tau^1,\,0) = {\textsl {do}} ,\,\, 
\nu(\tau^1,\,1) = {\textsl {sol}} \, . 
\monendstar
%
Pour une gamme \`a trois notes, 
$ \, {\cal G}^2_1 = ( {\textsl {do}} ,\, {\textsl {fa}} ,\, {\textsl {sol}} ,\,  {\textsl {do}}^* ) \, $ de type $ \, \tau^2 $,
on a 
\moneq \label{17aout}
\nu(\tau^2,\,0) = {\textsl {do}} ,\,\, 
\nu(\tau^2,\,1) = {\textsl {fa}} ,\,\, 
\nu(\tau^2,\,2) = {\textsl {sol}} \, . 
\monend 
Pour les gammes \`a 5 notes, on dispose d'abord du type
$ \, \tau^3_1 = ( \theta_3 ,\, \delta_3 ,\, \delta_3 ,\,   \delta_3 ,\, \theta_3 ) \, $
comme pour la gamme
$ \, {\cal G}^3_5 = ( {\textsl {do}} ,\,   {\textsl {mi}}\,\flat ,\,  {\textsl {fa}} ,\, {\textsl {sol}} ,\, {\textsl {la}} ,\,  {\textsl {do}}^* ) \, $
et on en d\'eduit simplement
\moneq \label{9.22}
\nu(\tau^3_1,\,0) = {\textsl {do}} ,\,\,
\nu(\tau^3_1,\,1) =  {\textsl {mi}}\,\flat  ,\,  
\nu(\tau^3_1,\,2) = {\textsl {fa}}  ,\,\, 
\nu(\tau^3_1,\,3) = {\textsl {sol}}  ,\,   
\nu(\tau^3_1,\,4) = {\textsl {la}} \, .
\monend
Pour le type 
$ \, \tau^3_2 = ( \delta_3 ,\,  \theta_3 ,\, \delta_3 ,\, \theta_3 ,\,  \delta_3 ) $,
le repr\'esentant naturel est 
$ \, {\cal G}^3_3 = ( {\textsl {do}} ,\,   {\textsl {r\'e}} ,\,  {\textsl {fa}} ,\, {\textsl {sol}} ,\,  {\textsl {si}}\,\flat  ,\, {\textsl {do}}^* ) \, $
et
\moneq \label{9.23}
\nu(\tau^3_2,\,0) = {\textsl {do}}  ,\,\,
\nu(\tau^3_2,\,1) =  {\textsl {r\'e}}  ,\,\,  
\nu(\tau^3_2,\,2) = {\textsl {fa}}  ,\,\,  
\nu(\tau^3_2,\,3) = {\textsl {sol}}  ,\,\,  
\nu(\tau^3_2,\,4) =  {\textsl {si}}\,\flat \, .
\monend

\bigskip \noindent
Avec les notes $ \, \nu (\tau,\, i) $, nous formons une nouvelle gamme simplement en changeant la premi\`ere
note et en laissant agir la permutation circulaire. Nous d\'efinissons de cette fa\c con un nouveau mode.
Par exemple, pour une gamme \`a deux notes, on obtient \`a partir de 
$ \, \nu(\tau^1,\,1) = {\textsl {sol}} $, la gamme $ \, (  {\textsl {sol}} ,\,  {\textsl {do}}^* ,\, {\textsl {sol}}^* ) $.
Ce ~\guillemotleft~mode de sol~\guillemotright~ d\'efinit la gamme
$ \, {\cal G}(\tau^1,\,1) \, $ qui est identique \`a la gamme
$ \, {\cal G}^1_2(sol) = {\cal G}^1_2(\Omega_1) \, $ obtenue \`a partir
de la gamme $ \, {\cal G}^1_2 =  (  {\textsl {do}} ,\,  {\textsl {fa}} ,\, {\textsl {do}}^* ) $.

\bigskip \noindent
Le ~\guillemotleft~mode~\guillemotright~  $\, i \in \{0,\cdots,\, p\!-\!1 \} \,$ du type $\, \tau \,$ est engendr\'e
par la gamme de premi\`ere note~$\, \nu(\tau,\, i) \,$ et de structure
$\, \bigl( \mu\ib{\sigma^{\scriptstyle i}(1)} ,\, \mu\ib{\sigma^{\scriptstyle i}(2)} \, ,\,
\cdots \, ,\, \mu\ib{\sigma^{\scriptstyle i}(p)}  \bigr) \,$ : on fait agir $i$
fois la permutation circulaire $\, \sigma \,$ de r\'ef\'erence. On pose en
cons\'equence :
\moneqstar
{\cal G}(\tau,\,i) = \bigl( \,  \nu(\tau,\, i)  ,\,
\nu(\tau,\, i+1) ,\, \cdots ,\, \nu(\tau,\, p-1)   ,\,
2 \, \Omega_0 ,\, 2 \, \nu(\tau,\, 1),\, \cdots ,\, 2 \, \nu(\tau,\, i-1)   \, \bigr)
\monendstar
et on a par d\'efinition m\^eme
\moneqstar
{\cal S}({\cal G}(\tau,\,i)) = \bigl( \, \mu\ib{\sigma^{\scriptstyle
i}(1)} ,\, \mu\ib{\sigma^{\scriptstyle i}(2)} \, ,\, \cdots \, ,\,
\mu\ib{\sigma^{\scriptstyle i}(p)} \, \bigr) \,.
\monendstar
%

\bigskip \noindent
Il peut \^etre utile de nommer le mode,  c'est-\`a-dire  le choix pr\'ecis de la
structure multiplicative  \`a travers le type et l'ordre de la permutation circulaire
gr\^ace au  nom de la premi\`ere note $\, \nu(\tau,\, i) $.  Pour les gammes
\`a deux notes, on a un seul type et deux modes, le mode de {\textsl {do}} qui d\'ebute par
une quinte et le mode de {\textsl {sol}} qui commence par une quarte :
\moneqstar \left \{ \begin {array}{l}
{\cal G}(\tau^1,\,0) = {\cal G}(\tau^1,\,{\textsl {do}}) = \bigl( \, {\textsl {do}} ,\, {\textsl {sol}} ,\, 
    {\textsl {do}}^{*} \, \bigr) = {\cal G}^1_1 \\
{\cal G}(\tau^1,\,1) = {\cal G}(\tau^1,\,{\textsl {sol}}) = \bigl( \, {\textsl {sol}} ,\, {\textsl {do}} ,\, 
{\textsl {sol}}^{*} \, \bigr) = {\cal G}^1_2({\textsl {sol}})  \,.
\end {array} \right.
\monendstar
Pour les gammes \`a trois notes, on a encore un seul type, donn\'e \`a la relation (\ref{17aout}) 
mais, outre le mode de {\textsl {do}}, on dispose des modes de {\textsl {fa}} et de {\textsl {sol}} :
\moneqstar \left \{ \begin {array}{l}
{\cal G}(\tau^2,\,0) =  {\cal G}(\tau^1,\,{\textsl {do}}) = \bigl( \, {\textsl {do}} ,\, {\textsl {fa}} ,\, {\textsl {sol}} 
,\, {\textsl {do}}^{*} \, \bigr) = {\cal G}_2^2  \\
{\cal G}(\tau^2,\,1) = {\cal G}(\tau^2,\,{\textsl {fa}}) = 
\bigl( \, {\textsl {fa}} ,\, {\textsl {sol}} ,\,
{\textsl {do}}^{*} ,\, {\textsl {fa}}^{*} \, \bigr) = {\cal G}_2^1(\Omega_{-1})   \\
{\cal G}(\tau^2,\,2) = {\cal G}(\tau^2,\,{\textsl {sol}}) =
\bigl( \,  {\textsl {sol}} ,\, {\textsl {do}}^{*} ,\,
{\textsl {fa}}^{*} ,\, {\textsl {sol}}^{*} \,  \bigr) = {\cal G}_2^3(\Omega_{1})  \,.
\end {array} \right.
\monendstar
Pour les gammes \`a cinq notes, la classification (type, mode) permet un nouveau
rangement des dix gammes vues au paragraphe 5 comme des gammes qui d\'ebutent par la
note {\textsl{do}}. On trouve 
pour le premier type  (relation  (\ref{9.22})), les cinq
modes de {\textsl {do}}, $\, {\textsl {mi}}\,\flat$, {\textsl {fa}}, {\textsl {sol}}, {\textsl {la}} :
\moneqstar \left \{ \begin {array}{l}
{\cal G}(\tau^3_1,\,0) =  {\cal G}(\tau^1,\,{\textsl {do}})  
  = \bigl( \, {\textsl {do}} ,\,  {\textsl {mi}}\,\flat  ,\,   {\textsl {fa}} ,\,    {\textsl {sol}} ,\,
 {\textsl {la}} ,\,    {\textsl {do}}^{*} \, \bigr)   =  {\cal G}_5^3 \\
{\cal G}(\tau^3_1,\,1) = {\cal G}(\tau^3_1 ,\,  {\textsl {mi}}\,\flat) = \bigl( \,
 {\textsl {mi}}\,\flat  ,\,  {\textsl {fa}} ,\,   {\textsl {sol}} ,\,   {\textsl {la}} ,\,
 {\textsl {do}}^{*} ,\,     {\textsl {mi}}\,\flat^{*}  \, \bigr)   =  {\cal G}_{10}^3 (\Omega_{-3}) \\
{\cal G}(\tau^3_1,\,2) = {\cal G}(\tau^3_1 ,\,  {\textsl {fa}}) = \bigl( \, 
{\textsl {fa}} ,\,    {\textsl {sol}} ,\,   {\textsl {la}} ,\,    {\textsl {do}}^{*} ,\,
{\textsl {mi}}\,\flat^{*}  ,\,    {\textsl {fa}}^{*}  \, \bigr)
=   {\cal G}_{4}^3 (\Omega_{-1})  \\
{\cal G}(\tau^3_1,\,3) = {\cal G}(\tau^3_1 ,\,  {\textsl {sol}}) = \bigl( \, 
{\textsl {sol}} ,\,   {\textsl {la}} ,\,
{\textsl {do}}^{*} ,\,    {\textsl {mi}}\,\flat^{*}  ,\,    {\textsl {fa}}^{*}   ,\,  {\textsl {sol}}^{*}
\, \bigr)    =   {\cal G}_{7}^3 (\Omega_{1}) \\
{\cal G}(\tau^3_1,\,4) = {\cal G}(\tau^3_1 ,\,  {\textsl {la}}) = \bigl( \, 
 {\textsl {la}} ,\,     {\textsl {do}}^{*} ,\,     {\textsl {mi}}\,\flat^{*}  ,\,    {\textsl {fa}}^{*} 
 ,\,    {\textsl {sol}}^{*} ,\,     {\textsl {la}}^{*} \, \bigr)     \,\,=
\,\,    {\cal G}_{9}^3 (\Omega_{3}) \,.
\end {array} \right.
\monendstar
Pour le second type (relation (\ref{9.23})), on dispose des cinq modes de {\textsl {do}}, {\textsl r\'e},
{\textsl {fa}}, {\textsl {sol}} et ${\textsl {si}}\,\flat \, : $
\moneqstar \left \{ \begin {array}{l}
{\cal G}(\tau^3_2,\,0)   =  {\cal G}(\tau^1,\,{\textsl {do}}) = \bigl( \,
{\textsl {do}} ,\, {\textsl {r\'e}}  ,\,  {\textsl {fa}}  ,\, {\textsl {sol}} ,\, {\textsl {si}}\,\flat  ,\,   {\textsl {do}}^{*} \, \bigr)
=  {\cal  G}_3^3 \\
{\cal G}(\tau^3_2,\,1) = {\cal G}(\tau^3_1 ,\,  {\textsl {r\'e}}) = \bigl( \,
{\textsl {r\'e}}  ,\,  {\textsl {fa}}  ,\, {\textsl {sol}}  ,\,
{\textsl {si}}\,\flat  ,\,   {\textsl {do}}^{*} ,\, {\textsl {r\'e}}^{*} \, \bigr)   =  {\cal G}_8^3 (\Omega_2) \\
{\cal G}(\tau^3_2,\,2) = {\cal G}(\tau^3_1 ,\,  {\textsl {fa}}) = \bigl( \,
{\textsl {fa}}  ,\, {\textsl {sol}}  ,\, {\textsl {si}}\,\flat  ,\,   {\textsl {do}}^{*} ,\,
{\textsl {r\'e}}^{*} ,\,  {\textsl {fa}}^{*}   \, \bigr)   =  {\cal G}_2^3 (\Omega_{-1}) \\
{\cal G}(\tau^3_2,\,3) = {\cal G}(\tau^3_1 ,\,  {\textsl {sol}}) = \bigl( \,
{\textsl {sol}}  ,\, {\textsl {si}}\,\flat  ,\,   {\textsl {do}}^{*} ,\, {\textsl {r\'e}}^{*}
,\,  {\textsl {fa}}^{*}   ,\, {\textsl {sol}}^{*}  \, \bigr)   =  {\cal G}_6^3 (\Omega_{1}) \\
{\cal G}(\tau^3_2,\,4) = {\cal G}(\tau^3_2 ,\,  {\textsl {si}}\,\flat) = \bigl( \,
{\textsl {si}}\,\flat  ,\,   {\textsl {do}}^{*} ,\, {\textsl {r\'e}}^{*} ,\,  {\textsl {fa}}^{*}  
  ,\, {\textsl {sol}}^{*}  ,\, {\textsl {si}}\,\flat^{*}  \, \bigr)   =   {\cal G}_1^3 (\Omega_{-2})  \,.
\end {array} \right.
\monendstar
On retrouve en particulier des transpositions de la gamme pentatonique majeure $\, {\cal G}_1^3 \, $
et de la  gamme pentatonique mineure $\, {\cal G}_6^3  $.

\bigskip \noindent
De fa\c con g\'en\'erale, la structure form\'ee avec le  couple (type, mode) = $ \, (\tau ,\, i) \, $  permet de classifier
les gammes d'une famille $ \, k \, $ donn\'ee
et qui commencent par la note {\textsl {do}}. On peut noter $ \, {\cal G}(\tau ,\,i) \, $
une gamme de cette famille ou \'eventuellement
$ \, {\cal G}(\tau ,\,{\textsl {note}}) \, $ si $ \, \nu (\tau,\, i) = {\textsl {note}} $. 
Ensuite le choix d'une  tonalit\'e permet de faire commencer la gamme par une
note arbitraire $ \, \Omega_\ell \, $ et on obtient ainsi la gamme
$ \, {\cal G}(\tau ,\,i) (\Omega_\ell) $.
Par exemple,
$ \, {\cal G}(\tau^3_2 ,\,4)  = ( \Omega_{-2} ,\,  \Omega_{0}^* ,\,  \Omega_{2}^* ,\,  \Omega_{-1}^* ,\,  \Omega_{1}^* ,\,  \Omega_{-2}^* ) \, $
et
\moneqstar
{\cal G}(\tau^3_2 ,\,4)(\Omega_0) = \big( \Omega_{0} ,\,  \Omega_{2} ,\,  \Omega_{4} ,\,  \Omega_{1} ,\,  \Omega_{3} ,\,  \Omega_{0}^* \big)
=  \big(  {\textsl {do}} ,\,  {\textsl {r\'e}} ,\,  {\textsl {mi}} ,\,  {\textsl {sol}}  ,\,  {\textsl {la}}  ,\,  {\textsl {do}}^* \big)
= {\cal G}^3_1 \, . 
\monendstar

\smallskip \noindent 
On a aussi 
$ \, {\cal G}(\tau^3_2 ,\,2)  = ( \Omega_{-1} ,\,  \Omega_{1} ,\,  \Omega_{-2} ,\,  \Omega_{0}^* ,\,  \Omega_{2}^* ,\,  \Omega_{-1}^* ) \, $
et on en d\'eduit
\moneqstar
{\cal G}(\tau^3_2 ,\,2)(\Omega_0) = \big( \Omega_{0} ,\,  \Omega_{2} ,\,  \Omega_{-1} ,\,  \Omega_{1} ,\,  \Omega_{3} ,\,  \Omega_{0}^* \big)
=  \big(  {\textsl {do}} ,\,  {\textsl {r\'e}} ,\,  {\textsl {fa}} ,\,  {\textsl {sol}}  ,\,  {\textsl {la}}  ,\,  {\textsl {do}}^* \big)
= {\cal G}^3_2 \, . 
\monendstar
%

\noindent 
Ayant maintenant une structure de la forme (type, mode, tonalit\'e) pour les diff\'erentes gammes, 
nous pouvons aborder l'\'etude d\'etaill\'ee des gammes \`a sept notes. 

\bigskip   \bigskip   \noindent {\bf \large    10) \quad  Gammes \`a sept notes }

\smallskip \noindent 
Avec les d\'efinitions introduites plus haut, il existe 21 gammes \`a 7 notes qui
commencent par la note {\textsl {do}}, compos\'ees de 5 tons $\, \theta_4 \,$ o\`u
\moneqstar
 \theta_4 = \xi_2 = {9\over8} 
\monendstar 
et de deux demi-tons $\, \delta_4 \,$ avec
\moneqstar 
 \delta_4 = \xi_{-5} = {{2^8}\over{3^5}}  \, . 
\monendstar

\smallskip \noindent 
Nous  proposons de les construire dans ce paragraphe, en utilisant la
classification par type et mode introduite dans les paragraphes pr\'ec\'edents. On
dispose de trois types pour ces~21 gammes et pour chacun de ces types, de sept modes.
Nous proposons ici, mais 
il s'agit d'une convention et d'autres choix \'equivalents sont possibles,
de caract\'eriser ces trois types  par la ~\guillemotleft~distance~\guillemotright~
entre les deux demi-tons : jointifs
pour le type~$\,  \tau_1^4  \, $ et  distants d'un ton pour le second type $\, \tau_2^4 $. 
De plus, si l'on garde en m\'emoire l'hypoth\`ese d'invariance d'un type par permutation circulaire,
les deux  demi-tons du troisi\` eme  type~$\,  \tau_3^4  \, $ sont
distants de trois tons, ou bien de  deux tons par p\'eriodicit\'e : 
\moneqstar  
\left \{ \begin {array}{l}
\tau_1^4 =  \bigl(  \theta_4 ,\, \theta_4  ,\, \delta_4  ,\,  \delta_4  ,\, \theta_4 ,\, \theta_4  ,\,  \theta_4  \bigr) \\
\tau_2^4 =  \bigl( \theta_4 ,\, \theta_4  ,\, \delta_4  ,\,  \theta_4  ,\, \delta_4  ,\,  \theta_4 ,\,  \theta_4  \bigr) \\
\tau_3^4 =  \bigl( \theta_4 ,\, \theta_4  ,\, \delta_4  ,\,  \theta_4  ,\,  \theta_4 ,\,   \theta_4  ,\, \delta_4   \bigr)  \, .
\end {array} \right.
\monendstar 
%
%
Le troisi\`eme type~$\,  \tau_3^4  \,$ est appel\'e  ~\guillemotleft~classique~\guillemotright~    dans la suite.

\smallskip  \noindent
Pour ces trois types, les gammes fondamentales,  c'est-\`a-dire en mode de {\textsl {do}}, 
s'\'ecrivent facilement. On a :
\moneqstar  \left \{ \begin {array}{l}
{\cal G}(\tau_1^4 ,\,0) = \big( {\textsl {do}} ,\, {\textsl {r\'e}}  ,\, {\textsl {mi}}  ,\,{\textsl {fa}}
,\, {\textsl {sol}}\,\flat   ,\, {\textsl {la}}\,\flat  \, ,\, {\textsl {si}}\,\flat \, ,\,  {\textsl {do}}^{*} \big) =
\big( \Omega_0 ,\,  \Omega_2 ,\,  \Omega_4 ,\, \Omega_{-1} ,\,
\Omega_{-6}  ,\,  \Omega_{-4} ,\,  \Omega_{-2}   ,\, 2 \,  \Omega_0 \big) \\ 
{\cal G}(\tau_2^4 ,\,0) = \big( {\textsl {do}} ,\, {\textsl {r\'e}}  ,\, {\textsl {mi}}  ,\,{\textsl {fa}}
,\, {\textsl {sol}}   ,\, {\textsl {la}}\,\flat  \, ,\, {\textsl {si}}\,\flat \, ,\,  {\textsl {do}}^{*}\big) =
\big( \Omega_0 ,\,  \Omega_2 ,\,  \Omega_4 ,\, \Omega_{-1} 
,\,\Omega_{1}  ,\,  \Omega_{-4} ,\,  \Omega_{-2}   ,\, 2 \,  \Omega_0 \big) \\
{\cal G}(\tau_3^4 ,\,0) = \big(  {\textsl {do}} ,\, {\textsl {r\'e}}  ,\, {\textsl {mi}}  ,\,{\textsl {fa}}
,\, {\textsl {sol}}   ,\, {\textsl {la}}  \, ,\, {\textsl {si}}    ,\,  {\textsl {do}}^{*} \big) =
\big( \Omega_0 ,\,  \Omega_2 ,\,  \Omega_4 ,\, \Omega_{-1}
,\,\Omega_{1}  ,\,  \Omega_{3} ,\,  \Omega_{5}   ,\, 2 \,  \Omega_0 \big) 
\end {array} \right.
\monendstar

\noindent
ce qui permet (enfin !) de retrouver la gamme la plus ~\guillemotleft~\'el\'ementaire~\guillemotright~
de {\textsl {do}} majeur avec le type classique. 
Pour chacun de ces types, on dispose de sept modes distincts. Leurs noms sont obtenus en lisant les
noms des notes des gammes pr\'ec\'edentes, choisies pour la classification. Il
vient~:
\moneq \label{10.9-10-11} \left\{ \begin{array} {l}
{\rm modes \,\, de  \,\, type}  \,\,  \tau_1^4 \,:\,  {\textsl {do}}  ,\,  {\textsl {r\'e}}  ,\,
{\textsl {mi}}  ,\,   {\textsl {fa}} ,\,   {\textsl {sol}}\,\flat  ,\,   {\textsl {la}}\,\flat   ,\,   {\textsl {si}}\,\flat \\
{\rm modes \,\, de  \,\, type}  \,\,  \tau_2^4 \,:\,  {\textsl {do}}  ,\,  {\textsl {r\'e}} ,\,    {\textsl {mi}}  ,\,   {\textsl {fa}}
,\,   {\textsl {sol}}  ,\,   {\textsl {la}}\,\flat   ,\,   {\textsl {si}}\,\flat \\
{\rm modes \,\,  (classiques) \,\, de  \,\, type}  \,\,  \tau_3^4 \,:\,    {\textsl {do}}  ,\,  {\textsl {r\'e}} ,\,    {\textsl {mi}}  ,\,   {\textsl {fa}}
,\,   {\textsl {sol}}  ,\,   {\textsl {la}}   ,\,   {\textsl {si}} \, .
\end{array} \right. \monend

\bigskip \noindent
Pour chacun des modes, nous d\'etaillons la gamme de tonalit\'e  {\textsl {do}},
qui par d\'efinition commence par {\textsl {do}}, et ce processus fournit une suite possible des 21
gammes \`a sept notes.
Nous notons parfois les gammes $ \, {\cal G}(\tau^4 ,\,i) \, $
en rempla\c cant le num\'ero du mode par la note correspondante dans l'une des
listes pr\'ec\'edentes. On a bien s\^ur toujours  $ \, {\cal G}(\tau^4 ,\,0) = {\cal G}(\tau^4 ,\,{\textsl {do}}) \, $
et par exemple $ \, {\cal G}(\tau^4_3 ,\,5) =  {\cal G}(\tau^4_3 ,\, {\textsl {la}}) $. 
Nous explicitons maintenant les gammes qui caract\'erisent les sept modes pour chacun des trois types.
Pour le premier type $\,  \tau_1^4 ,\,$ la suite des modes est donn\'ee
\`a la premi\`ere des relations (\ref{10.9-10-11}). On retrouve d'abord~:
\moneqstar 
{\cal G}^4_1 =  {\cal G}(\tau^4_1 ,\,0) =  {\cal G}(\tau^4_1 ,\, {\textsl {do}})   =   \left\{ \begin {array}{l}
\bigl( {\textsl {do}} ,\, {\textsl {r\'e}}  ,\, {\textsl {mi}}  ,\,{\textsl {fa}}
,\, {\textsl {sol}}\,\flat   ,\, {\textsl {la}}\,\flat  ,\, {\textsl {si}}\,\flat
,\,  {\textsl {do}}^{*} \bigr) \\
\bigl( \Omega_0 ,\,  \Omega_2 ,\,  \Omega_4 ,\, \Omega_{-1} ,\,
\Omega_{-6}  ,\,  \Omega_{-4} ,\,  \Omega_{-2} ,\, 2 \,  \Omega_0 \bigr) \, . 
\end {array} \right. 
\monendstar
Ce mode est parfois appel\'e   ~\guillemotleft~mode arabe~\guillemotright. 
%
%
%
Les six autres modes s'obtiennent en permuttant petit \`a petit les 7 notes de cette gamme.
Nous les pr\'esentons avec une transposition en  {\textsl{do}} afin de pouvoir les comparer facilement.
Pour la premi\`ere permutation, on trouve 
\moneqstar
{\cal G}^4_2 =  {\cal G}(\tau^4_1 ,\,1) =  {\cal G}(\tau^4_1 ,\, {\textsl {r\'e}})   =   \left\{ \begin {array}{l}
\bigl( {\textsl {do}} ,\, {\textsl {r\'e}}  ,\, {\textsl {mi}}\,\flat  ,\,{\textsl {fa}}\,\flat
,\, {\textsl {sol}}\,\flat   ,\, {\textsl {la}}\,\flat  ,\, {\textsl {si}}\,\flat ,\,  {\textsl {do}}^{*} \bigr) \\
\bigl( \Omega_0 ,\,  \Omega_2 ,\,  \Omega_{-3} ,\, \Omega_{-8} ,\,
\Omega_{-6}  ,\,  \Omega_{-4} ,\,\Omega_{-2} ,\, 2 \,  \Omega_0 \bigr) \, . 
\end {array} \right. 
\monendstar
On recommence le m\^eme exercice pour les cinq modes de ce premier type \`a demi-tons jointifs.
Le r\'esultat de ce calcul est explicit\'e dans la liste qui suit. On a
\moneqstar
{\cal G}^4_3 =  {\cal G}(\tau^4_1 ,\,2) =  {\cal G}(\tau^4_1 ,\, {\textsl {mi}})   =  \left\{ \begin {array}{l}
\bigl( {\textsl {do}} ,\, {\textsl {r\'e}}\,\flat  ,\, {\textsl {mi}}\,\flat\flat
,\,{\textsl {fa}}\,\flat    ,\, {\textsl {sol}}\,\flat   ,\, {\textsl {la}}\,\flat
,\, {\textsl {si}}\,\flat   ,\,  {\textsl {do}}^{*} \bigr) \\
\bigl( \Omega_0 ,\,  \Omega_{-5} ,\,  \Omega_{-10} ,\, \Omega_{-8} ,\,
\Omega_{-6}  ,\,  \Omega_{-4} ,\,\Omega_{-2} ,\, 2 \,  \Omega_0 \bigr) \,, 
\end {array} \right. 
\monendstar
\moneqstar
{\cal G}^4_4 =  {\cal G}(\tau^4_1 ,\,3) = {\cal G}(\tau^4_1 ,\, {\textsl {fa}})   =  \left\{ \begin {array}{l}
\bigl( {\textsl {do}} ,\, {\textsl {r\'e}}\,\flat  ,\, {\textsl {mi}}\,\flat
,\,{\textsl {fa}}    ,\, {\textsl {sol}}  ,\, {\textsl {la}} ,\, {\textsl {si}}   ,\,  {\textsl {do}}^{*} \bigr) \\
\bigl( \Omega_0 ,\,  \Omega_{-5} ,\,  \Omega_{-3} ,\, \Omega_{-1} ,\,
\Omega_{1}  ,\,  \Omega_{3} ,\,\Omega_{5} ,\, 2 \,  \Omega_0 \bigr)
\end {array} \right.
\monendstar
et ce  mode est  appel\'e quelquefois   ~\guillemotleft~mode napolitain~\guillemotright.
On a ensuite 
\moneqstar
{\cal G}^4_5 =  {\cal G}(\tau^4_1 ,\,4) = {\cal G}(\tau^4_1 ,\, {\textsl {sol}}\,\flat)   =  \left\{ \begin {array}{l}
\bigl( {\textsl {do}} ,\, {\textsl {r\'e}}  ,\, {\textsl {mi}} ,\,{\textsl {fa}}\,\sharp ,\,
{\textsl {sol}}\,\sharp  ,\, {\textsl {la}}\,\sharp ,\, {\textsl {si}} ,\,  {\textsl {do}}^{*} \bigr) \\
\bigl( \Omega_0 ,\,  \Omega_{2} ,\, \Omega_{4} ,\, \Omega_{6} ,\,
\Omega_{8}  ,\,  \Omega_{10} ,\,\Omega_{5} ,\, 2 \,  \Omega_0 \bigr) \,, 
\end {array} \right.
\monendstar
dite  ~\guillemotleft~gamme par tons et sensible~\guillemotright, puis 
\moneqstar
{\cal G}^4_6 =  {\cal G}(\tau^4_1 ,\,5) = {\cal G}(\tau^4_1 ,\, {\textsl {la}}\,\flat)   =  \left\{ \begin {array}{l}
\bigl( {\textsl {do}} ,\, {\textsl {r\'e}}  ,\, {\textsl {mi}} ,\,{\textsl {fa}}\,\sharp ,\,
{\textsl {sol}}\,\sharp  ,\, {\textsl {la}} ,\, {\textsl {si}}\,\flat ,\,  {\textsl {do}}^{*} \bigr) \\
\bigl( \Omega_0 ,\,  \Omega_{2} ,\,  \Omega_{4} ,\, \Omega_{6} ,\,
\Omega_{8}  ,\,  \Omega_{3} ,\,\Omega_{-2} ,\, 2 \,  \Omega_0 \bigr) \,,
\end {array} \right.
\monendstar
ou mode  ~\guillemotleft~hypolydien~\guillemotright~  et 
\moneqstar
{\cal G}^4_7 =  {\cal G}(\tau^4_1 ,\,6) = {\cal G}(\tau^4_1 ,\, {\textsl {si}}\,\flat)   =  \left\{ \begin {array}{l}
\bigl( {\textsl {do}} ,\, {\textsl {r\'e}}  ,\, {\textsl {mi}} ,\,{\textsl {fa}}\,\sharp ,\,
{\textsl {sol}}  ,\, {\textsl {la}}\,\flat ,\, {\textsl {si}}\,\flat   \,,\,  {\textsl {do}}^{*} \bigr) \\
\bigl( \Omega_0 ,\,  \Omega_{2} ,\,  \Omega_{4} ,\, \Omega_{6} ,\,
\Omega_{1}  ,\,  \Omega_{-4} ,\,\Omega_{-2} ,\, 2 \,  \Omega_0 \bigr) \,, 
\end {array} \right.
\monendstar
le    mode ~\guillemotleft~lydien-phrygien~\guillemotright. 
Ces modes ne semblent pas compl\`etement  formalis\'es \`a notre connaissance.
Bien entendu, on peut ensuite transposer ces diff\'erentes gammes en une tonalit\'e $ \, \Omega_\ell \, $
arbitraire et construire par exemple une gamme $ \, {\cal G}^4_7 ({\textsl {sol}}) \, $;  
il est facile de d\'ecaler les num\'eros des notes d'une unit\'e puis de les nommer en prenant
la note correspondante \`a l'octave sup\'erieure en cas de besoin. Nous obtenons ainsi 
\moneqstar
 {\cal G}^4_7 ({\textsl {sol}}) = 
\bigl( \Omega_1 ,\,  \Omega_{3} ,\,  \Omega_{5} ,\, \Omega_{7} ,\,
\Omega_{2}  ,\,  \Omega_{-3} ,\,\Omega_{-1} ,\, 2 \,  \Omega_0 \bigr) =
\bigl( {\textsl {sol}} ,\,  {\textsl {la}} ,\,  {\textsl {si}} ,\,  {\textsl {do}} \,\sharp^* ,\,
 {\textsl {r\'e}}^* ,\,  {\textsl {mi}}\,\flat^* ,\, {\textsl {fa}}^* ,\,  {\textsl {sol}}^* \bigr) \, . 
\monendstar

\bigskip \noindent
Pour le second type $\,  \tau_2^4 ,\,$ caract\'eris\'e par le fait qu'un ton s\'epare les deux demi-tons, 
la suite des modes est donn\'ee par  la seconde relation (\ref{10.9-10-11}). 
Nous poursuivons aussi une num\'erotation globale des gammes \`a 7 notes. 
Nous choisissons  toujours de faire op\'erer ces modes sur la tonalit\'e de {\textsl{do}}.
Nous avons pour le premier mode
\moneqstar  
{\cal G}_8^4 = {\cal G}(\tau_2^4 ,\,0) =  {\cal G}(\tau_2^4 ,\, {\textsl {do}}) = \left\{  \begin {array}{l}
\bigl( {\textsl {do}} ,\, {\textsl {r\'e}}  ,\, {\textsl {mi}}  ,\,{\textsl {fa}} ,\,
{\textsl {sol}} ,\, {\textsl {la}}\,\flat ,\, {\textsl {si}}\,\flat ,\,  {\textsl {do}}^{*} \bigr) \\
\bigl( \Omega_0 ,\,  \Omega_2 ,\,  \Omega_4 ,\, \Omega_{-1} ,\,
\Omega_{1}  ,\,  \Omega_{-4} ,\,  \Omega_{-2}   ,\, 2 \,  \Omega_0 \bigr) \, .
\end {array} \right.
\monendstar
Notons que les gammes $ \, {\cal G}_1^4 \, $  et $ \, {\cal G}_8^4 \, $ sont diff\'erentes~: 
le  $ \, {\textsl {sol}}\,\flat \, $ de   $ \, {\cal G}_1^4 \, $ devient un sol pour $ \, {\cal G}_8^4  $.
On passe ensuite au mode de {\textsl {r\'e}} avec 
\moneqstar  
{\cal G}_9^4 = {\cal G}(\tau_2^4 ,\,1) =  {\cal G}(\tau_2^4 ,\, {\textsl {r\'e}}) = \left\{  \begin {array}{l}  
\bigl( {\textsl {do}} ,\, {\textsl {r\'e}}  ,\, {\textsl {mi}}\,\flat ,\,
{\textsl {fa}}  ,\, {\textsl {sol}} \,\flat  ,\, {\textsl {la}}\,\flat  ,\, {\textsl {si}}\,\flat ,\, {\textsl {do}}^{*} \bigr) \\
\bigl( \Omega_0 ,\,  \Omega_2 ,\,  \Omega_{-3} ,\, \Omega_{-1} ,\,
\Omega_{-6}  ,\,  \Omega_{-4} ,\,  \Omega_{-2}   ,\, 2 \,  \Omega_0 \bigr) \, . 
\end {array} \right.
\monendstar
Les autres modes de ce second type $ \, \tau^4_2 \, $ se construisent de la m\^eme mani\`ere~:
\moneqstar   
{\cal G}_{10}^4 = {\cal G}(\tau_2^4 ,\,2) =  {\cal G}(\tau_2^4 ,\, {\textsl {mi}}) = \left\{  \begin {array}{l}
\bigl( {\textsl {do}} ,\, {\textsl {r\'e}} \,\flat ,\, {\textsl {mi}}\,\flat ,\, {\textsl {fa}} \,\flat ,\,
{\textsl {sol}} \,\flat ,\, {\textsl {la}}\,\flat ,\, {\textsl {si}}\,\flat \,  {\textsl {do}}^{*} \bigr) \\
\bigl( \Omega_0 ,\,  \Omega_{-5} ,\,  \Omega_{-3} ,\, \Omega_{-8} ,\,
\Omega_{-6}  ,\,  \Omega_{-4} ,\,  \Omega_{-2}   ,\, 2 \,  \Omega_0 \bigr)
\end {array} \right. 
\monendstar
\moneqstar   
{\cal G}_{11}^4 = {\cal G}(\tau_2^4 ,\,3) =  {\cal G}(\tau_2^4 ,\, {\textsl {fa}}) = \left\{  \begin {array}{l}
\bigl( {\textsl {do}} ,\, {\textsl {r\'e}},\, {\textsl {mi}}\,\flat ,\,
{\textsl {fa}} ,\, {\textsl {sol}} ,\, {\textsl {la}} ,\, {\textsl {si}}   ,\,  {\textsl {do}}^{*} \bigr) \\
\bigl( \Omega_0 ,\,  \Omega_{2} ,\,  \Omega_{-3} ,\, \Omega_{-1} ,\,
\Omega_{1}  ,\,  \Omega_{3} ,\,  \Omega_{5}   ,\, 2 \,  \Omega_0 \bigr)
\end {array} \right.
\monendstar
\moneqstar   
{\cal G}_{12}^4 = {\cal G}(\tau_2^4 ,\,4) =  {\cal G}(\tau_2^4 ,\, {\textsl {sol}}) = \left\{  \begin {array}{l}
\bigl( {\textsl {do}} ,\, {\textsl {r\'e}}\,\flat ,\, {\textsl {mi}}\,\flat ,\,
{\textsl {fa}} ,\, {\textsl {sol}} ,\, {\textsl {la}} ,\, {\textsl {si}}\,\flat ,\,  {\textsl {do}}^{*} \bigr) \\
\bigl( \Omega_0 ,\,  \Omega_{-5} ,\,  \Omega_{-3} ,\, \Omega_{-1} ,\,
\Omega_{1}  ,\,  \Omega_{3} ,\,  \Omega_{-2}   ,\, 2 \,  \Omega_0 \bigr)
\end {array} \right.
\monendstar
\moneqstar  
{\cal G}_{13}^4 = {\cal G}(\tau_2^4 ,\,5) =  {\cal G}(\tau_2^4 ,\, {\textsl {la}}\,\flat)  = \left\{  \begin {array}{l}
\bigl( {\textsl {do}} ,\, {\textsl {r\'e}} ,\, {\textsl {mi}} ,\,
{\textsl {fa}} \,\sharp ,\, {\textsl {sol}}\,\sharp ,\, {\textsl {la}} ,\, {\textsl {si}}  ,\, {\textsl {do}}^{*} \bigr) \\
\bigl( \Omega_0 ,\,  \Omega_{2} ,\,  \Omega_{4} ,\, \Omega_{6} ,\,
\Omega_{8}  \,\,,\,  \Omega_{3} ,\,  \Omega_{5}   ,\, 2 \, \Omega_0 \bigr)
\end {array} \right.
\monendstar
\moneqstar  
{\cal G}_{14}^4 =  {\cal G}(\tau_2^4 ,\,6) =  {\cal G}(\tau_2^4 ,\, {\textsl {si}}\,\flat)  = \left\{  \begin {array}{l}
\bigl( {\textsl {do}} ,\, {\textsl {r\'e}} ,\, {\textsl {mi}} ,\,
{\textsl {fa}} \,\sharp ,\, {\textsl {sol}} ,\, {\textsl {la}} ,\, {\textsl {si}}\,\flat ,\,  {\textsl {do}}^{*} \bigr) \\
\bigl( \Omega_0 ,\,  \Omega_{2} ,\,  \Omega_{4} ,\, \Omega_{6} ,\,
\Omega_{1}  \,\,,\,  \Omega_{3} ,\,  \Omega_{-2}   ,\, 2 \, \Omega_0 \bigr) \, .
\end {array} \right.
\monendstar

\smallskip \noindent 
Ce  dernier mode est parfois appel\'e    ~\guillemotleft~mode de Bartok~\guillemotright~
et nous renvoyons au livre de  Ern\"o Lendvai (1971) 
pour une analyse des diff\'erentes structures musicales cr\'e\'ees par B\'el\`a Bart\'ok.
On  retrouve aussi ce mode dans  la musique indienne sous le nom  de  ~\guillemotleft~Raga Vachaspati~\guillemotright. 
Signalons aussi l'\oe uvre de  Christian Clav\`ere  (2021),  cr\'eation  musicale
au violon inspir\'ee par la gamme $\, {\cal G}_{14}^4 $.

\bigskip \noindent
Pour le type classique $\,  \tau_3^4 ,\,$ la suite des modes est donn\'ee par  la  troisi\`eme
relation~(\ref{10.9-10-11}), 
la suite plus courante des notes
~\guillemotleft~sans alt\'eration~\guillemotright.
Les tonalit\'es de la gamme de {\textsl {do}} s'\'ecrivent dans ces sept  modes en suivant la m\^eme approche. 
On a d'abord
\moneqstar  
{\cal G}_{15}^4 = {\cal G}(\tau_3^4 ,\,0) =  {\cal G}(\tau_3^4 ,\, {\textsl {do}}) = \left\{  \begin {array}{l}
\bigl( {\textsl {do}} ,\, {\textsl {r\'e}}  ,\, {\textsl {mi}}  ,\,{\textsl {fa}} ,\,
{\textsl {sol}} ,\, {\textsl {la}} ,\, {\textsl {si}} ,\,  {\textsl {do}}^{*} \bigr) \\
\bigl( \Omega_0 ,\,  \Omega_2 ,\,  \Omega_4 ,\, \Omega_{-1} ,\,
\Omega_{1}  ,\,  \Omega_{3} ,\,  \Omega_{5}   ,\, 2 \,  \Omega_0 \bigr) \, .
\end {array} \right.
\monendstar
C'est le mode ionien, ou mode ~\guillemotleft~majeur~\guillemotright~  de la musique occidentale,
d\'efini avec la gamme de  {\textsl{do}} majeur. C'est bien s\^ur  ~\guillemotleft~la~\guillemotright~  gamme de r\'ef\'erence ! 
On peut ensuite transposer ce mode de  {\textsl{do}}, par exemple en $ \, {\textsl {si}}\,\flat $. 
On retranche deux unit\'es aux notes $ \, \Omega_\ell \, $ de la gamme $ \, {\cal G}_{15}^4 \, $
initiale et on trouve $ \,  {\cal G}_{15}^4 (\Omega_{-2}) =  {\cal G}_{15}^4 ({\textsl {si}}\,\flat) \, $: 
\moneqstar
{\cal G}_{15}^4 ({\textsl {si}}\,\flat) \!=\!
\bigl( \Omega_{-2} ,\,  \Omega_0 ,\,  \Omega_2 ,\, \Omega_{-3} ,\,
\Omega_{-1}  ,\,  \Omega_{1} ,\,  \Omega_{3}   ,\, 2 \,  \Omega_{-2} \bigr) = 
\bigl( {\textsl {si}}\,\flat ,\, {\textsl {do}}^* ,\, {\textsl {r\'e}}^* ,\,   {\textsl {mi}}\,\flat^* ,\, {\textsl {fa}}^* ,\, {\textsl {sol}}^* ,\,
{\textsl {la}}^* ,\,  {\textsl {si}}\,\flat^* \bigr) .
\monendstar
On obtient tout simplement la tonalit\'e de de $ \, {\textsl {si}}\,\flat \, $ majeur.
On poursuit l'\'enum\'eration des  diff\'erents modes classiques associ\'es au type~$ \, \tau^4_3 $. La note suivante est un  {\textsl {r\'e}} ; on
obtient le mode de  {\textsl {r\'e}} ou mode dorien.  Il a le charme de la musique du moyen-\^age
et est utilis\'e par Jean-S\'ebastien Bach dans les ann\'ees 1715 pour 
sa {\it Toccata et fugue dorienne} (en {\textsl {r\'e}} mineur),  BWV~538. 
Nous l'\'ecrivons en commen\c cant par la note  {\textsl {do}}~: 
\moneqstar  
{\cal G}_{16}^4 = {\cal G}(\tau_3^4 ,\,1) =  {\cal G}(\tau_3^4 ,\, {\textsl {r\'e}}) = \left\{  \begin {array}{l}
\bigl( {\textsl {do}} ,\, {\textsl {r\'e}}  ,\, {\textsl {mi}}\,\flat ,\
{\textsl {fa}} \,,\, {\textsl {sol}} ,\, {\textsl {la}}  ,\, {\textsl {si}}\,\flat ,\,  {\textsl {do}}^{*} \bigr)  \\
\bigl( \Omega_0 ,\,  \Omega_2 ,\,  \Omega_{-3} ,\, \Omega_{-1} ,\,
\Omega_{1}  ,\,  \Omega_{3} ,\,  \Omega_{-2}   ,\, 2 \,  \Omega_0 \bigr) \, . 
\end {array} \right. 
\monendstar
On ne  confondra pas $\, {\cal G}_{16}^4 ({\textsl {si}}\,\flat) \, $
donn\'ee par la suite
$ \, \bigl( \Omega_{-2} ,\,  \Omega_0 ,\,  \Omega_{-5} ,\, \Omega_{-3} ,\,
\Omega_{-1}  ,\,  \Omega_{1} ,\,  \Omega_{-4}   ,\, 2 \,  \Omega_{-2} \bigr) $, 
 c'est-\`a-dire 
\moneqstar 
 {\cal G}_{16}^4 ({\textsl {si}}\,\flat) = \bigl( {\textsl {si}}\,\flat  ,\, {\textsl {do}}^*  ,\, {\textsl {r\'e}}\,\flat^{*}   ,\, {\textsl {mi}}\,\flat^{*}  , 
{\textsl {fa}}^{*}  \,,\, {\textsl {sol}}^{*}  ,\, {\textsl {la}}\,\flat^{*}   ,\, {\textsl {si}}\,\flat^{*} \bigr) ,  
\monendstar
 avec la gamme de $ \, {\textsl {si}}\,\flat \, $ majeur rappel\'ee quelques lignes ci-dessus ! 
Le mode de  {\textsl {mi}}, ou mode phrygien, commence la gamme par un demi-ton. Il est associ\'e \`a la gamme qui suit~: 
\moneqstar  
{\cal G}_{17}^4 = {\cal G}(\tau_3^4 ,\,2) =  {\cal G}(\tau_3^4 ,\, {\textsl {mi}}) = \left\{  \begin {array}{l}
\bigl( {\textsl {do}} ,\,  {\textsl {r\'e}}\,\flat   ,\, {\textsl {mi}}\,\flat ,\,
{\textsl {fa}}  ,\, {\textsl {sol}}   ,\, {\textsl {la}}\,\flat  ,\, {\textsl {si}}\,\flat     ,\,  {\textsl {do}}^{*} \bigr) \\
\bigl( \Omega_0 ,\,  \Omega_{-5} ,\,  \Omega_{-3} ,\, \Omega_{-1} ,\,
\Omega_{1}  ,\,  \Omega_{-4} ,\,  \Omega_{-2}   ,\, 2 \,  \Omega_0 \bigr) \,. 
\end {array} \right.
\monendstar
Il est explicitement utilis\'e par Jehan Alain (1930) dans sa {\it Ballade en mode phrygien} pour orgue.
Le mode suivant, ou mode de  {\textsl {fa}},   mode lydien~:  
\moneqstar  
{\cal G}_{18}^4 = {\cal G}(\tau_3^4 ,\,3) =  {\cal G}(\tau_3^4 ,\, {\textsl {fa}}) = \left\{  \begin {array}{l}
\bigl( {\textsl {do}} ,\,  {\textsl {r\'e}}  ,\,  {\textsl {mi}} ,\,
{\textsl {fa}}\,\sharp ,\, {\textsl {sol}} ,\, {\textsl {la}}  ,\, {\textsl {si}} ,\,  {\textsl {do}}^{*} \bigr) \\
\bigl( \Omega_0 ,\,  \Omega_{2} ,\,  \Omega_{4} ,\, \Omega_{6} ,\,
\Omega_{1} ,\, \Omega_{3} ,\, \Omega_{5}   ,\, 2 \, \Omega_0 \bigr) \, . 
\end {array} \right.
\monendstar
Un exemple c\'el\`ebre du mode  {\textsl {fa}}  de date de 1825~: c'est le 
  troisi\`eme mouvement du {\it Quatuor \`a  cordes num\'ero 15} 
de Ludwig van Beethoven (voir par exemple Daniel K. L. Chua, 1995). 
%
On ne doit bien s\^ur pas confondre le mode de  {\textsl {fa}}  avec la gamme de {\textsl {fa}} majeur
qui est \'egale \`a $ \, {\cal G}_{15}^4 ({\textsl {fa}}) \, $ et poss\`ede un $ \, {\textsl {si}}\,\flat $. 
On le retrouve aussi dans la musique populaire.  
%

\smallskip \noindent 
Le mode de {\textsl {sol}} est \'eclatant ; on le nomme aussi  mode  mixolydien~: 
\moneqstar  
{\cal G}_{19}^4 = {\cal G}(\tau_3^4 ,\,4) =  {\cal G}(\tau_3^4 ,\, {\textsl {sol}}) = \left\{  \begin {array}{l}
\bigl( {\textsl {do}} ,\,  {\textsl {r\'e}}  ,\,  {\textsl {mi}} ,\,
{\textsl {fa}}  ,\,  {\textsl {sol}}   ,\, {\textsl {la}}  ,\, {\textsl {si}}\,\flat   ,\,  {\textsl {do}}^{*} \bigr) \\
\bigl( \Omega_0 ,\,  \Omega_{2} ,\,  \Omega_{4} ,\, \Omega_{-1} ,\,
\Omega_{1} ,\, \Omega_{3} ,\, \Omega_{-2} ,\, 2 \, \Omega_0 \bigr) \, . 
\end {array} \right.
\monendstar
On le retrouve par exemple  dans la pi\`ece pour piano {\it L'Isle joyeuse}, de Claude  Debussy
(Dmitri Tymoczk, 2004) 
ou dans la musique populaire avec  plusieurs chansons des Beatles 
comme {\it Norwegian wood} ou {\it If I needed someone} 
 (voir par exemple l'ouvrage \'edit\'e en 2006 par  Kenneth Womack and Todd~F.~Davis). 
%
Il est bien entendu diff\'erent de la gamme de  {\textsl {sol}}  majeur,  c'est-\`a-dire $ \, {\cal G}_{15}^4 ({\textsl {sol}}) \, $
qui poss\`ede un $ \, {\textsl {fa}}\,\sharp \, $ ce qui n'est pas le cas de
$ \, {\cal G}_{19}^4 ({\textsl {sol}}) = \bigl( {\textsl {sol}} ,\,  {\textsl {la}}  ,\,  {\textsl {si}} ,\,
{\textsl {do}}^*  ,\,  {\textsl {r\'e}}^*   ,\, {\textsl {mi}}^*  ,\, {\textsl {fa}}^*  ,\,  {\textsl {sol}}^{*} \bigr) $.

\smallskip \noindent 
Le mode de  {\textsl {la}} ou mode  \'eolien, est aussi le mode  ~\guillemotleft~mineur~\guillemotright~
de la musique occidentale. Nous le donnons, comme les autres modes, dans la tonalit\'e de   {\textsl {do}},
 c'est-\`a-dire  avec la gamme de  {\textsl {do}} mineur~: 
\moneqstar  
{\cal G}_{20}^4 = {\cal G}(\tau_3^4 ,\,5) =  {\cal G}(\tau_3^4 ,\, {\textsl {la}}) = \left\{  \begin {array}{l}
\bigl( {\textsl {do}} \,\,  ,\, {\textsl {r\'e}}   \,,\, {\textsl {mi}}\,\flat ,\,
{\textsl {fa}} ,\,  {\textsl {sol}}   ,\, {\textsl {la}} \,\flat ,\,  {\textsl {si}}\,\flat   ,\,  {\textsl {do}}^{*} \bigr) \\
\bigl( \Omega_0 ,\,  \Omega_{2} ,\,  \Omega_{-3} ,\, \Omega_{-1} ,\,
\Omega_{1}  ,\,  \Omega_{-4} ,\,  \Omega_{-2}   ,\, 2 \,  \Omega_0 \bigr) \, . 
\end {array} \right.
\monendstar
Enfin, le mode de  {\textsl {si}}, ou mode locrien, commence, comme le mode de  {\textsl {mi}}, la gamme par un demi-ton :
\moneqstar  
{\cal G}_{21}^4 = {\cal G}(\tau_3^4 ,\,6) =  {\cal G}(\tau_3^4 ,\, {\textsl {si}}) = \left\{  \begin {array}{l}
\bigl( {\textsl {do}} ,\,  {\textsl {r\'e}} \,\flat ,\, {\textsl {mi}}\,\flat ,\,
{\textsl {fa}} ,\,  {\textsl {sol}}\,\flat ,\, {\textsl {la}} \,\flat ,\,  {\textsl {si}}\,\flat , \,  {\textsl {do}}^{*} \bigr) \\
\bigl( \Omega_0 ,\,  \Omega_{-5} ,\,  \Omega_{-3} ,\, \Omega_{-1} ,\,
\Omega_{-6} ,\,  \Omega_{-4} ,\, \Omega_{-2} ,\, 2 \,  \Omega_0 \bigr) \,.
\end {array} \right.
\monendstar
%
%
Il est beaucoup plus rare que les modes pr\'ec\'edents.
L'ouvrage de Vincent Persichetti (1961) 
donne toute une liste de compositeurs du 20$^{\rm e}$ si\`ecle qui ont utilis\'e le mode locrien ;
nous retenons par exemple l'\oe uvre {\it Ludus tonalis}
de Paul Hindemith (1942). 

\bigskip \noindent
Les 21 modes qui caract\'erisent les  gammes \`a 7 notes de la famille expos\'ee 
dans ce travail font  appara\^itre les modes majeurs et mineurs, modes dominants de la musique
occidentale, ainsi que des modes issus des  anciens modes grecs.
Les autres modes qui ne sont pas toujours fond\'es sur le t\'etracorde~;
ils semblent avoir \'et\'e d\'ecouverts ou reconstruits plus tard, comme les noms
qu'ils portent l'indiquent.
Enfin,  les structures \'etudi\'ees ici sont contraintes par 
deux sortes d'intervalles seulement, ton  et demi-ton, ce qui ne correspond pas toujours
\`a la pratique des musiciens et des compositeurs,  
comme  en t\'emoigne par exemple l'existence de la note sensible au sein du mode mineur.

\bigskip \smallskip  \noindent {\bf \large  Conclusion }

\smallskip \noindent
Dans cette contribution, nous sommes partis d'une  division tr\`es simple de l'octave
en une quinte et une quarte  pour
former des gammes primitives \`a deux notes. Nous avons montr\'e ensuite qu'avec
cette notion \'el\'ementaire,  
l'invariance multiplicative des tonalit\'es induit
naturellement un nombre de notes 
{\it a priori} infini. Nous
avons alors  propos\'e une d\'efinition r\'ecursive d'une famille de gammes qui
permet de reconstruire la gamme de Pythagore et ses diverses  variantes. Une
classification en type et mode utilise la structure cyclique du motif form\'e par la
suite des tons et demi-tons. Nous avons d\'etaill\'e l'exemple des  gammes \`a sept notes, en
insistant sur le fait que deux structures seulement  sur les vingt et une possibles
 sont utilis\'ees dans la musique classique occidentale.
Observons qu'une  telle   construction th\'eorique sur un sujet tr\`es classique  a 
\'et\'e quasi  formalis\'ee 
dans des travaux ant\'erieurs
(voir par exemple Franck  Jedrzejewski, 2007). 

\bigskip \noindent
Avec une approche qui introduit une grande vari\'et\'e de gammes,
nous pouvons envisager la musique globalement, sans se limiter aux cons\-tructions
occidentales, en particulier pour les gammes pentatoniques abord\'ees  dans cette contribution.
Nous renvoyons le lecteur \`a l'article d'Alain Boudet (2006) pour une pr\'esentation
des gammes occidentales, arabes,  chinoises et indiennes.
Pour une \'etude approfondie de la musique arabe du 20e si\`ecle,
nous nous r\'ef\'erons par exemple \`a Ahmed et Mohamed Elhabib Hachlef (1993)
et pour la  musique chinoise, \`a l'ouvrage
de Fran\c cois Picard (1991).
Le jazz utilise aussi des gammes pentatoniques et nous sugg\'erons au lecteur
les ouvrages de Jacques Siron (1992) et Mark Levine (1995).
Enfin,  une introduction \`a  la musique indienne est propos\'ee par Patrick Moutal (1987).
%

\bigskip \noindent
Nous avons abord\'e plusieurs aspects math\'ematiques
avec une approche volontairement \'el\'emen\-taire. 
Il est bien s\^ur possible d'enrichir la construction pr\'ec\'edente
{\it via} par un algorithmique plus complexe prenant comme base de
construction non seulement le troisi\`eme harmonique 
de la vibration d'une corde mais aussi
celui associ\'e au  cinqui\`eme de la longueur, comme le propose
la gamme de Zarlino.
Le lecteur int\'eress\'e pourra approfondir le\br
sujet avec
les ``gammes bien form\'ees'' (Norman Carey et David Clampitt, 1989),
le lien avec l'informatique th\'eorique (Marc Chemillier, 2004), 
l'importance de la transform\'ee de Fourier (Emmanuel Amiot, 2016),
les cycles hamiltoniens pour l'analyse des musiques popu\-laires (Moreno Andratta et Gilles Baroin, 2016) 
ou les applications musicales de la th\'eorie des n\oe uds (Franck Jedrzejewski, 2019). 
%
Enfin, il semble utile \`a ce niveau de poursuivre l'exploration 
des musiques qui  utilisent les gammes  propos\'ees dans ce m\'emoire
%
%
et celles associ\'ees \`a 
une \'ecriture  qui rel\`eve  d'un tout autre type de construction formelle.

\newpage 
\bigskip \bigskip   \noindent {\bf \large    Annexe A) \quad   Une gamme \`a 8 notes ?}

\smallskip \noindent
Nous d\'eveloppons dans cette annexe une tentative de brisure de ton d'une gamme pentatonique
dans le cas o\`u l'on n'aurait pas pris garde au fait que lors de l'\'etape pr\'ec\'edente de l'algorithme,
on brise une quarte par une seconde,  
le nouveau demi-ton $ \, \delta = {4\over3} /  {9\over8} = {32\over27} \, $ est plus grand
que le nouveau ton $ \, \theta = {9\over8} \, $
puisque $ \, {32\over27} > {9\over8} $. 
On repart pour fixer les id\'ees d'une gamme pentatonique majeure 
\moneqstar
{\cal G}^3_1 =  \bigl( {\textsl {do}} ,\,  {\textsl {r\'e}} ,\,
 {\textsl {mi}} ,\, {\textsl {sol}} ,\, {\textsl {la}} ,\, {\textsl {do}}^*   \bigr) \,\,
 = \,\,  \bigl( \Omega_0 ,\, \Omega_2 ,\, \Omega_4 ,\, \Omega_1 ,\, \Omega_3
,\, 2 \, \Omega_0 \bigr) 
 \monendstar
 qui comporte avec notre hypoth\`ese trois tons $ \, \theta =  {9\over8} \, $ et deux demi-tons
 $ \, \delta = {32\over27} $.
La structure multiplicative $ \, {\cal S} ({\cal G}^3_1) \, $ de cette gamme
est donn\'ee par la relation
$ \,  {\cal S} ({\cal G}^3_1) = (\theta ,\, \theta ,\, \delta ,\, \theta ,\, \delta ) $.

\bigskip \noindent
On poursuit l'algorithme de construction lin\'eaire  et on coupe le ton $\, \theta \, $ avec le demi-ton $ \, \delta $.
Le nouveau ton $ \, \theta^* \, $ est donn\'e par la relation
$ \,  \theta^* = \delta =  {32\over27} \, $ et le nouveau demi-ton $ \, \delta^* \,$ s'obtient avec une division~:
$ \, \delta^* = {{\theta}\over\delta} = {{3^2}\over{2^3}} / {{2^5}\over{3^3}} = {{3^5}\over{2^8}} $.
Remarquons qu'il est inf\'erieur \`a 1 ! 

\bigskip \noindent
On construit une nouvelle gamme d'une part en rempla\c cant  chaque ton $\, \theta \, $ de la structure  $ \, {\cal S} ({\cal G}^3_1) \, $
par un nouveau ton $\, \theta^* \, $ suivi d'un nouveau  demi-ton $ \, \delta^* \, $
et d'autre part en rempla\c cant chaque demi-ton $ \, \delta \, $ par un ton $ \, \theta^* $.
On obtient ainsi une nouvelle structure
\moneqstar
{\cal S}' = (\theta^* ,\, \delta^* ,\, \theta^* ,\, \delta^* ,\, \theta^* ,\, \theta^* ,\, \delta^*  ,\, \theta^* )  
\monendstar
avec 5 tons $ \, \theta^* \, $ et 3 demi-tons $\, \delta^* $.
\`A partir du {\it do} initial, on construit une ~\guillemotleft~gamme \`a 8 notes~\guillemotright~
qui s'\'ecrit
$\, {\cal G}^8  =   \bigl( \Omega_0 ,\, \Omega_{-3} ,\, \Omega_2 ,\,  \Omega_{-1} ,\, \Omega_4 ,\, \Omega_1 ,\, \Omega_{-2} ,\, 
\Omega_3 ,\, 2 \, \Omega_0 \bigr) $, 
 c'est-\`a-dire  
\moneqstar
{\cal G}^8  = \bigl( {\textsl {do}} ,\,  {\textsl {mi}}\,\flat ,\,  {\textsl {r\'e}}  ,\, {\textsl {fa}}  ,\, 
 {\textsl {mi}} ,\, {\textsl {sol}} ,\,  {\textsl {si}}\,\flat ,\, {\textsl {la}} ,\, {\textsl {do}}^*   \bigr) \, .
\monendstar
Comme le demi-ton $ \, \delta^* \, $ est inf\'erieur \`a 1, cette suite de notes n'est pas de fr\'equences
croissantes et l'hypoth\`ese initiale (\ref{2.2}) de monotonie n'est plus satisfaite.
Est-ce encore une gamme~?

\bigskip  \bigskip      \noindent {\bf \large    Annexe B) \quad   Preuve de la proposition 3}

\smallskip \noindent 
La Proposition 3 s'\'etablit par
r\'ecurrence sur $ \, k$. Pour $\, k=1 ,\,$    
la premi\`ere des relations~(\ref{7.15-a-7.17}) est claire et montre que
les relations (\ref{7.10-7.11}) ont bien 
lieu, avec de plus
%
$ \, \varepsilon_1 = +1 $.
Dans le cas g\'en\'eral, on a toujours
\moneq \label{7.19}
\delta_k < \theta_k
\monend
et les relations (\ref{6.4}) \`a (\ref{6.7}) ont \'et\'e construites pour maintenir cette
propri\'et\'e. On d\'eduit donc de (\ref{7.19}) et de l'hypoth\`ese de r\'ecurrence
l'in\'egalit\'e suivante
$ \,  {{3^{^{\,- \varepsilon_k \, T_k}}} \over{2^{^{ \ell(\,- \varepsilon_k \, T_k)}}}}  < {{3^{^{ \,\varepsilon_k \, D_k}}}
\over{2^{^{ \ell(\, \varepsilon_k \, D_k)}}}} $, soit
\moneq \label{7.20}
\ell(\varepsilon_k \, D_k) -  \ell(- \varepsilon_k \, T_k) < \varepsilon_k \,(T_k + D_k
) \, {{{\rm log}\,3}\over{{\rm log}\,2}} \,.
\monend
De plus, compte tenu de la relation (\ref{7.9}), on a
\moneq \label{7.21}
\ell(\varepsilon_k \, D_k) \leq \varepsilon_k \,D_k \,
{{{\rm log}\,3}\over{{\rm log}\,2}} < \ell(\varepsilon_k \, D_k) + 1
\monend
et
\moneqstar
\ell(-\varepsilon_k \, T_k) \leq -\varepsilon_k \,T_k \,
{{{\rm log}\,3}\over{{\rm log}\,2}} < \ell(-\varepsilon_k \, T_k) + 1 \,
\monendstar
dont on d\'eduit par changement de signe 
\moneq \label{7.21bis}
-\ell(-\varepsilon_k \, T_k)  - 1  < \varepsilon_k \,T_k \, {{{\rm
log}3}\over{{\rm log}\,2}} \,\,\leq \,\, -\ell(-\varepsilon_k \, T_k) \, .
\monend
%
En additionnant les  in\'egalit\'es de droite des relations (\ref{7.21}) et (\ref{7.21bis}), il vient
\moneq \label{7.22}
 \varepsilon_k \,(T_k + D_k ) \, {{{\rm log}\,3}\over{{\rm log}\,2}} <
\ell(\varepsilon_k \, D_k) -  \ell(- \varepsilon_k \, T_k) \,+\, 1 \,.
\monend
Les relations (\ref{7.20}) et (\ref{7.22}) \'etablissent donc, compte tenu de la d\'efinition
(\ref{7.9}), la propri\'et\'e qui suit :
\moneq \label{7.23}
\ell(\varepsilon_k \, D_k) -  \ell(- \varepsilon_k \, T_k)  = \ell \bigl(
 \varepsilon_k \,(T_k + D_k ) \bigr) \,,\quad k \geq 1 \,.
\monend
%

\bigskip \noindent
Nous supposons dans ce paragraphe que la famille de gammes $\, {\cal G}^k \,$ de
num\'ero $ \, k \, $ satisfait l'hypoth\`ese $ \, \theta_k < \delta_k^2 $.  On a alors,
compte tenu des relations (\ref{6.4}) \`a (\ref{6.7})~:

\noindent
$  T_{k+1} = T_k + D_k $,  $ \,  D_{k+1} = T_k $, $\, \theta_{k+1} = \delta_k \, $
et $ \, \delta_{k+1} = {{\theta_k}\over{\delta_k}} $.
Par suite, $\, \smash { \theta_{k+1}  = \xi\ib{ \,-\varepsilon_k \, T_k} } =
\xi\ib{ \,-\varepsilon_k \, D_{k+1}} \, $  
et la premi\`ere relation   (\ref{7.10-7.11})  est vraie
\`a l'ordre $ \, k + 1 $,  en remarquant que
\moneq \label{7.24}
\varepsilon_{k+1} = - \varepsilon_k \,\,\,\,\, {\rm si} \,\,\,\,\, \theta_k < \delta_k^2 \, .
\monend
On a alors le calcul suivant 

\smallskip \smallskip \noindent  $
\delta_{k+1} = {{\xi\ib{ \varepsilon_k \,D_k}}\over{\xi\ib{ -\varepsilon_k \, T_k}}} =
{{3^{^{ \varepsilon_k \, D_k}}}\over{2^{^{ \ell( \varepsilon_k \, D_k)}}}} \,\,
{{2^{^{ \ell(- \varepsilon_k \, T_k)}}}\over{3^{^{ -\varepsilon_k \, T_k}}}} =
{{3^{^{ \varepsilon_k (T_k + D_k)}}} \over {2^{^{\ell( \varepsilon_k \, D_k) - \ell(- \varepsilon_k \, T_k)}}}} $

\smallskip \noindent  $ \qquad
=  {{3^{^{ \varepsilon_k (T_k + D_k)}}} \over {2^{^{\ell( \varepsilon_k (T_k + D_k) ) }}}} \, $
\quad compte tenu de la relation (\ref{7.23})

\smallskip \noindent  $ \qquad
=  {{3^{^{ -\varepsilon_{k+1} T_{k+1}}}} \over {2^{^{\ell( -\varepsilon_{k+1} \, T_{k+1}) }}}} = \xi\ib{ -\varepsilon_{k+1}
\,T_{k+1}} \, $  \quad au vu de la relation (\ref{7.24})

\smallskip \noindent
et la seconde relation de (\ref{7.10-7.11}) passe \`a l'ordre sup\'erieur.

\bigskip \noindent
Dans l'autre cas de figure o\`u l'on a  $\, \theta_k > \delta_k^2 $,  on d\'eduit
des relations (\ref{6.4}) \`a (\ref{6.7})~:    $  \, T_{k+1} = T_k $,
$ \, D_{k+1} = T_k + D_k $, $ \, \theta_{k+1} = {{\theta_k}\over{\delta_k}} \, $ et
$ \, \delta_{k+1} = \delta_k $.
Il r\'esulte alors  
des  relations (\ref{7.10-7.11}) :

\smallskip \noindent
$ \delta_{k+1} = \delta_k = \xi\ib{ \,-\varepsilon_k \,T_{k}}
= \xi\ib{ \,-\varepsilon_k \,T_{k+1}} \, $  et 
la seconde relation  (\ref{7.10-7.11})
est encore vraie \`a l'ordre $\, k\!+\!1 , \,$ avec
\moneqstar
\varepsilon_{k+1} =  \varepsilon_k \,\,\,\,\, {\rm si} \,\,\,\,\,  \theta_k > \delta_k^2  \,.
\monendstar
On a alors $ \, \theta_{k+1} = {{\xi\ib{ \, \varepsilon_k \,D_k}}
\over{\xi\ib{\, -\varepsilon_k \, T_k}}} =
\xi\ib{ \varepsilon_{k} (T_{k} + D_k ) } \,  $  compte tenu de (\ref{7.23}) et
du calcul men\'e plus haut, donc
$ \,  \theta_{k+1} = \xi\ib{\, \varepsilon_{k+1} D_{k+1} } \, $ et la
premi\`ere relation   (\ref{7.10-7.11})  est vraie
\`a l'ordre suivant $\, k\!+\!1 .\,$ La proposition~3 est \'etablie.   \hfill $ \square $

\bigskip \smallskip  \noindent {\bf \large     Remerciements }

\noindent   
Sans les remarques de St\'ephane Dubois,  Jean-Fran\c{c}ois Gonzales, Maurice
Rosset et  Michel Rosset, ce  texte 
n'aurait pas vu le jour fin 1999. 
Cette \'edition   a  b\'en\'efici\'e
des encou\-ragements de Thomas H\'elie, Pierre Mazet  et Roger Ohayon et de 
la lecture  tr\`es constructive
de Johan Broekaert, Andr\'e Calvet  et Christian Clav\`ere.
Certaines de leurs suggestions ont parfois \'et\'e int\'egr\'ees dans le texte. 
Ce travail doit beaucoup aux encou\-ragements de Marc Chemiller qui a remarqu\'e le lien entre 
la brisure de ton et le processus d'accordage des instruments dans la M\'esopotamie ancienne.
Nous remercions \'egalement Emmanuel Amiot et Franck Jedrzejewski qui m'ont transmis des
r\'ef\'erences bibliographiques r\'ecentes. 
Sophie Mougel a relu et annot\'e une version quasi-termin\'ee de ce m\'emoire.
Enfin, Maxime Chupin a point\'e plusieurs corrections orthographiques et typographiques. 
Qu'ils soient tous remerci\'es pour leur aide~!

\bigskip \smallskip  \noindent {\bf \large R\'ef\'erences bibliographiques }

\hangindent=7.3mm \hangafter=1 \noindent
C. Abromont, E. de Montalembert, {\it Guide de la th\'eorie de la musique}, Fayard-Lemoine, Paris, 2001.

\hangindent=7.3mm \hangafter=1 \noindent 
E. Amiot,
{\it Music through Fourier space: discrete Fourier transform in music theory},
Springer, 2016. 

\smallskip \hangindent=7.3mm \hangafter=1 \noindent
 M. Andreatta, {\it M\'ethodes alg\'ebriques en musique et musicologie du 20e si\`ecle :
   aspects th\'eoriques, analytiques et compositionnels}, Th\`ese,
 \'Ecole des Hautes \'Etudes en Sciences Sociales, 2003.

\smallskip \hangindent=7.3mm \hangafter=1 \noindent
M. Andreatta, G. Baroin, 
An introduction on formal and computational models in popular music analysis and generation,
in {\it Aesthetics and neuroscience, scientific and artistic perspectives},
\'editeurs Z. Kapoula et M. Vernet, pages 257-269, Springer, 2016.  
 
\smallskip \hangindent=7.3mm \hangafter=1 \noindent 
Aristox\`ene de Tarente,  {\it \'El\'ements  harmoniques},  -330,
traduction C.-\'E. Ruelle, Pottier de Lalaine,   
Paris, 1870.

\smallskip \hangindent=7.3mm \hangafter=1 \noindent 
V. Arnold, A. Avez,  {\it Probl\`emes ergodiques de la m\'ecanique classique}, Gauthier-Villars, Paris, 1967.

\smallskip \hangindent=7.3mm \hangafter=1 \noindent 
R. Bosanquet, {\it An elementary treatise on musical intervals,
  with an account of an enharmonic harmonium exhibited in the loan collection of scientific instruments},
Macmillan and Co, Londres, 1876.

\smallskip \hangindent=7.3mm \hangafter=1 \noindent
A. Boudet, Gammes et modes musicaux, academia.edu/44898805
et spirit-science.fr, 2006.

\smallskip \hangindent=7.3mm \hangafter=1 \noindent
J. Broekaert, Le m\'esotonique temp\'er\'e de Bach,
relatif \`a~\guillemotleft~Das wohltemperirte Clavier~\guillemotright,
academia.edu/44820384, 2021.

\smallskip \hangindent=7.3mm \hangafter=1 \noindent 
M. Brou\'e, 
Les tonalit\'es musicales vues par un math\'ematicien, 
in {\it Le code}, D. Rousseau et M.  Morvan \'editeurs, pages 41-82, 
 Odile Jacob, Paris, 2002. 

\smallskip \hangindent=7.3mm \hangafter=1 \noindent   
A. Calvet, {\it Le clavier bien obtemp\'er\'e}, piano e forte \'editions, Montpellier, 2020.

\smallskip \hangindent=7.3mm \hangafter=1 \noindent   
N. Carey, D. Clampitt, 
Aspects of well-formed scales, 
{\it Music Theory Spectrum},  volume 11, pages 187-206, 1989.

\smallskip \hangindent=7.3mm \hangafter=1 \noindent
J. Chailley, {\it L'Imbroglio des modes}, Leduc, Paris, 1960.

\smallskip \hangindent=7.3mm \hangafter=1 \noindent  
M. Chemiller, 
Grammaires, automates et musique, 
in {\it Informatique musicale : du signal au signe musical}, \'editeurs 
F. Pachet et J.-P. Briot, Lavoisier, Paris, 2004. 

\smallskip \hangindent=7.3mm \hangafter=1 \noindent  
M. Chemiller, {\it Les math\'ematiques naturelles}, Odile Jacob, Paris, 2007. 

\smallskip \hangindent=7.3mm \hangafter=1 \noindent   
D. K. L. Chua,
{\it The galitzin quartets of Beethoven: opp. 127, 132, 130}, 	
Princeton university Press, 1995. 

\smallskip \hangindent=7.3mm \hangafter=1 \noindent
C.  Clav\`ere, {\it Gamme G4-14}, 
https://www.youtube.com/watch?v=y3POXTg1Afk, 2021.

\smallskip \hangindent=7.3mm \hangafter=1 \noindent
A. Dani\'elou, {\it Trait\'e de musicologie compar\'ee}, Hermann, Paris, 1959.


\smallskip \hangindent=7.3mm \hangafter=1 \noindent 
M. Demazure,  {\it Cours d'alg\`ebre ; primalit\'e, divisibilit\'e, codes},
Nouvelle biblioth\`eque math\'ematique, Cassini, Paris, 1997.

\smallskip \hangindent=7.3mm \hangafter=1 \noindent 
M. Duchesne-Guillemin, 
Sur la restitution de la musique hourrite,
{\it Revue de Musicologie}, volume 66, pages 5-26, 1980. 


\smallskip \hangindent=7.3mm \hangafter=1 \noindent
A. J. Ellis,  On the musical scales of various nations,
{\it Journal of the Society of Arts},  volume~33, pages  485-527, 1885. 

\smallskip \hangindent=7.3mm \hangafter=1 \noindent   
A. Hachlef et M.-E. Hachlef, {\it Anthologie de la musique arabe, 1906-1960},
Publisud \'editions, Paris, 1993.

\smallskip \hangindent=7.3mm \hangafter=1 \noindent  
T. H\'elie, C. Picasso, A. Calvet,
The Snail : un nouveau proc\'ed\'e d'analyse et de visualisation du son,
Congr\`es Europiano France, Tours, 2016, hal.archives-ouvertes.fr/hal-01467014, 2016.

\smallskip \hangindent=7.3mm \hangafter=1 \noindent 
Y. Hellegouarch,  Gammes naturelles, {\it La gazette des math\'ematiciens},
n$^{\rm o}$81, pages~25-39, juillet 1999 et  n$^{\rm o}$82, pages~13-25, octobre 1999.

\smallskip \hangindent=7.3mm \hangafter=1 \noindent  
H. von Helmholtz,
{\it Die Lehre von den Tonempfindungen als physiologische Grundlage f\"ur die Theorie der Musik}, 
Friedrich Vieweg \& Sohn, Braunschweig,  1863. Traduction   Fran\c caise de G. Gu\'eroult et A. Wolff,
{\it Th\'eorie physiologique de la musique, fond\'ee sur l'\'etude  des sensations auditives},
Victor Masson et Fils, Paris, 1868. 

\smallskip \hangindent=7.3mm \hangafter=1 \noindent  
W. Holder,
{\it  A treatise of the natural grounds, and principles of harmony}, 
 J. Heptinstall, London, 1694, 1702 ; W. Pearson, London, 1731. 

\smallskip \hangindent=7.3mm \hangafter=1 \noindent 
P. von Jank\'o, 
Ueber mehr als zw\"olfstufige gleichswebende Temperaturen, pages 6-12,
{\it Beitr\"age zur Akustik und Musikwissenschaft}, volume 3, pages 6-12, \'edit\'e par C. Stumpf, Leipzig, 1901. 

\smallskip \hangindent=7.3mm \hangafter=1 \noindent  
J.-R. Jannot, La lyre et la cithare : les instruments \`a cordes de la musique \'etrusque,
{\it L'antiquit\'e classique}, volume  48, pages 469-507, 1979. 

\smallskip \hangindent=7.3mm \hangafter=1 \noindent   
F. Jedrzejewski,
Tresses n\'eoriemanniennes ; quelques applications musicales de la th\'eorie des n\oe uds,
{\it L'Ouvert}, volume 114, pages 37-49, Universit\'e de Strasbourg, 2007. 

\smallskip \hangindent=7.3mm \hangafter=1 \noindent   
F. Jedrzejewski,
{\it H\'et\'erotopies musicales ; mod\`eles math\'ematiques de la musique},
Hermann, Paris, 2019.

\smallskip \hangindent=7.3mm \hangafter=1 \noindent 
E. Lendvai, 
{\it B\'el\`a Bart\'ok: an analysis of his music},  with an introduction by A. Bush,
Kahn \& Averill, Londres, 1971.

\smallskip \hangindent=7.3mm \hangafter=1 \noindent
M. Levine,  {\it The jazz theory book},  Sher  Music Company, Petaluma, Californie, 1995.

\smallskip \hangindent=7.3mm \hangafter=1 \noindent 
F. Marpurg, {\it Versuch \"uber die musikalische Temperatur
  nebst einem Anhang \"uber den Rameau und Kirnbergerschen Grundba\ss}, Johann Friedrich Korn, Breslau, 1776.

\smallskip \hangindent=7.3mm \hangafter=1 \noindent  
M. Mersenne, {\it L'Harmonie universelle, contenant la th\'eorie et la pratique de la musique},
S\'ebastien Cramoisy, ru\"e S. Jacques, aux Cigognes, Paris, 1636.

\smallskip \hangindent=7.3mm \hangafter=1 \noindent 
P. Moutal, {\it Hindustani Raga Sangeet ; m\'ecanismes de base de la musique classique du nord de l'Inde},
Centre d'\'Etudes de Musique Orientale, Paris,  1987.

\smallskip \hangindent=7.3mm \hangafter=1 \noindent 
V. Persichetti, {\it Twentieth century harmony; creative aspects and practice},
W. W. Norton \& Company, New York, 1961.

\smallskip \hangindent=7.3mm \hangafter=1 \noindent
F. Picard, {\it La musique chinoise}, Minerve, Paris, 1991.

\smallskip \hangindent=7.3mm \hangafter=1 \noindent   
F. Picard,
Modalit\'e et pentatonisme,
{\it  Analyse Musicale}, Soci\'et\'e fran\c caise d'analyse musicale,
pages 37-46, 2001. 

\smallskip \hangindent=7.3mm \hangafter=1 \noindent
J.-P. Rameau, {\it Trait\'e de l'harmonie r\'eduite \`a ses principes naturels}, Jean-Baptiste Christophe
Ballard, rue Saint-Jean-de-Beauvais, au Mont-Parnasse, Paris, 1722.

\smallskip \hangindent=7.3mm \hangafter=1 \noindent  
H. Riemann, {\it Riemanns Musiklexikon. Theorie und Geschichte der Musik,
  die Tonk\"unstler alter und neuer Zeit mit Angabe ihrer Werke, nebst einer vollst\"andigen Instrumentenkunde},
Verlag des Bibliographischen Instituts,  Leipzig, 1882.
Traduction Fran\c caise G. Humbert, 
{\it Dictionnaire de musique de  Hugo Riemann}, Perrin et Compagnie, Paris, 1900. 

\smallskip \hangindent=7.3mm \hangafter=1 \noindent
A. Sch\"onberg, {\it Harmonielehre}, Universal Edition, Vienne, 1911.

\smallskip \hangindent=7.3mm \hangafter=1 \noindent  
L. Schwartz, {\it M\'ethodes math\'ematiques
pour les sciences physiques}, Hermann, Paris, 1965.

\smallskip \hangindent=7.3mm \hangafter=1 \noindent
J. Siron, {\it La partition int\'erieure ; jazz, musiques improvis\'ees}, Outre Mesure, Paris, 1992.

\smallskip \hangindent=7.3mm \hangafter=1 \noindent
S. Stevin, {\it Van de Spiegheling der Singconst},
manuscrit de 1605, David Bierens de Haan, Amsterdam, 1884.

\smallskip \hangindent=7.3mm \hangafter=1 \noindent  
D. Tymoczko, Scale Networks and Debussy,
{\it Journal of Music Theory}, volume~48, pages 219-294, 2004.

\smallskip \hangindent=7.3mm \hangafter=1 \noindent
A. Werckmeister, {\it Musichalische Temperatur},  Quedlinburg,  Allemagne, 1691.

\smallskip \hangindent=7.3mm \hangafter=1 \noindent  
K.  Womack, T. F. Davis (Eds),
{\it Reading the Beatles: cultural studies, literary criticism, and the fab four},
State University of New York Press, 2006.

\smallskip \hangindent=7.3mm \hangafter=1 \noindent   
S. Woo,  {\it The ceremonial music of Zhu Zaiyu}, Doctorat de Philosophie
de l'Universit\'e Rutgers, New Jersey, 2017.

\smallskip \hangindent=7.3mm \hangafter=1 \noindent
I. Wyschnegradsky, {\it Manuel d'harmonie \`a quart de ton},
 La Sir\`ene Musicale, Paris, 1932.

\smallskip \hangindent=7.3mm \hangafter=1 \noindent
I. Xenakis, {\it Musiques formelles}, La Revue Musicale, Paris, 1963.

\smallskip \hangindent=7.3mm \hangafter=1 \noindent  
G. Zarlino da Chioggia,  {\it  Le istitutioni harmoniche
del reverendo Gioseffo Zarlino da Chioggia Nelle quali~; oltra le materie appartenenti alla mvsica~;
si trouano diciarati molti luoghi di peti, historici, \& di filosofi}, Venise, 1558.

\end{document}